# Nonparametric Testing under Random Projection

Meimei Liu*, Zuofeng Shang†, Guang Cheng‡

**Abstract**

A common challenge in nonparametric inference is its high computational complexity when data volume is large. In this paper, we develop computationally efficient nonparametric testing by employing a random projection strategy. In the specific kernel ridge regression setup, a simple distance-based test statistic is proposed. Notably, we derive the minimum number of random projections that is sufficient for achieving testing optimality in terms of the minimax rate. An adaptive testing procedure is further established without prior knowledge of regularity. One technical contribution is to establish upper bounds for a range of tail sums of empirical kernel eigenvalues. Simulations and real data analysis are conducted to support our theory.

**Key Words:** Computational limit, kernel ridge regression, minimax optimality, nonparametric testing, random projection.

## 1 Introduction

A number of computationally efficient statistical methods have been proposed for analyzing massive data sets. Examples include divide-and-conquer approaches [46, 21, 10, 35]; low-rank approximations: random projection methods [29, 26, 45, 19], subsampling methods [22, 27, 1], Nyström approximations [15, 31]; and online learning methods [5, 32, 13].

An interesting question arising from these new methods is the minimum computational cost required for obtaining statistically satisfactory solutions. This might be viewed as a type of "computational limit" from a statistical perspective. Such an issue has been addressed in certain situations. For divide-and-conquer approaches, [35] derived a *sharp* upper bound for the number of distributed computing units in the smoothing spline setup, while [41] estimated the quantile regression process under an additional *sharp* lower bound on the number of quantile levels. For random projection methods, the literature nonetheless only focused on parametric cases such as compressed sensing. For example, [9] showed that the minimum number of random projections is $s \log n$ for signal recovery, where $n$ is the number of measurements and $s$ is the number of nonzero components in the true signal. To our knowledge, the computational limit for random projection methods remains unknown in nonparametric models.

*Ph.D student, Department of Statistics, Purdue University, West Lafayette, IN 47906. E-mail: liu1197@purdue.edu.

†Assistant Professor, Department of Mathematical Sciences, Indiana University-Purdue University at Indianapolis. Email: shangzf@iu.edu. Research sponsored by NSF DMS-1764280 and a startup grant from IUPUI.

‡Professor, Department of Statistics, Purdue University, West Lafayette, IN 47906. E-mail: chengg@purdue.edu. Research Sponsored by NSF CAREER Award DMS-1151692, DMS-1712907, and Office of Naval Research (ONR N00014-15-1-2331).

There are two purposes in this paper: (i) develop an optimal nonparametric testing procedure based on random projection; (ii) explore its computational limit in the kernel ridge regression setup. We remark that classical nonparametric testing methods, e.g., the locally most powerful test, the generalized/penalized likelihood ratio test and the distance-based test [11, 25, 14, 34, 2], may not be directly applied to big data due to their high computational costs.

Specifically, we consider the following nonparametric model

$$y_i = f(x_i) + \epsilon_i, \quad i = 1, \cdots, n, \tag{1.1}$$

where $x_i \in \mathcal{X} \subseteq \mathbb{R}^d$ for a fixed $d \geq 1$ are random design points, and $\epsilon_i$ are random noise with mean zero. The regression function $f$ belongs to a reproducing kernel Hilbert space (RKHS) $\mathcal{H}$. The hypothesis of interest is

$$H_0 : f = f_0 \ \text{ v.s. } H_1 : f \in \mathcal{H} \setminus \{f_0\}, \tag{1.2}$$

where $f_0$ is a hypothesized function. We construct a distance-based test statistic $T_{n,\lambda} = \|\widehat{f}_R - f_0\|_n^2$ for testing (1.2), where $\widehat{f}_R$ is a random projection version of the kernel ridge regression (KRR) estimator $\widehat{f}_n$ ([36]) defined as

$$\widehat{f}_n := \underset{f \in \mathcal{H}}{\operatorname{argmin}} \left\{ \frac{1}{n} \sum_{i=1}^{n} (y_i - f(x_i))^2 + \lambda \|f\|_{\mathcal{H}}^2 \right\}, \tag{1.3}$$

where $\|f\|_{\mathcal{H}}^2 = \langle f, f \rangle_{\mathcal{H}}$ with $\langle \cdot, \cdot \rangle_{\mathcal{H}}$ the inner product of $\mathcal{H}$, $\lambda > 0$ is a smoothing parameter. The computational cost and storage occupation of $\widehat{f}_n$ are of orders $\mathcal{O}(n^3)$ and $\mathcal{O}(n^2)$, respectively. However, computing $\widehat{f}_R$ reduces these costs to $\mathcal{O}(s^3)$ and $\mathcal{O}(s^2)$ under $s (\ll n)$ random projections; see Section 2. After $\widehat{f}_R$ is obtained, $T_{n,\lambda}$ can be computed in a parallel fashion. Hence, $s$ can be viewed as a simple proxy for computing and storage costs.

In this paper, we reveal a phase transition phenomenon in terms of $s$. Specifically, a sharp lower bound for $s$ is established: when $s$ is above this bound, $T_{n,\lambda}$ is minimax optimal; otherwise, minimax optimality becomes impossible even when the best possible $\lambda$ is chosen. We next illustrate more subtle details using Figure 1, where the strength of the weakest detectable signals (SWDS) is characterized given any $s$ and $\lambda$. In general, we require $s \geq s_\lambda$ for any $\lambda$, where $s_\lambda$ is determined by kernel eigenvalues and $\lambda$. An important observation is that the smallest SWDS can be achieved at $\lambda = \lambda^*$ and $s \geq s_{\lambda^*} := s^*$ (note that when $s \ll s^*$, our testing procedure under a proper $\lambda$ is still powerful as long as SWDS becomes sufficiently large). Both $\lambda^*$ and $s^*$ have precise orders in specific situations. For example, in an $m$-order polynomial decay kernel, the smallest SWDS achieves the minimax optimal rate $n^{-\frac{2m}{4m+1}}$ ([20]) when $\lambda^* = n^{-\frac{4m}{4m+1}}$ and $s^* = n^{\frac{2}{4m+1}}$. As a by-product, we also derive a sharp lower bound for $s$ for obtaining the minimax optimal estimation. Our results hold for a general class of random projection matrix, such as the sub-Gaussian matrix or certain data-dependent matrix.

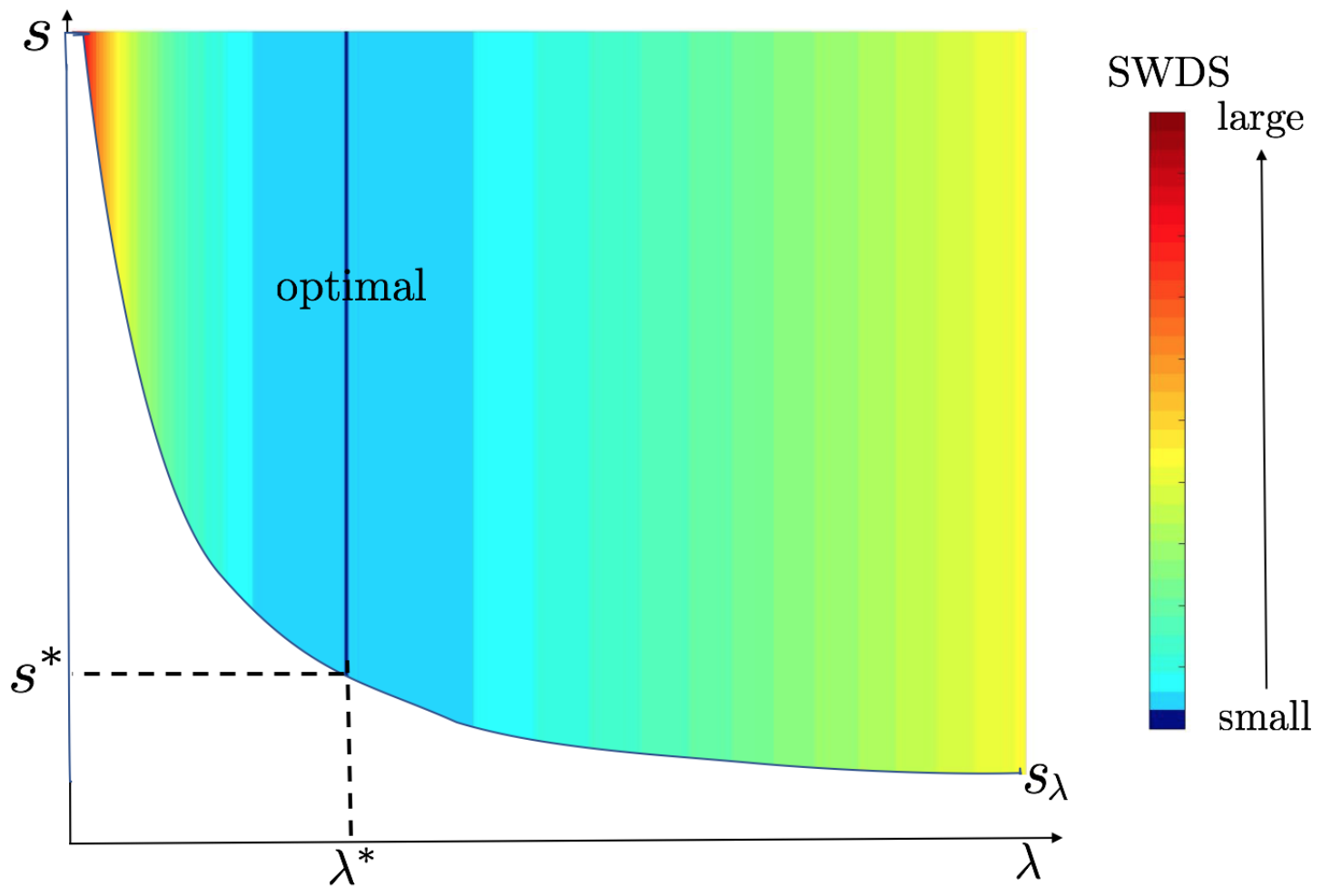

Figure 1: Phase transition in $(\lambda, s)$ for signal detection. The horizontal axis is the smoothing parameter $\lambda$, and the vertical axis is the projection dimension $s$. The shade indicates the values of SWDS: dark red corresponds to greater values of SWDS than light blue. The vertical line labeled by "optimal" indicates the choices of $\lambda$ that achieve the smallest SWDS.

It is worth mentioning that the construction of $T_{n,\lambda}$ crucially relies on the regularity of $\mathcal{H}$, which is often unavailable in practice. Hence, we propose an adaptive test statistic based on the maximum of a sequence of (standardized) non-adaptive test statistics corresponding to various regularities. Based on a recent Gaussian approximation result in [24], we prove that the null limit distribution is an extreme value distribution.

The proofs of main results rely on the behavior of the tail sum of empirical kernel eigenvalues. One technical contribution of this work is to derive upper bounds for a range of tail sums such that nonparametric estimation and testing can now be analyzed in a unified framework; see Section 3. This is obtained by flexibly adjusting the size of the function class associated with the Rademacher average in the local Rademacher complexity theory ([4]).

In simulation studies, we find that the size and power of the proposed non-adaptive and adaptive test statistics are both satisfactory. In particular, the power cannot be further improved as the number of random projections grows beyond some threshold, as predicted by our theory. For an illustration purpose, we also demonstrate that when $n = 2^{12}$, conducting testing based on $\widehat{f}_R$ only takes 3.2 seconds in comparison with 42 seconds based on $\widehat{f}_n$. In practice, the smoothing parameter $\lambda$ can be directly selected via generalized cross validation. We would like to point out that this is an advantage of the random projection method over the divide-and-conquer method [46], where the selection of the smoothing parameter is nontrivial; see [44].

The rest of this paper is organized as follows. Section 2 introduces kernel ridge regression together with its approximation based on random projection, while Section 3 develops an upper bound for the tail sum of empirical eigenvalues. Our main results are presented in Section 4: Section 4.1 introduces one primary assumption on random projection; Sections 4.2 and 4.3 study testing consistency and power behaviors in terms of the projection dimension $s$ and the smoothing parameter $\lambda$, with specific situations considered in Section 4.4; Section 4.5 proves the lower bound on $s$ given in Section 4.4 to be sharp. An adaptive testing procedure is developed in Section 5. Section 6 includes numerical studies based on simulated and real data sets. All technical details are deferred to either the Appendix or the online supplementary.

**Notation:** Denote $\delta_{jk}$ the Kronecker delta: $\delta_{jk} = 1$ if $j = k$ and $\delta_{jk} = 0$ if $j \neq k$. For positive sequences $a_n$ and $b_n$, put $a_n \lesssim b_n$ if there exists a constant $c > 0$ such that $a_n \leq c b_n$ for all $n \in \mathbb{N}$; $a_n \gtrsim b_n$ if there exists a constant $c > 0$ such that $a_n \leq c b_n$. Put $a_n \asymp b_n$ if $a_n \lesssim b_n$ and $a_n \gtrsim b_n$. Frequently, we use $a_n \lesssim b_n$ and $a_n = \mathcal{O}(b_n)$ interchangeably. $Pf^2 \equiv \mathrm{E}\, f(X)^2$, $\|f\|_n^2 \equiv P_n f^2 \equiv \frac{1}{n}\sum_{i=1}^n f(X_i)^2$. For a matrix $A \in \mathbb{R}^{m \times n}$, its operator norm is defined as $\|A\|_{\mathrm{op}} = \max_{x \in \mathbb{R}^n \setminus \{0\}} \frac{\|Ax\|_2}{\|x\|_2}$. A random variable $X$ is said to be sub-Gaussian if there exists a constant $\sigma^2 > 0$ such that for any $t \geq 0$, $\mathrm{P}[|X| \geq t] \leq 2\exp(-t^2/(2\sigma^2))$. The sub-Gaussian norm of $X$ is defined as $\|X\|_{\psi_2} = \inf\{t > 0 : \mathrm{E}\exp(X^2/t^2) \leq 1\}$. We will use $c, c_1, c_2, C$ to denote generic absolute constants, whose values may vary from line to line.

## 2 Kernel Ridge Regression via Random Projection

In this section, we review kernel ridge regression and its variant based on random projection. Suppose that we have $n$ i.i.d. observations $\{(x_i, y_i)\}_{i=1}^n$ from (1.1). Throughout assume that $f \in \mathcal{H}$, where $\mathcal{H} \subset L^2(P_X)$ is a reproducing kernel Hilbert space (RKHS) associated with an inner product $\langle \cdot, \cdot \rangle_{\mathcal{H}}$ and a reproducing kernel function $K(\cdot, \cdot) : \mathcal{X} \times \mathcal{X} \to \mathbb{R}$. By Mercer's theorem, $K$ has the following spectral expansion:

$$K(x, x') = \sum_{i=1}^{\infty} \mu_i \phi_i(x) \phi_i(x'), \quad x, x' \in \mathcal{X}, \tag{2.1}$$

where $\mu_1 \geq \mu_2 \geq \cdots \geq 0$ is a sequence of ordered eigenvalues and the eigenfunctions $\{\phi_i\}_{i=1}^{\infty}$ form a basis in $L^2(P_X)$. Moreover, for any $i, j \in \mathbb{N}$,

$$\langle \phi_i, \phi_j \rangle_{L^2(P_X)} = \delta_{ij} \qquad \text{and} \qquad \langle \phi_i, \phi_j \rangle_{\mathcal{H}} = \delta_{ij}/\mu_i.$$

Throughout this paper, assume that $\phi_i$'s are uniformly bounded, a common condition in literature, e.g., [17], and $\mu_i$'s satisfy certain tail sum property.

**Assumption A1.** $c_K := \sup_{i \geq 1} \|\phi_i\|_{\sup} < \infty$ *and* $\sup_{k \geq 1} \frac{\sum_{i=k+1}^{\infty} \mu_i}{k \mu_k} < \infty$.

Assumption A1 is satisfied in two types of commonly used kernels, categorized by the eigenvalue decay rates. The first is $\mu_i \asymp i^{-2m}$ for a constant $m > 0$, called as polynomial decay kernel (PDK) of order $m$. The second is $\mu_i \asymp \exp(-\gamma i^p)$ for constants $\gamma, p > 0$, called as exponential decay kernel (EDK) of order $p$. Verification of Assumption A1 is deferred to Section S.3 of the supplement document. Examples of PDK include kernels of Sobolev space and periodic Sobolev space (see [42]). Examples of EDK include Gaussian kernel $K(x_1, x_2) = \exp(-(x_1 - x_2)^2/2)$ (see [33]).

Recall the KRR estimator $\widehat{f}_n$ from (1.3). By representer theorem, it has an expression $\widehat{f}_n(\cdot) = \sum_{i=1}^{n} \widehat{\omega}_i K(\cdot, x_i)$, where $\widehat{\omega} = (\widehat{\omega}_1, \ldots, \widehat{\omega}_n)^\top$ is a real vector determined by

$$\begin{aligned} \widehat{\omega} &= \operatorname*{argmin}_{\omega \in \mathbb{R}^n} \left\{ \omega^\top \boldsymbol{K}^2 \omega - \frac{2}{n} \omega^\top \boldsymbol{K} \boldsymbol{y} + \lambda \omega^\top \boldsymbol{K} \omega \right\} \\ &= \frac{1}{n} (\boldsymbol{K} + \lambda I)^{-1} \boldsymbol{y}, \end{aligned} \tag{2.2}$$

$\boldsymbol{y} = (y_1, \cdots, y_n)^\top$, $\boldsymbol{K} = [n^{-1} K(x_i, x_j)]_{1 \leq i,j \leq n}$, and $I \in \mathbb{R}^{n \times n}$ is identity. This standard procedure requires storing $(\boldsymbol{K}^2, \boldsymbol{K}, \boldsymbol{K}\boldsymbol{y})$ and inverting $\boldsymbol{K} + \lambda I$, which requires $\mathcal{O}(n^2)$ memory usage and $\mathcal{O}(n^3)$ floating operations.

The above computational and storage constraints become severe for a large sample size, and thus motivate the random projection approach proposed by [45]. Specifically, $\omega$ in (2.2) is substituted with $S^\top \beta$, where $\beta \in \mathbb{R}^s$ and $S$ is an $s \times n$ real-valued random matrix; see Section 4.1. Then, $\beta$ is solved as:

$$\begin{aligned} \widehat{\beta} &= \operatorname*{argmin}_{\beta \in \mathbb{R}^s} \left\{ \beta^\top (S\boldsymbol{K})(\boldsymbol{K}S^\top)\beta - \frac{2}{n} \beta^\top S \boldsymbol{K} \boldsymbol{y} + \lambda \beta^\top S \boldsymbol{K} S^\top \beta \right\}, \\ &= \frac{1}{n} (S \boldsymbol{K}^2 S^\top + \lambda S \boldsymbol{K} S)^{-1} S \boldsymbol{K} \boldsymbol{y}. \end{aligned} \tag{2.3}$$

Hence, the resulting estimator of $f$ becomes

$$\widehat{f}_R(\cdot) = \sum_{i=1}^{n} (S^\top \widehat{\beta})_i K(\cdot, x_i), \tag{2.4}$$

which requires storing $(S\boldsymbol{K}^2S^\top, S\boldsymbol{K}S^\top, S\boldsymbol{K}\boldsymbol{y})$ and inverting an $s \times s$ matrix. Hence, the memory usage and floating operations are reduced to $\mathcal{O}(s^2)$ and $\mathcal{O}(s^3)$, respectively, when $s = o(n)$. On the other hand, $s$ cannot be too small in order to maintain sufficient data information for achieving statistical optimality. Critical lower bounds for $s$ will be derived in Section 4.5.

# 3 Tail Sum of Empirical Eigenvalues

An accurate upper bound for the tail sum of empirical eigenvalues is needed for studying nonparametric testing and estimation. However, this bound was often *assumed* to hold in the kernel learning literature, e.g., [7, 43]. And, the application of concentration inequalities of individual eigenvalues ([37, 8]) only provides a very loose bound due to accumulative errors. Recently, the local Rademacher complexity (LRC) theory ([4]) was employed by [45] to derive a more accurate upper bound that is useful in studying nonparametric estimation. However, this upper bound no longer works for testing problems, due to the improper size of the function class defining Rademacher average.

In this section, we establish upper bounds, i.e., Lemma 3.1, for a range of tail sums of empirical eigenvalues that can be applied to both nonparametric estimation and testing. This result may be of independent interest. Consider the singular value decomposition $\boldsymbol{K} = UDU^\top$, where $UU^\top = I_n$ and $D = \mathrm{diag}(\widehat{\mu}_1, \widehat{\mu}_2, \ldots, \widehat{\mu}_n)$ with $\widehat{\mu}_1 \geq \widehat{\mu}_2 \geq \cdots \widehat{\mu}_n \geq 0$. For any $\lambda > 0$, define $\widehat{s}_\lambda$ (or $s_\lambda$) to be the number of $\widehat{\mu}_i$'s (or $\mu_i$'s) greater than $\lambda$, i.e.,

$$\widehat{s}_\lambda = \mathrm{argmin}\{i : \widehat{\mu}_i \leq \lambda\} - 1, \quad s_\lambda = \mathrm{argmin}\{i : \mu_i \leq \lambda\} - 1. \tag{3.1}$$

For a range of $\lambda$, Lemma 3.1 below provides an upper bound for the tail sum of $\widehat{\mu}_i$ in terms of population quantities $s_\lambda$ and $\mu_{s_\lambda}$, with known orders.

**Lemma 3.1.** *If $1/n < \lambda \to 0$, then with probability at least $1 - 4e^{-s_\lambda}$, $\sum_{i=\widehat{s}_{\lambda+1}}^{n} \widehat{\mu}_i \leq C s_\lambda \mu_{s_\lambda}$, where $C > 0$ is an absolute constant.*

Clearly, Lemma 3.1 is a sample analog to the tail sum assumption for $\mu_i$ in Assumption A1. The proof of Lemma 3.1 is based on an adaptation of the classical LRC theory as explained below.

In Section 4, it will be shown that $\lambda$ and $s_\lambda/n$ correspond to (squared-)bias and variance of $\widehat{f}_R$, respectively. We then define the variance-to-bias ratio as

$$\kappa_\lambda = \frac{s_\lambda}{n\lambda}, \tag{3.2}$$

for any $\lambda > 0$. Consider a bundle of function classes indexed by $\kappa_\lambda$:

$$\mathcal{F}_\lambda = \{f \in \mathcal{H} : f \text{ maps } \mathcal{X} \text{ to } [-1,1], \|f\|_{\mathcal{H}}^2 \le \kappa_\lambda\}, \lambda > 0.$$

To characterize the complexity of $\mathcal{F}_\lambda$, we introduce a generalized version of local Rademacher complexity function:

$$\Psi_\lambda(r) = \mathrm{E}\Big\{\sup_{\substack{f\in\mathcal{F}_\lambda\\ Pf^2\le r}} \frac{1}{n}\sum_{i=1}^n \sigma_i f(x_i)\Big\}, r \ge 0,$$

where $\sigma_1, \dots, \sigma_n$ are independent Rademacher random variables, i.e., $\mathrm{P}(\sigma_i = 1) = \mathrm{P}(\sigma_i = -1) = 1/2$. Let $\widehat{\Psi}_\lambda(\cdot)$ be an empirical version of $\Psi_\lambda(\cdot)$ defined as

$$\widehat{\Psi}_\lambda(r) = \mathrm{E}\Big\{\sup_{\substack{f\in\mathcal{F}_\lambda\\ P_nf^2\le r}} \frac{1}{n}\sum_{i=1}^n \sigma_i f(x_i)\Big| x_1, \cdots, x_n\Big\}, r \ge 0.$$

When $\kappa_\lambda \asymp 1$, $\Psi_\lambda(\cdot)$ and $\widehat{\Psi}_\lambda(\cdot)$ become the original LRC functions introduced in [4]. Note that $\kappa_\lambda \asymp 1$ actually corresponds to the optimal bias vs. variance trade-off required for estimation. Rather, a different type of trade-off is needed for optimal testing as revealed by [20, 34], which corresponds to a different choice of $\kappa_\lambda$ in $\mathcal{F}_\lambda$ as demonstrated later in Section 4.

Lemma 3.2 below says that both $\Psi_\lambda$ and $\widehat{\Psi}_\lambda$ possess unique (positive) fixed points. This fixed point property is crucial in proving Lemma 3.1. Interestingly, we find that the fixed points turn out to be proportional to the estimation variance asymptotically.

**Lemma 3.2.** *There exist uniquely positive $r_\lambda$ and $\widehat{r}_\lambda$ such that $\Psi_\lambda(r_\lambda) = r_\lambda$ and $\widehat{\Psi}_\lambda(\widehat{r}_\lambda) = \widehat{r}_\lambda$. Furthermore, if $\lambda > 1/n$, then $r_\lambda \asymp s_\lambda/n$, and there exists an absolute constant $c > 0$ such that, with probability at least $1 - e^{-cs_\lambda}$, $\widehat{r}_\lambda \asymp s_\lambda/n$.*

We are now ready to sketch the proof of Lemma 3.1. Detailed proofs are deferred to Appendix S.5. First, note that

$$\sum_{i=\widehat{s}_\lambda+1}^n \widehat{\mu}_i = \sum_{i=\widehat{s}_\lambda+1}^n \min\{\lambda, \widehat{\mu}_i\} \le \sum_{i=1}^n \min\{\lambda, \widehat{\mu}_i\}.$$

By Lemma 3.2, we have $\widehat{r}_\lambda/\kappa_\lambda \asymp \lambda$ with high probability. Then,

$$\begin{aligned}\sum_{i=1}^n \min\{\lambda, \widehat{\mu}_i\} &\asymp \sum_{i=1}^n \min\Big\{\frac{\widehat{r}_\lambda}{\kappa_\lambda}, \widehat{\mu}_i\Big\}\\ &\asymp \frac{n}{\kappa_\lambda}\widehat{\Psi}_\lambda(\widehat{r}_\lambda)^2 = \frac{n\widehat{r}_\lambda^2}{\kappa_\lambda} \asymp \lambda s_\lambda \le s_\lambda \mu_{s_\lambda},\end{aligned}$$

where the second step is by Lemma S.1 that

$$\widehat{\Psi}_\lambda(\widehat{r}_\lambda) \asymp \sqrt{\frac{\kappa_\lambda}{n}\sum_{i=1}^n \min\{\frac{\widehat{r}_\lambda}{\kappa_\lambda}, \widehat{\mu}_i\}},$$

the third step follows from the fixed point property stated in Lemma 3.2, and the last step follows from the definition of $\mu_{s_\lambda}$ given in (3.1).

# 4 Main Results

Consider the nonparametric testing problem (1.2). For convenience, assume $f_0 = 0$, i.e., we will test

$$H_0 : f = 0 \quad \text{vs.} \quad H_1 : f \in \mathcal{H} \setminus \{0\}. \tag{4.1}$$

In general, testing $f = f_0$ (for an arbitrary known $f_0$) is equivalent to testing $f_* \equiv f - f_0 = 0$. So, (4.1) has no loss of generality. Based on $\widehat{f}_R$, we propose the following distance-based test statistic:

$$T_{n,\lambda} = \|\widehat{f}_R\|_n^2. \tag{4.2}$$

In the subsequent sections, we will derive the null limit distribution of $T_{n,\lambda}$ (Theorems 4.2 and 4.5), and further provide a sufficient and necessary condition in terms of $s$ such that $T_{n,\lambda}$ is minimax optimal (Section 4.5). As a byproduct, we derive a critical bound in terms of $s$ such that $\widehat{f}_R$ is minimax optimal. Proof of such results rely on an exact analysis on the kernel and projection matrices which requires an accurate estimate of the tail sum of the empirical eigenvalues by Lemma 3.1. Our results hold for a general choice of projection matrix which will be discussed in Section 4.1.

## 4.1 Choice of Projection Matrix

Recall the singular value decomposition $\boldsymbol{K} = UDU^\top$. Put $U = (U_1, U_2)$ with $U_1$ consisting of the first $\widehat{s}_\lambda$ columns of $U$ and $U_2$ consisting of the rest $n - \widehat{s}_\lambda$ columns; $D = \text{diag}(D_1, D_2)$, with $D_1 = \text{diag}(\widehat{\mu}_1, \ldots, \widehat{\mu}_{\widehat{s}_\lambda})$, $D_2 = \text{diag}(\widehat{\mu}_{\widehat{s}_\lambda+1}, \ldots, \widehat{\mu}_n)$.

The following definition of "$\boldsymbol{K}$-satisfiability" describes a class of matrices that preserve the principal components of the kernel matrix.

**Definition 4.1.** ($\boldsymbol{K}$-satisfiability) A matrix $S \in \mathbb{R}^{s\times n}$ is said to be $\boldsymbol{K}$-satisfiable if there exists a constant $c > 0$ such that

$$\|(SU_1)^\top SU_1 - I_{\widehat{s}_\lambda}\|_{\text{op}} \le 1/2, \quad \|SU_2D_2^{1/2}\|_{\text{op}} \le c\lambda^{1/2}.$$

By Definition 4.1, a $\boldsymbol{K}$-satisfiable $S$ will make $(SU_1)^\top SU_1$ "nearly" identity as well as downweight the tail eigenvalues. Such a matrix will be able to extract the principle information from the kernel matrix; see also [45].

Besides, we need the following definition which will make the statement of our assumptions concise.

**Definition 4.2.** An event $\mathcal{E}$ is said to be of $(a,b)$-type for $a,b \in (0,\infty]$, if $\mathrm{P}(\mathrm{P}(\mathcal{E}|x_1,\cdots,x_n) \geq 1-\exp(-a)) \geq 1-\exp(-b)$.

Definition 4.2 describes events whose probabilities have exponential type lower bounds. It is easy to see that, if $\mathcal{E}$ is of $(a,b)$-type, then $\mathrm{P}(\mathcal{E}) \geq (1-\exp(-a))(1-\exp(-b))$. In particular, $\mathcal{E}$ is of $(\infty,\infty)$-type if and only if $\mathcal{E}$ occurs almost surely.

Throughout the rest of this paper, assume the following condition on $S$.

**Assumption A2.** *(a) $s \geq ds_\lambda$ for a sufficiently large constant $d > 0$.*

*(b) There exist $c_1, c_2 \in (0,\infty]$ such that the event "$S$ is $\boldsymbol{K}$-satisfiable" is of $(c_1 s, c_2 s_\lambda)$-type.*

Assumption A2 (a) requires a sufficient amount of random projections to preserve data information. Assumption A2 (b) requires $S$ to be $\boldsymbol{K}$-satisfiable with high probability which holds in a broad range of situations such as certain data dependent matrix (Example 4.3) or matrix of sub-Gaussian entries (Example 4.4).

**Example 4.3.** Let $S = U_s^\top$, where $U_s$ is an $n \times s$ matrix consisting of the first $s$ columns of $U$. Then it trivially holds that, almost surely, $(SU_1)^\top SU_1 = I_{\widehat{s}_\lambda}$ and $\|SU_2 D_2^{1/2}\|_{\mathrm{op}} = 0$, i.e., Assumption A2 (b) holds. Such a construction of $S$ has also been utilized in [3].

**Example 4.4.** Let $S$ be an $s \times n$ random matrix of entries $S_{ij}/\sqrt{s}$, $i = 1,\ldots,s$, $j = 1,\ldots,n$, where $S_{ij}$ are independent (not necessarily identically distributed) sub-Gaussian variables. Examples of such sub-Gaussian variables include Gaussian variables, bounded variables such as Bernoulli, multinomial, uniform, variables with strictly log-concave density, or mixtures of sub-Gaussian variables. The following lemma shows that Assumption A2 (b) holds in all these situations.

**Lemma 4.1.** *Let $S_{ij} : 1 \leq i \leq s$, $1 \leq j \leq n$ be independent sub-Gaussian of mean zero and variance one, and $\lambda \in (1/n, 1)$. If $s \geq ds_\lambda$ for a sufficiently large constant $d$, then Assumption A2 (b) holds for $S = [S_{ij}/\sqrt{s}]_{1\leq i\leq s, 1\leq j\leq n}$.*

## 4.2 Testing Consistency

In this section, we derive the null limit distribution of (standardized) $T_{n,\lambda}$ as standard Gaussian, and then extend our result to the case of composite hypothesis testing.

**Theorem 4.2.** *Suppose that $\lambda \to 0$ and $s \to \infty$ as $n \to \infty$. Then under $H_0$, we have*

$$\frac{T_{n,\lambda} - \mu_{n,\lambda}}{\sigma_{n,\lambda}} \xrightarrow{d} N(0,1), \quad \text{as } n \to \infty.$$

*Here, $\mu_{n,\lambda} := \mathrm{E}_{H_0}\{T_{n,\lambda}|\boldsymbol{x}, S\} = \mathrm{tr}(\Delta^2)/n$, $\sigma_{n,\lambda}^2 := \mathrm{Var}_{H_0}\{T_{n,\lambda}|\boldsymbol{x}, S\} = 2\,\mathrm{tr}(\Delta^4)/n^2$ with $\boldsymbol{x} = (x_1,\cdots,x_n)$ and $\Delta = \boldsymbol{K}S^\top(S\boldsymbol{K}^2S^\top + \lambda S\boldsymbol{K}S^\top)^{-1}S\boldsymbol{K}$.*

Theorem 4.2 holds once $s$ diverges (no matter how slowly). Theorem 4.2 implies the following testing rule at significance level $\alpha$:

$$\phi_{n,\lambda} = I(|T_{n,\lambda} - \mu_{n,\lambda}| \geq z_{1-\alpha/2}\sigma_{n,\lambda}) \tag{4.3}$$

where $z_{1-\alpha/2}$ is the $100 \times (1-\alpha/2)$th percentile of $N(0,1)$.

An important consequence of Theorem 4.2 is the following estimation rate

$$\|\widehat{f}_R - f_0\|_n^2 = O_P(r_{n,\lambda}), \tag{4.4}$$

where $r_{n,\lambda} = \lambda + \mu_{n,\lambda}$. The proof of (4.4) is sketched as follows. Suppose that $f_0 \in \mathcal{H}$ is the "true" function in (1.1). Note that $\|\widehat{f}_R - f_0\|_n^2$ has a trivial upper bound

$$\|\widehat{f}_R - f_0\|_n^2 \leq 2\|\widehat{f}_R - \mathrm{E}_\epsilon \widehat{f}_R\|_n^2 + 2\|\mathrm{E}_\epsilon \widehat{f}_R - f_0\|_n^2, \tag{4.5}$$

where $\mathrm{E}_\epsilon$ is the expectation w.r.t. $\epsilon$. By direct examinations, it can be shown that $\|\widehat{f}_R - \mathrm{E}_\epsilon \widehat{f}_R\|_n^2 = \epsilon^\top \Delta^2 \epsilon / n$, hence, $\mathrm{E}_\epsilon \|\widehat{f}_R - \mathrm{E}_\epsilon \widehat{f}_R\|_n^2 = \mathrm{tr}(\Delta^2)/n = \mu_{n,\lambda}$. This leads to $\|\widehat{f}_R - \mathrm{E}_\epsilon \widehat{f}_R\|_n^2 = O_P(\mu_{n,\lambda})$. Meanwhile, it follows from Lemma 4.3 below that $\|\mathrm{E}_\epsilon \widehat{f}_R - f_0\|_n^2 = O_P(\lambda)$. This completes the proof of (4.4).

**Lemma 4.3.** *Suppose that $1/n < \lambda < 1$ and Assumption A2 holds with $c_1, c_2 \in (0, \infty]$. Then with probability greater than $1 - e^{-c_1 s} - e^{-c_2 s_\lambda}$,*

$$\|\mathrm{E}_\epsilon \widehat{f}_R - f_0\|_n^2 \leq C\lambda,$$

*where $C$ is a positive absolute constant.*

The above discussions are summarized in the following corollary.

**Corollary 4.4.** *Suppose that $1/n < \lambda < 1$ and Assumption A2 holds. Then with probability approaching one, it holds that*

$$\|\widehat{f}_R - f_0\|_n^2 \leq C r_{n,\lambda},$$

*where $r_{n,\lambda} = \lambda + \mu_{n,\lambda}$ and $C$ is an absolute constant.*

From Corollary 4.4, the best upper bound can be obtained through balancing $\lambda$ and $\mu_{n,\lambda}$. Denote $\lambda^\dagger$ the optimizer. This in turn provides a lower bound $s^\dagger$ for $s$ according to (3.1), i.e., $s^\dagger = s_{\lambda^\dagger}$. In Section 4.4, we will show that the upper bound under $\lambda^\dagger$ is minimax optimal, and further provide explicit orders for $s^\dagger$ in concrete settings.

In practice, it is often of interest to test certain structure of $f$, e.g., linearity,

$$H_0^{\mathrm{linear}} : f \in \mathcal{L}(\mathcal{X}) \text{ vs. } H_1^{\mathrm{linear}} : f \notin \mathcal{L}(\mathcal{X}),$$

where $\mathcal{L}(\mathcal{X})$ is the class of linear functions over $\mathcal{X} \subseteq \mathbb{R}^d$. Testing $H_0^{\text{linear}}$ can be easily converted into simple hypothesis testing problem as follows. Suppose that $f_0(x) = \beta_0 + \beta_1^\top x$ is the "true" function under $H_0^{\text{linear}}$. The corresponding MLE is $\widehat{f}_0(x) = \widehat{\beta}_0 + \widehat{\beta}_1^\top x$, where $\widehat{\beta} = (\mathbf{X}\mathbf{X}^\top)^{-1}\mathbf{X}\boldsymbol{y} \equiv (\widehat{\beta}_0, \widehat{\beta}_1)$. By defining $f^* = f - \widehat{f}_0$, it amounts to testing $f^* = 0$. Correspondingly, we define $\boldsymbol{y}^* = \boldsymbol{y} - \widehat{\boldsymbol{y}}_0$, where $\widehat{\boldsymbol{y}}_0 = (\widehat{f}_0(x_1), \ldots, \widehat{f}_0(x_n))^\top = H\boldsymbol{y}$ and $H = \mathbf{X}^\top(\mathbf{X}\mathbf{X}^\top)^{-1}\mathbf{X}$. This leads to $\widehat{f}_R^*$ and $T_{n,\lambda}^* = \|\widehat{f}_R^*\|_n^2$, whose null limit distribution is given in the following theorem.

**Theorem 4.5.** *Suppose that $\lambda \to 0$ and $s \to \infty$ as $n \to \infty$. Under $H_0^{linear}$, we have*

$$\frac{T_{n,\lambda}^* - \mu_{n,\lambda}^*}{\sigma_{n,\lambda}^*} \xrightarrow{d} N(0,1), \quad \text{as } n \to \infty,$$

*where $\mu_{n,\lambda}^* = \mathrm{E}_{H_0^{linear}}\{T_{n,\lambda}^* | \boldsymbol{x}, S\} = \mathrm{tr}((I-H)\Delta^2(I-H))/n$ and $\{\sigma_{n,\lambda}^*\}^2 = \mathrm{Var}_{H_0^{linear}}(T_{n,\lambda}^* | \boldsymbol{x}, S) = 2\,\mathrm{tr}((I-H)\Delta^4(I-H))/n^2$.*

Clearly, our testing procedure and theory can be easily generalized to polynomial testing such as $H_0^{\text{poly}}$ : $f$ is polynomial of order $q$.

### 4.3 Power Analysis

In this section, we investigate the power of $T_{n,\lambda}$ under a sequence of local alternatives. The following result shows that $T_{n,\lambda}$ can achieve high power provided that $s$ diverges fast enough *and* the local alternative is separated from the null by at least an amount of $d_{n,\lambda}$.

**Theorem 4.6.** *Suppose that $1/n < \lambda \to 0$ as $n \to \infty$, and Assumption A2 holds for $c_1, c_2 \in (0, \infty]$. Then for any $\varepsilon > 0$, there exist positive constants $C_\varepsilon$ and $N_\varepsilon$ such that, with probability greater than $1 - e^{-c_1 s} - e^{-c_2 s_\lambda}$,*

$$\inf_{n \geq N_\varepsilon} \inf_{\substack{f \in \mathcal{B} \\ \|f\|_n \geq C_\varepsilon d_{n,\lambda}}} P_f(\phi_{n,\lambda} = 1 | \boldsymbol{x}, S) \geq 1 - \varepsilon,$$

*where $d_{n,\lambda} := \sqrt{\lambda + \sigma_{n,\lambda}}$ and $\mathcal{B} = \{f \in \mathcal{H} : \|f\|_{\mathcal{H}} \leq C\}$ for a constant $C$ and $P_f(\cdot | \boldsymbol{x}, S)$ is the conditional probability measure under $f$ given $\boldsymbol{x}, S$.*

In view of Theorem 4.6, to maximize the power of $T_{n,\lambda}$, one needs to minimize $d_{n,\lambda} = \sqrt{\lambda + \sigma_{n,\lambda}}$ through balancing $\lambda$ and $\sigma_{n,\lambda}$. Denote $\lambda^*$ the optimizer. The lower bound $s^*$ for $s$ is obtained via (3.1), i.e., $s^* = s_{\lambda^*}$. The explicit forms of $\lambda^*$ and $s^*$ will be provided in Section 4.4.

### 4.4 Examples

As an application of our main results, i.e., Corollary 4.4 and Theorem 4.6, we will derive the lower bounds for $s$ to achieve optimal estimation and testing in two spacial cases.

**Example 1: PDK**

Suppose that $\mathcal{H}$ is generated by an $m$-order PDK. The following Lemma characterizes the orders of $\mu_{n,\lambda}$ and $\sigma^2_{n,\lambda}$.

**Lemma 4.7.** *Suppose that $1/n < \lambda \to 0$ as $n \to \infty$. Meanwhile, Assumption A2 holds with $c_1, c_2 \in (0, \infty]$ and $m > 3/2$. Then with probability at least $1 - e^{-c_m n^{(2m-3)/(2m-1)}} - e^{-c_1 s} - e^{-c_2 s_\lambda}$, it holds that $\mu_{n,\lambda} \asymp s_\lambda/n$ and $\sigma^2_{n,\lambda} \asymp s_\lambda/n^2$, where $c_m$ is an absolute constant depending on $m$ only.*

It follows from Lemma 4.7 and Corollary 4.4 that $\widehat{f}_R$ has the convergence rate $r_{n,\lambda} = \lambda + s_\lambda/n$. Note that $\lambda$ is the bias of $\widehat{f}_R$ by Lemma 4.3, and $s_\lambda/n$ is the variance of $\widehat{f}_R$ by (4.5) and Lemma 4.7. Hence, the optimal estimation rate $r^\dagger_{n,\lambda}$ is achieved as follows

$$r^\dagger_{n,\lambda} = \operatorname{argmin}\Big\{\lambda \,:\, \lambda > s_\lambda/n\Big\}.$$

To be concrete, by (3.1) and $\mu_i \asymp i^{-2m}$, it can be shown that $s_\lambda \asymp \lambda^{-\frac{1}{2m}}$. Hence, $\lambda^\dagger = n^{-\frac{2m}{2m+1}}$ and $s^\dagger \asymp (\lambda^\dagger)^{-\frac{1}{2m}} = n^{\frac{1}{2m+1}}$. In summary, $\widehat{f}_R$ achieves the minimax rate of estimation $n^{-\frac{2m}{2m+1}}$ when $s \gtrsim n^{\frac{1}{2m+1}}$ and $\lambda \asymp n^{-\frac{2m}{2m+1}}$.

**Corollary 4.8.** *Suppose that Assumption A2 holds for $\lambda \asymp \lambda^\dagger$ and $m > 3/2$. Then $\|\widehat{f}_R - f_0\|_n^2 = \mathcal{O}_P(n^{-\frac{2m}{2m+1}})$.*

We next proceed to find another lower bound for $s$ to achieve optimal testing. By Lemma 4.7, $\sigma_{n,\lambda} \asymp \sqrt{s_\lambda}/n$, which leads to $d_{n,\lambda} \asymp \sqrt{\lambda + \sqrt{s_\lambda}/n}$. The optimal separation rate $d^*_n$ can be achieved by another type of trade-off, i.e., the bias of $\widehat{f}_R$ v.s. the standard derivation of $T_{n,\lambda}$, as follows

$$d^{*\,2}_n = \operatorname{argmin}\Big\{\lambda \,:\, \lambda > \sqrt{s_\lambda}/n\Big\}.$$

Hence, $\lambda^* = n^{-\frac{4m}{4m+1}}$, $s^* \asymp (\lambda^*)^{-\frac{1}{2m}} = n^{\frac{2}{4m+1}}$, and the corresponding $\kappa_{\lambda^*} \asymp n^{\frac{1}{4m+1}}$ as defined in (3.2). In summary, $T_{n,\lambda}$ achieves the minimax optimal rate of testing $n^{-\frac{2m}{4m+1}}$ ([20]) when $s \gtrsim n^{\frac{2}{4m+1}}$ and $\lambda \asymp n^{-\frac{4m}{4m+1}}$.

**Corollary 4.9.** *Suppose that Assumption A2 holds with $\lambda \asymp \lambda^*$, $c_1, c_2 \in (0, \infty]$ and $m > 3/2$. Then for any $\varepsilon > 0$, there exist constants $C_\varepsilon$ and $N_\varepsilon$ such that, with probability greater than $1 - e^{-c_m n^{(2m-3)/(2m-1)}} - e^{-c_1 s} - e^{-c_2 s^*}$,*

$$\inf_{n \geq N_\varepsilon} \inf_{\substack{f \in \mathcal{B} \\ \|f\|_n \geq C_\varepsilon n^{-\frac{2m}{4m+1}}}} P_f(\phi_{n,\lambda} = 1 | \boldsymbol{x}, S) \geq 1 - \varepsilon.$$

It is worth emphasizing that $\lambda^\dagger, s^\dagger$ are different from $\lambda^*, s^*$, indicating a fundamental difference between estimation and testing. A more explicit reason for such a difference in minimax rate is due to two different types of trade-off, as illustrated in Figure 2.

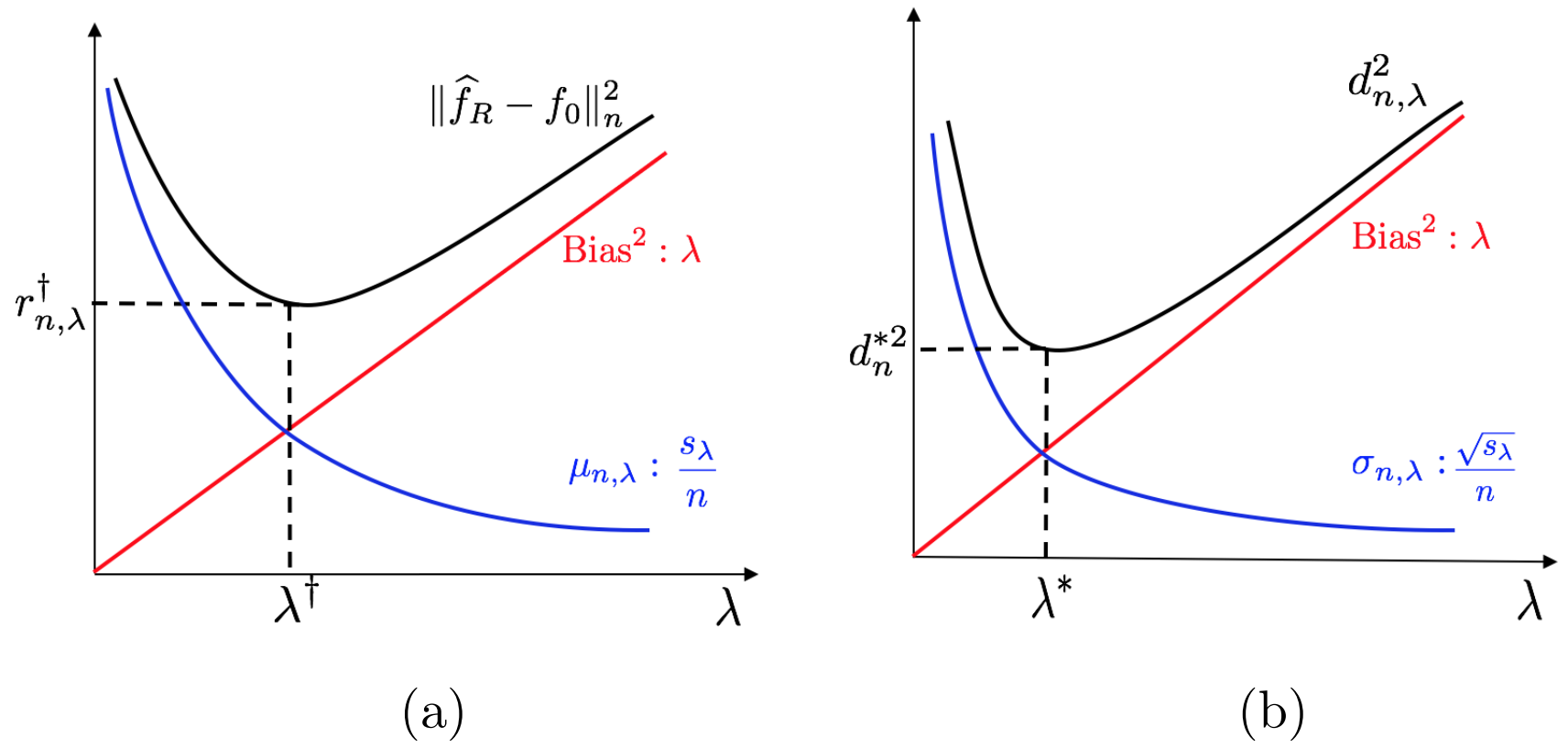

Figure 2: Trade-offs for achieving $(a)$ optimal estimation rate; $(b)$ optimal testing rate.

**Example 2: EDK**

Suppose that $\mathcal{H}$ is generated by EDK with $\gamma > 0$, $p \geq 1$. Parallel to Lemma 4.7, we have the following technical result.

**Lemma 4.10.** *Suppose that Assumption A2 holds with $c_1, c_2 \in (0, \infty]$. Then with probability at least $1 - e^{-c_{\gamma,p} n (\log n)^{-2/p}} - e^{-c_1 s} - e^{-c_2 s_\lambda}$, it holds that $\mu_{n,\lambda} \asymp s_\lambda / n$ and $\sigma^2_{n,\lambda} \asymp s_\lambda / n^2$, where $c_{\gamma,p}$ is an absolute constant depending on $\gamma, p$.*

Similar to Corollaries 4.8 and 4.9, one can prove the following result based on Corollary 4.4 and Theorem 4.6, which show that $s^\dagger = (\log n)^{1/p}$ or $s^* = (\log n)^{1/p}$[1] are lower bounds of $s$ for optimal estimation and testing, respectively.

**Corollary 4.11.** *Suppose that Assumption A2 holds for $\lambda \asymp \lambda^\dagger = (\log n)^{1/p} n^{-1}$ and $c_1, c_2 \in (0, \infty]$. Then $\|\widehat{f}_R - f_0\|_n^2 = \mathcal{O}_P(n^{-1} (\log n)^{1/p})$.*

**Corollary 4.12.** *Suppose that Assumption A2 holds for $\lambda \asymp \lambda^* = (\log n)^{1/(2p)} n^{-1}$ and $c_1, c_2 \in (0, \infty]$. Then for any $\varepsilon > 0$, there exist constants $C_\varepsilon$ and $N_\varepsilon$ such that, with probability greater than $1 - e^{-c_{\gamma,p} n (\log n)^{-2/p}} - e^{-c_1 s} - e^{-c_2 s_\lambda}$,*

$$\inf_{n \geq N_\varepsilon} \inf_{\substack{f \in \mathcal{B} \\ \|f\|_n \geq C_\varepsilon n^{-\frac{1}{2}} (\log n)^{\frac{1}{4p}}}} P_f(\phi_{n,\lambda} = 1 | \boldsymbol{x}, S) \geq 1 - \varepsilon.$$

We conclude our findings of this section in the following Table 1.

[1] In fact, for EDK, $s^*_\lambda \asymp (\log n - \frac{1}{2p} \log\log n)^{1/p}$. For simplicity, we keep the main term $s^*_\lambda \asymp (\log n)^{1/p}$.

| | Estimation | | | Testing | | |
|---|---|---|---|---|---|---|
| | $\lambda^\dagger$ | $s^\dagger$ | $r^\dagger_{n,\lambda}$ | $\lambda^*$ | $s^*$ | $d_n^{*2}$ |
| PDK | $n^{-\frac{2m}{2m+1}}$ | $n^{\frac{1}{2m+1}}$ | $n^{-\frac{2m}{2m+1}}$ | $n^{-\frac{4m}{4m+1}}$ | $n^{\frac{2}{4m+1}}$ | $n^{-\frac{4m}{4m+1}}$ |
| EDK | $(\log n)^{1/p} n^{-1}$ | $(\log n)^{1/p}$ | $(\log n)^{1/p} n^{-1}$ | $(\log n)^{1/2p} n^{-1}$ | $(\log n)^{1/p}$ | $(\log n)^{1/(2p)} n^{-1}$ |

Table 1: Lower bound of $s$ and choice of $\lambda$ for optimal estimation or testing in PDK and EDK.

### 4.5 Sharpness of $s^\dagger$ and $s^*$

In this section, we will show that $s^*$ and $s^\dagger$ derived in PDK and EDK are actually sharp. For technical convenience, define

$$\delta_n = \begin{cases} n^{-\frac{2}{2m-1}}, & K \text{ is PDK} \\ (\log n)^{-2/p}, & K \text{ is EDK} \end{cases}$$

Our first result is about the sharpness of $s^\dagger$. Theorem 4.13 shows that when $s \ll s^\dagger$, there exists a true function $f$ such that $\|\widehat{f}_R - f\|_n^2$ is substantially slower than the optimal estimation rate. Our proof is constructive in the sense we set $S = U_s^\top$ (as in Example 4.3), and construct the above true function as $\sum_{i=1}^n K(x_i, \cdot) w_i$ with $w_i$ being defined in (7.7).

**Theorem 4.13.** *Suppose $s = o(s^\dagger)$. Then there exists an $s \times n$ random matrix $S$ satisfying Assumption A2, such that with probability greater than $1 - e^{-cn\delta_n} - e^{-c_1 s} - e^{-c_2 s_\lambda}$, it holds that*

$$\sup_{f \in \mathcal{B}} \|\widehat{f}_R - f\|_n^2 \gg r^\dagger_{n,\lambda},$$

*where $c$ is a constant independent of $n$, and $c_1, c_2 \in (0, \infty]$ are given in Assumption A2 (b).*

Our second result is about the sharpness of $s^*$. Theorem 4.14 shows that when $s \ll s^*$, there exists a local alternative $f$ that is not detectable by $T_{n,\lambda}$ even when it is separated from zero by $d_n^*$. In this case, the asymptotic testing power is actually smaller than $\alpha$. The proof of Theorem 4.14 is similar as that of Theorem 4.13, except that a different true function (as defined in (7.9)) is constructed.

**Theorem 4.14.** *Suppose $s = o(s^*)$. Then there exists an $s \times n$ projection matrix $S$ satisfying Assumption A2 and a positive nonrandom sequence $\beta_{n,\lambda}$ satisfying $\lim_{n\to\infty} \beta_{n,\lambda} = \infty$ such that, with probability at least $1 - e^{-cn\delta_n} - e^{-c_1 s} - e^{-c_2 s_\lambda}$,*

$$\limsup_{n\to\infty} \inf_{\substack{f \in \mathcal{B} \\ \|f\|_n \ge \beta_{n,\lambda} d_n^*}} P_f(\phi_{n,\lambda} = 1 | \boldsymbol{x}, S) \le \alpha,$$

*where $c$ is a constant independent of $n$, and $c_1, c_2 \in (0, \infty]$ are given in Assumption A2 (b). Recall $1 - \alpha$ is the significance level.*

In view of Theorems 4.6 and 4.14, we observe a subtle phase transition phenomenon for testing signals as shown in Figure 1.

## 5 Adaptive Testing

In this section, we focus on the case of PDK as a leading example, and construct an adaptive testing procedure that does not require any exact prior knowledge on $m$ except for $m \geq 2$. The adaptive procedure is proven to achieve the minimax rate of testing established by [39] (up to an iterative-logarithmic term).

Consider an RKHS generated by a PDK of order $m_* \geq 2$, i.e., $\mathcal{H} = \mathcal{H}_{m_*}$. To reflect the role of $m$, we modify all previous notation by adding a subscript $m$. For example, let $K_m(\cdot, \cdot)$ be the reproducing kernel function associated with $\mathcal{H}_m$, and $\boldsymbol{K}_m = \frac{1}{n}[K_m(x_i, x_j)]_{1\leq i,j\leq n}$ be the corresponding empirical kernel matrix. Let $S_m$ be an $s_m \times n$ projection matrix. We will construct the corresponding $\widehat{f}_{R,m}(\cdot)$ based on (2.4) under $S_m$ and $\lambda_m$. Here

$$\lambda_m = cn^{-4m/(4m+1)}(\log\log n)^{2m/(4m+1)},$$

and the corresponding projection dimension $s_m$ is an integer satisfying

$$s_m \geq dn^{2/(4m+1)}(\log\log n)^{-1/(4m+1)}, \tag{5.1}$$

where $d > 0$ is a sufficiently large constant.

Given each $m$, the test statistic is defined as

$$T_{n,m} \equiv \|\widehat{f}_{R,m}\|_n^2 = \frac{1}{n}\boldsymbol{y}^\top \Delta_m^2 \boldsymbol{y}. \tag{5.2}$$

Based on $T_{n,m}$, our adaptive testing procedure is constructed as follows.

Step 1. For any $2 \leq m \leq m_n \to \infty$, standardize $T_{n,m}$ as

$$\tau_m = \frac{nT_{n,m} - \mathrm{tr}(\Delta_m^2)}{\sqrt{2\,\mathrm{tr}(\Delta_m^4)}}.$$

Step 2. Calculate $\tau_n^* = \max_{1\leq m\leq m_n} \tau_m$.

Step 3. Find $\tau_{n,m_n} = B_n(\tau_n^* - B_n)$, where $B_n$[2] satisfies

$$2\pi B_n^2 \exp(B_n^2) = m_n^2. \tag{5.3}$$

[2]According to [12], $B_n$ satisfying (5.3) has an approximation

$$\begin{aligned} B_n &= \sqrt{2\log m_n} - \frac{1}{2}(\log\log m_n + \log 4\pi)/\sqrt{2\log m_n} + O(1/\log m_n) \\ &\asymp \sqrt{2\log m_n}. \end{aligned}$$

By allowing $m_n \to \infty$, the unknown $m_*$ will be eventually covered over a sequence of test statistics. Under the null hypothesis (4.1), $T_{n,m} = \frac{1}{n}\epsilon^\top \Delta_m^2 \epsilon$, and thus $\tau_m$ is of a standardized quadratic form. Then, $\tau_n^*$ is the maxima of a sequence of *dependent* $\tau_m$'s. Based on a recent Gaussian approximation result in [24], i.e., Lemma S.3, we prove in the following Theorem 5.1 that the null limit distribution of $\tau_{n,m_n}$ is some extreme value distribution.

**Theorem 5.1.** *Suppose that* $m_n \asymp (\log n)^{d_0}$ *for a constant* $d_0 \in (0, 1/2)$*, and, for* $2 \le m \le m_n$*,* $S_m$ *satisfies Assumption A2 (b) with projection dimension* $s_m$*. Then, under* $H_0$ *in (4.1), for any* $\alpha \in (0, 1)$*, it holds that*

$$P(\tau_{n,m_n} \le c_\alpha) \to 1 - \alpha, \quad \text{as } n \to \infty,$$

*where* $c_\alpha = -\log(-\log(1-\alpha))$*.*

Our next result states that the above adaptive testing procedure is asymptotically minimax optimal. Specifically, Theorem 5.2 shows that $\tau_{n,m_n}$ achieves high power if the local alternative is separated from zero by an order $\delta(n, m_*)$ defined as

$$\delta(n, m_*) \equiv n^{-2m_*/(4m_*+1)} (\log\log n)^{m_*/(4m_*+1)}. \tag{5.4}$$

And, [39] showed that $\delta(n, m_*)$ is minimax optimal rate for adaptive testing.

**Theorem 5.2.** *Suppose that* $m_n \asymp (\log n)^{d_0}$ *for a constant* $d_0 \in (0, 1/2)$*, and* $S_m$ *satisfies Assumption A2 (b) with projection dimension* $s_m$*. Then, for any* $\varepsilon > 0$*, there exist positive constants* $C_\varepsilon, N_\varepsilon$ *for any* $n \ge N_\varepsilon$*, with probability approaching 1,*

$$\inf_{\substack{f \in \mathcal{B}_{n,m_*} \\ \|f\|_n \ge C_\varepsilon \delta(n,m_*)}} P_f(\tau_{n,m_n} \ge c_{\bar\alpha} | \boldsymbol{x}, S) \ge 1 - \varepsilon,$$

*where* $\mathcal{B}_{n,m_*} = \{f \in \mathcal{H}_{m_*} : (\boldsymbol{f})^\top \boldsymbol{K}_{m^*}^{-1} \boldsymbol{f} \le 1\}$ *and* $\boldsymbol{f} = (f(x_1), \dots, f(x_n))^\top$*.*

In the end, we point out that the lower bound for $s_m$ given in (5.1) is slightly smaller than the sharp lower bound for $s$ derived in the non-adaptive case; see Table 1. This is not surprising since the corresponding minimax rate $\delta(n, m_*)$, i.e., (5.4), is larger than the non-adaptive rate, i.e., $n^{-2m_*/(4m_*+1)}$.

# 6 Numerical Study

In this section, we examine the performance of the proposed testing procedure through simulation studies in Sections 6.1 and 6.2, and through a real data set in Section 6.3.

### 6.1 Simulation Study I: PDK

Data were generated from the regression model (1.1) with $f(x) = c(3\beta_{30,17}(x) + 2\beta_{3,11}(x))$, where $\beta_{a,b}$ is the density function for Beta$(a, b)$, $x_i \overset{iid}{\sim} \text{Unif}[0, 1]$, $\epsilon_i \overset{iid}{\sim} N(0, 1)$ and $c$ is a constant. To fit the model, we consider the periodic Sobolev kernel with eigenvalues $\mu_{2i} = \mu_{2i-1} = (2\pi i)^{-2m}$ for $i \geq 1$; see [16] for details. Set $n = 2^9, 2^{10}, 2^{11}, 2^{12}$, and $H_0 : f = 0$. The significance level was chosen as 0.05 and the Gaussian random projection matrix was applied in this setting.

We examined the empirical performance of the distance-based test (DT) $T_{n,\lambda}$, and adaptive test (AT) $\tau_{n,m_n}$. For DT, the projection dimension $s$ was chosen as $2n^\gamma$ for $\gamma = 1/(4m + 1), 2/(4m + 1), 3/(4m + 1)$, with $m = 2$ corresponding to cubic splines. For AT, the projection dimensions $s_m$ was chosen as $2n^\gamma(\log\log n)^{-\frac{1}{4m+1}}$ for $m = 2, \cdots, \sqrt{\log n}$. The smoothing parameter $\lambda$ was chosen by a projection version of Wahba's GCV score ([42]) based on $\widehat{f}_R$. Specifically, our new GCV score is defined as follows

$$V_S(\lambda) = \frac{\frac{1}{n}\boldsymbol{Y}^\top(I - A_S(\lambda))^2\boldsymbol{Y}}{\{\frac{1}{n}\operatorname{tr}(I - A_S(\lambda))\}^2}, \ \lambda > 0,$$

where $A_S = \boldsymbol{K}S^\top(S\boldsymbol{K}^2S^\top + \lambda S\boldsymbol{K}S^\top)^{-1}S\boldsymbol{K}$ is a projection version of the classical smoothing matrix. Our new GCV score enjoys much computational advantage than the classical one in our empirical study.

Empirical size was evaluated at $c = 0$, and power was evaluated at $c = 0.01, 0.02, 0.03$. Both size and power were calculated based on 500 independent replications. Figure 3 shows that the size of both DT and AT approach the nominal level 0.05 under various choices of $(s, n)$, demonstrating the validity of the proposed testing procedure. Figure 4 displays the power of DT and AT. Under various choices of $c$ and $\gamma$, it is not surprising to see from Figure 4 $(a)$, $(c)$, and $(e)$ that the power of DT approaches one as $n$ or $c$ increases. Rather, a key observation is that the power cannot be further improved as $\gamma$ grows beyond the critical point $2/(4m + 1)$ when $c \geq 0.02$. This is consistent with our theoretical result; see Theorem 4.14. Similar patterns have been observed for the power of AT in Figure 4 $(b)$, $(d)$, and $(f)$. Of course, the power of AT is usually lower than that of DT under the same setup, especially when the signal strength is weak. This is the price paid for adaptivity.

### 6.2 Simulation Study II: EDK

In this section, we consider a multivariate case and test $H_0 : f = 0$. Data were generated from

$$y_i = c(x_{i1}^2 + 2x_{i1}x_{i2} + 4x_{i1}x_{i2}x_{i3}) + \epsilon_i, \ i = 1, \cdots, n,$$

where $(x_{i1}, x_{i2}, x_{i3})$ follows from $N(\mu, I_3)$ with $\mu = (0, 0, 0)$, $\epsilon_i \sim N(0, 1)$, and $c \in \{0, 0.05, 0.1, 0.15\}$. Specifically, we chose the Gaussian kernel

$$K(\mathbf{x}, \mathbf{x}') = e^{-\frac{1}{2}\sum_{i=1}^3(x_i - x_i')^2}.$$

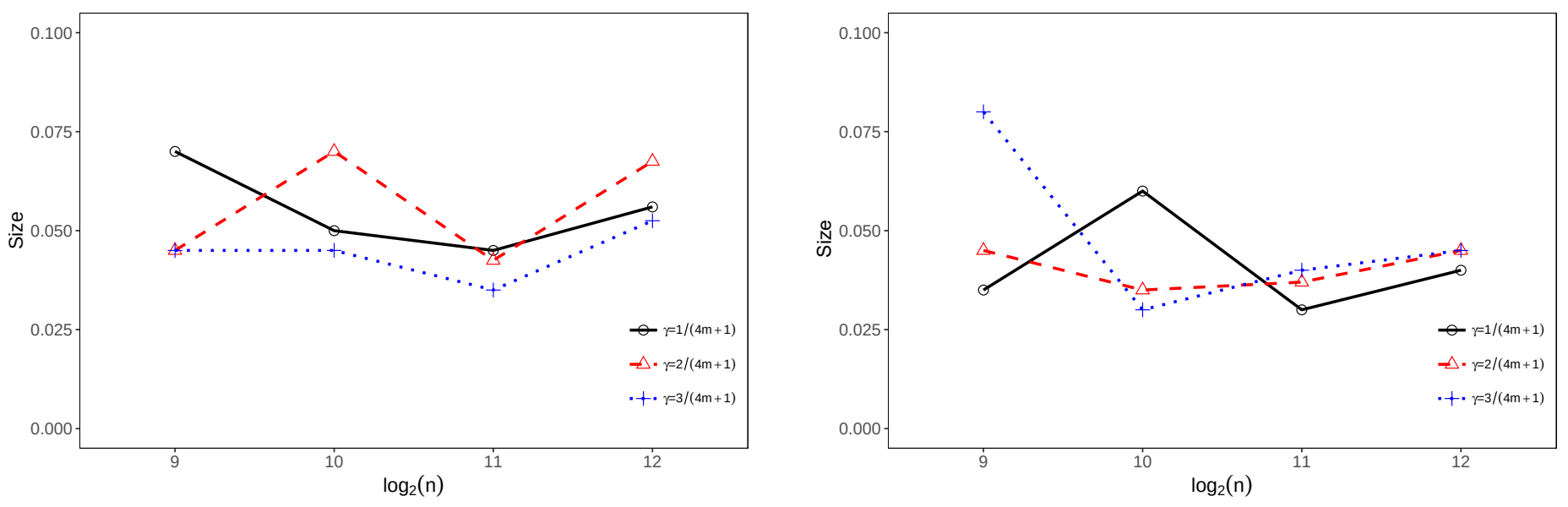

$(a)$DT $(b)$AT

Figure 3: Size for $(a)$ DT and $(b)$ AT with projection dimension varies. Signal strength $c = 0$.

We considered sample sizes $n = 2^9$ to $n = 2^{12}$ and sketch dimensions $s = 1.2\log(n), 1.2(\log n)^{3/2}, 1.2(\log n)^2$. For each pair $(n, s)$, experiments were independently repeated 500 times for calculating the size and power.

Interpretations for Figure 5 about the size and power are similar to those for Figures 3 and 4. Interestingly, we observe that the power increases dramatically as $\gamma$ increases from 1 to 1.5, while becomes stable near one as $\gamma \geq 1.5$. This is consistent with Corollary 4.12. Figure 6 demonstrates the significant computational advantage of DT (corresponding to $\gamma < 1$) over the testing procedure based on standard KRR (corresponding to $\gamma = 1$).

In the supplementary, we conduct additional synthetic experiments under the same simulation setup as Sections 6.1 and 6.2 except for using the Bernoulli random matrix. As shown in Figures 8-10, the interpretations remain the same.

### 6.3 Real Data Analysis: Air Quality Data

Urban air pollution is listed as one of the world's worst toxic pollution problems, which occurs when harmful substances including particulates and biological molecules are introduced into the Earth's atmosphere. In this section, we aim to analyze the association between environmental factors and pollutants, particularly, PM 2.5, which are airborne particles with aerodynamic diameters less than $2.5\mu m$.

Our analysis uses hourly PM2.5 readings taken from the US Embassy in Beijing located at (116.47E, 39.95N), in conjunction with hourly meteorological measurements at Beijing Capital International Airport (BCIA) in 2010-2012, which can be obtained from `weather.nocrew.org`. The data set contains information of 16,123 records with 9 variables such as PM2.5 concentration ($Y$), Hours ($X_1$), Dew Point ($X_2$), Temperature ($X_3$), Pressure ($X_4$), Cumulated wind speed ($X_5$), Cumulated hours of snow ($X_6$) and Cumulated hours of rain ($X_7$). We first tested the overall effect based on the multivariate Gaussian kernel and found that the p-value is less than $10^{-10}$. To further

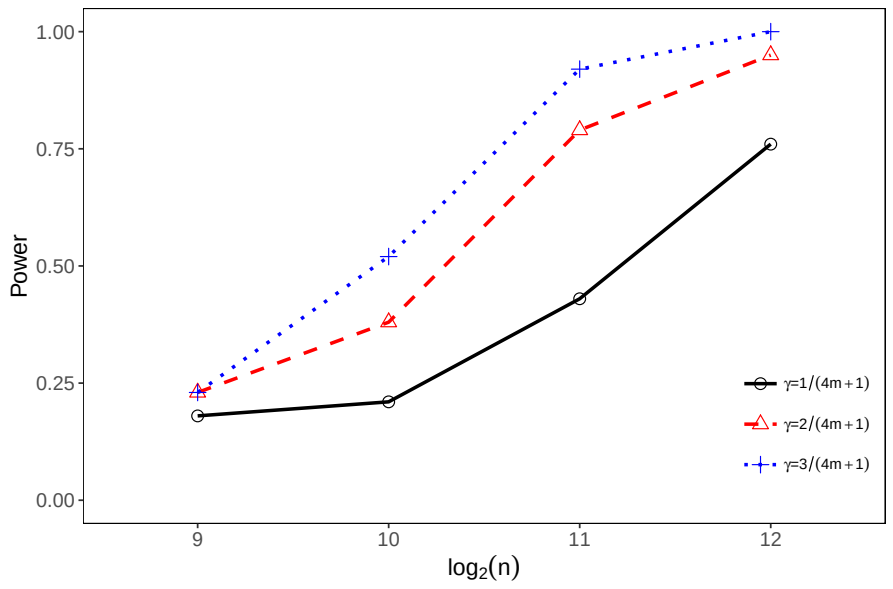

$(a)$ DT

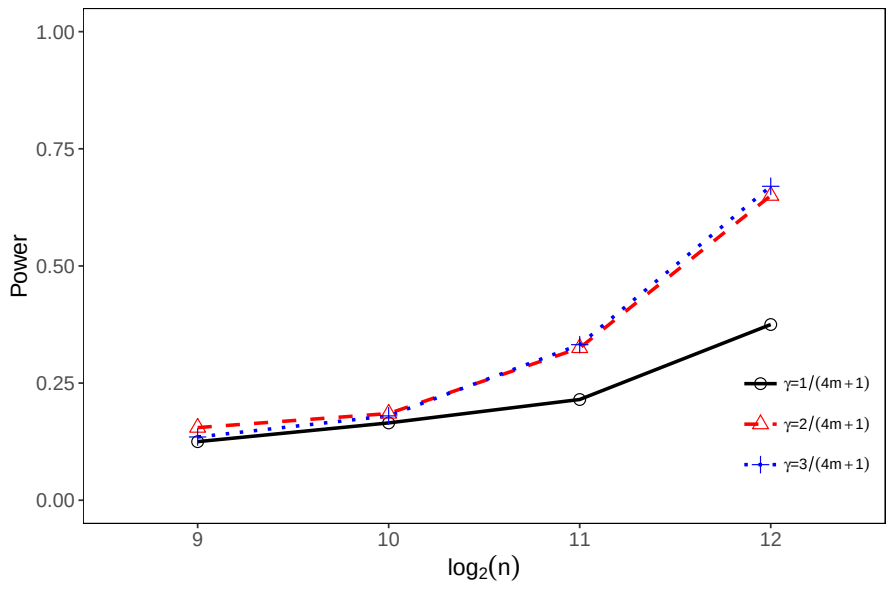

$(b)$AT

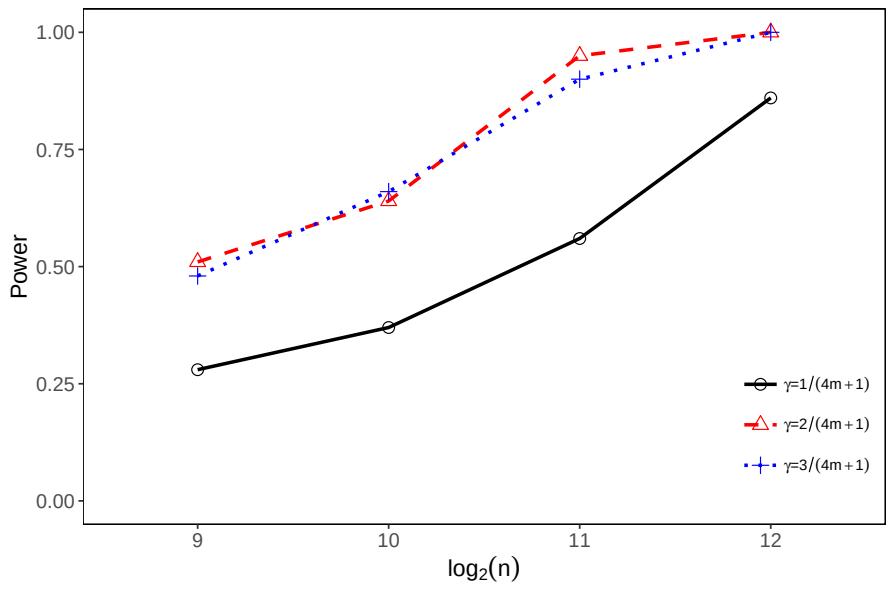

$(c)$DT

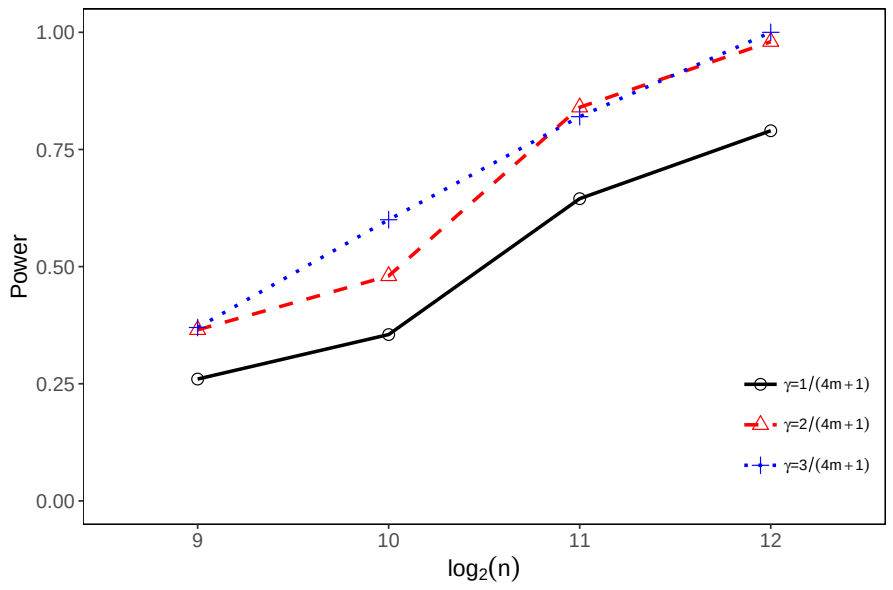

$(d)$AT

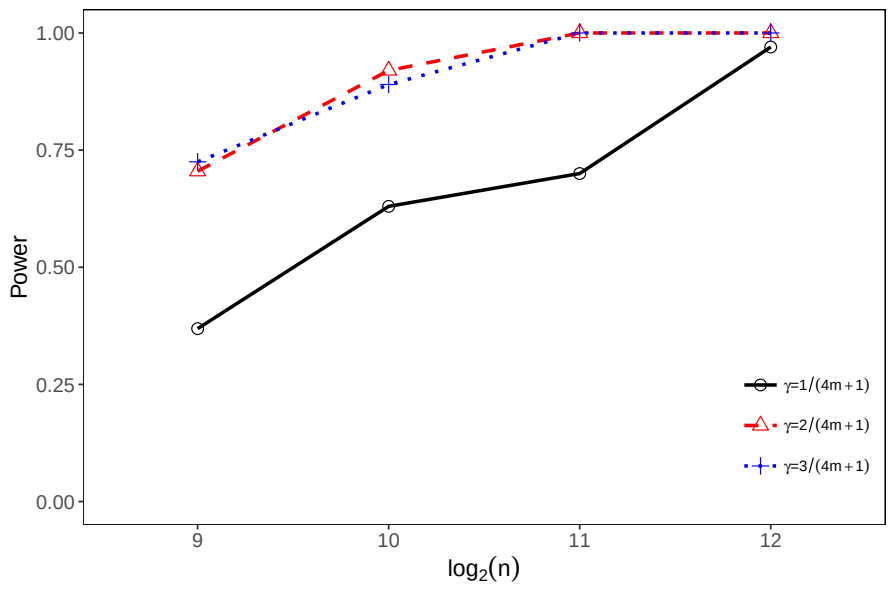

$(e)$DT

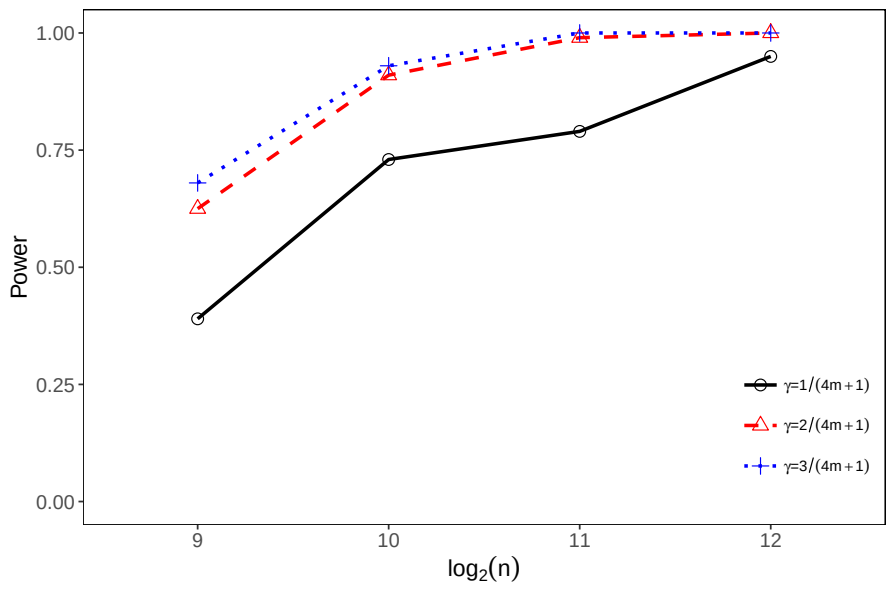

$(f)$AT

Figure 4: Power for DT and AT with varying projection dimension. Signal strength $c = 0.01$ for $(a)$ and $(b)$; $c = 0.02$ for $(c)$ and $(d)$; $c = 0.03$ for $(e)$ and $(f)$.

evaluate the association between each environmental factor and the concentration of PM2.5, we applied the proposed testing procedure with a smoothing spline kernel for each covariate. We found that all the p-values are less than $10^{-16}$.

To better understand the relationship between $Y$ and each covariate, we further tested whether the effect of $X_i$ on $Y$ is linear or nonlinear based on Theorem 4.5. As shown in Figure 7, PM2.5 reduces linearly as the cumulated hours of rain or snow increases. And, the other environmental factors all show strong nonlinear effect on the concentration of PM2.5. Hence, we suggest further fitting a semi-parametric model with $X_6, X_7$ the linear components to explore the relationship between PM2.5 and the 9 covariates.

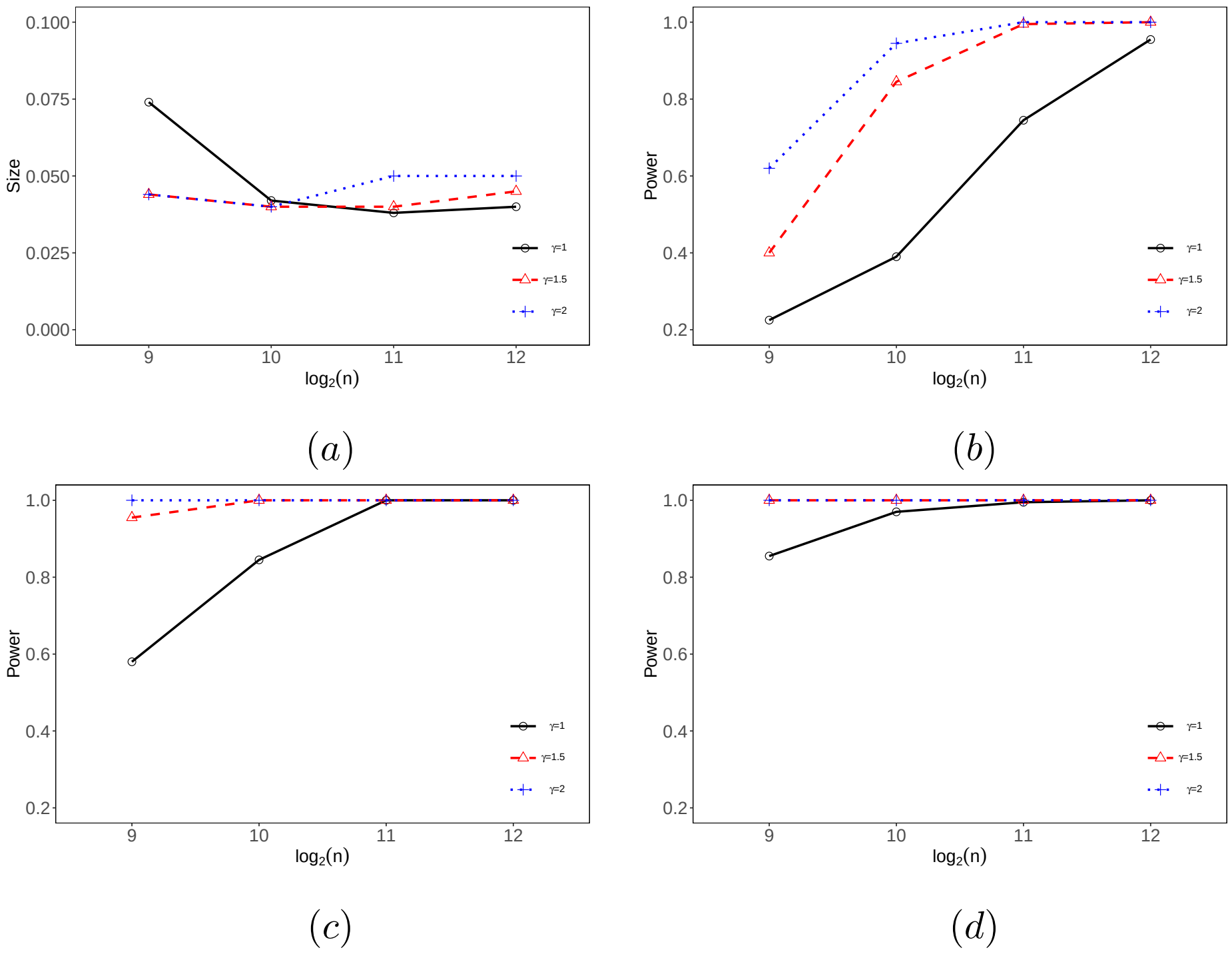

Figure 5: Size and power for DT with varying projection dimensions. Signal strength $c = 0$ for $(a)$; $c = 0.05$ for $(b)$; $c = 0.1$ for $(c)$; $c = 0.15$ for $(d)$.

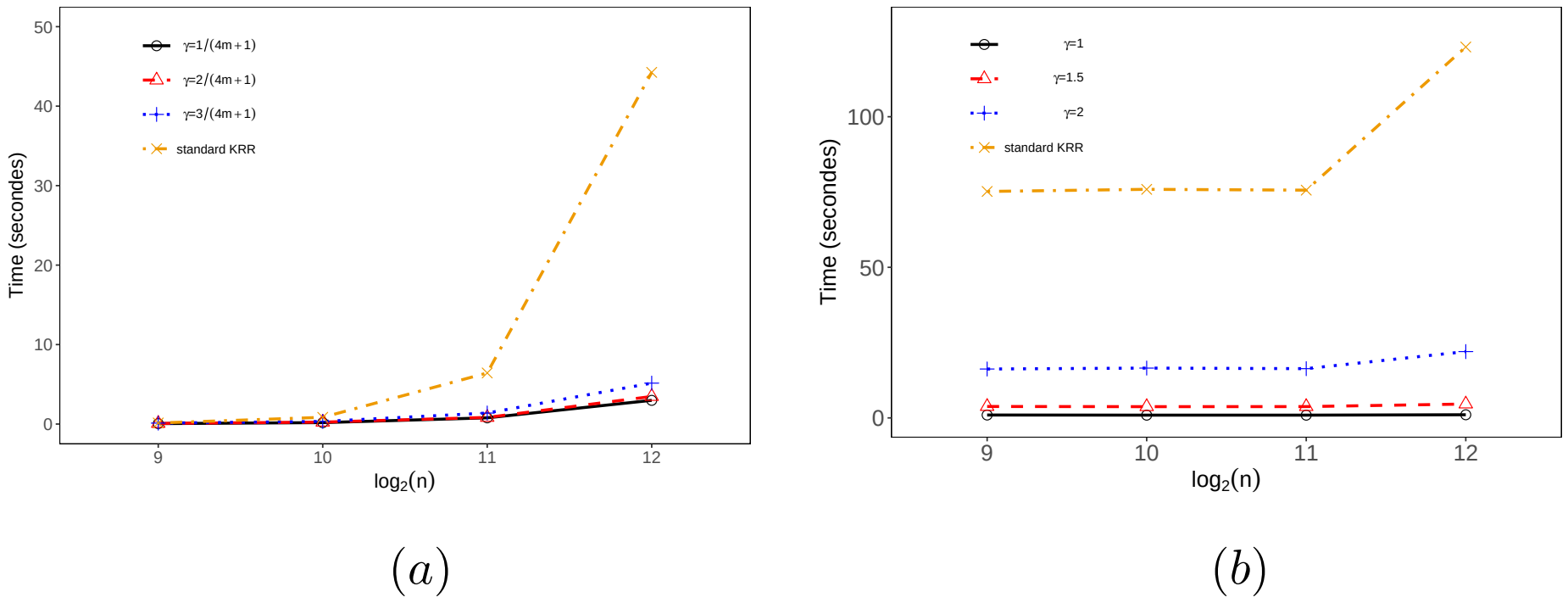

Figure 6: Computing time for DT with varying projection dimensions: $(a)$ is polynomially decay kernels; $(b)$ is exponentially decay kernels.

# 7 Proof of main results

In this section, we present main proofs of Theorem 4.2, Theorem 4.5, Theorem 4.6, Lemma 4.7, Lemma 4.10, Theorem 4.13, Theorem 4.14 in the main text.

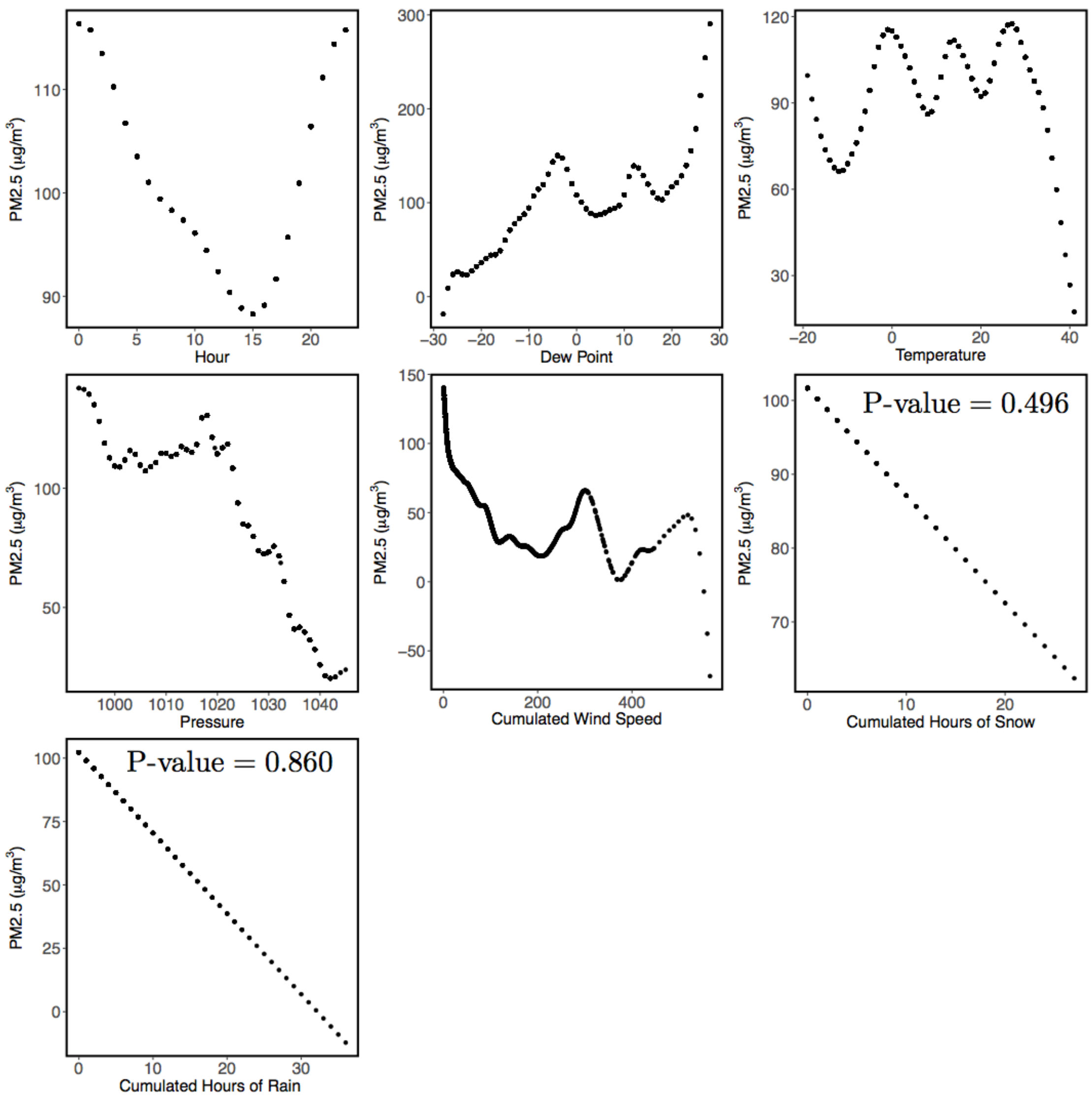

Figure 7: *P-values of the testing procedure for* 6 *predictors and the fitted curve, respectively.*

### 7.1 Proof of Theorem 4.2

*Proof.* Let $\Delta = \boldsymbol{K}S^\top(S\boldsymbol{K}^2S^\top + \lambda S\boldsymbol{K}S^\top)^{-1}S\boldsymbol{K}$, under the null hypothesis, $T_{n,\lambda} = \frac{1}{n}\epsilon^\top\Delta^2\epsilon$. We first derive the testing consistency of $T_{n,\lambda}$ conditional on $\boldsymbol{x} = (x_1, \cdots, x_n)$ and $S$. By the Gaussian assumption of $\epsilon$, we have $\mu_{n,\lambda} \equiv \mathrm{E}\left(T_{n,\lambda}|\boldsymbol{x}, S\right) = \frac{\mathrm{tr}(\Delta^2)}{n}$ and $\sigma^2_{n,\lambda} \equiv \mathrm{Var}\left(T_{n,\lambda}|\boldsymbol{x}, S\right) = 2\,\mathrm{tr}(\Delta^4)/n^2$.

Define $U = \frac{T_{n,\lambda}-\mu_{n,\lambda}}{\sigma_{n,\lambda}}$, then for any $t \in (-1/2, 1/2)$, we have

$$
\begin{aligned}
&\log \mathrm{E}_\epsilon \left( \exp(itU) \right) \\
=&\log \mathrm{E}_\epsilon \left( \exp(it\epsilon^\top \Delta^2 \epsilon/(n\sigma_{n,\lambda})) \right) - it\mu_{n,\lambda}/(n\sigma_{n,\lambda}) \\
=&-\frac{1}{2}\log\det(I_n - 2it\Delta^2/(n\sigma_{n,\lambda})) - it\mu_{n,\lambda}/(n\sigma_{n,\lambda}) \\
=&it \cdot \mathrm{tr}(\Delta^2)/(n\sigma_{n,\lambda}) - t^2\,\mathrm{tr}(\Delta^4)/(n^2\sigma_{n,\lambda}^2) + \mathcal{O}(t^3\,\mathrm{tr}(\Delta^6)/(n^3\sigma_{n,\lambda}^3)) \\
&- it\mu_{n,\lambda}/(n\sigma_{n,\lambda}) \\
=&-t^2/2 + \mathcal{O}(t^3\,\mathrm{tr}(\Delta^6)/(n^3\sigma_{n,\lambda}^3)),
\end{aligned}
$$

where $i = \sqrt{-1}$, $\mathrm{E}_\epsilon$ is the expectation with respect to $\epsilon$, and $I_n$ is $n \times n$ identity matrix. Therefore, to prove the normality of $U$, we need to show $\mathrm{tr}(\Delta^6)/(n^3\sigma_{n,\lambda}^3) = o(1)$. Note that

$$
\frac{\mathrm{tr}(\Delta^6)}{(n^3\sigma_{n,\lambda}^3)} \asymp \frac{\mathrm{tr}(\Delta^6)}{\mathrm{tr}(\Delta^4)} \cdot \frac{1}{\sqrt{\mathrm{tr}(\Delta^4))}}
$$

where $\mathrm{tr}(\Delta^6) = \mathrm{tr}\left((I + \lambda(S\boldsymbol{K}^2S^\top)^{-1}S\boldsymbol{K}S^\top)^{-6}\right)$ and $\mathrm{tr}(\Delta^4) = \mathrm{tr}\left((I + \lambda(S\boldsymbol{K}^2S^\top)^{-1}S\boldsymbol{K}S^\top)^{-4}\right)$. Since $\mathrm{tr}(\Delta^6)/\,\mathrm{tr}(\Delta^4) < 1$, it is sufficient to prove $\frac{1}{\mathrm{tr}(\Delta^4)} = o(1)$ as $n \to \infty$.

Let $(S\boldsymbol{K}^2S^\top)^{-1}S\boldsymbol{K}S^\top = P\Lambda P^{-1}$, where $\Lambda$ is an $s \times s$ diagonal matrix, then

$$
\mathrm{tr}(\Delta^4) = \mathrm{tr}\left((I + \lambda\Lambda)^{-4}\right) = \sum_{i=1}^{s}(1 + \lambda\Lambda_i)^{-4},
$$

where $\Lambda_i$ is the $i$th diagonal element in $\Lambda$. Next we show $\lambda\Lambda$ has at least $\min\{s, \widehat{s}_\lambda\}$ bounded eigenvalues. Notice that $(S\boldsymbol{K}^2S^\top)^{-1}S\boldsymbol{K}S^\top$ has the same non-zero eigenvalues as $\boldsymbol{K}^{1/2}S^\top(S\boldsymbol{K}^2S^\top)^{-1}S\boldsymbol{K}^{1/2}$. For $\boldsymbol{K} = UDU^\top$, let $U = (U_s, U_{n-s})$, $D = (D_s, D_{n-s})$ with $D_s = \mathrm{diag}\{\widehat{\mu}_1, \cdots, \widehat{\mu}_s\}$, $D_{n-s} = \mathrm{diag}\{\widehat{\mu}_{s+1}, \cdots, \widehat{\mu}_n\}$. Let $\widetilde{S}_1 = SU_s$, $\widetilde{S}_2 = SU_{n-s}$, $\boldsymbol{K}^{1/2}S^\top(S\boldsymbol{K}^2S^\top)^{-1}S\boldsymbol{K}^{1/2}$ can be rewritten as the block matrix:

$$
\begin{aligned}
\boldsymbol{K}^{1/2}S^\top(S\boldsymbol{K}^2S^\top)^{-1}S\boldsymbol{K}^{1/2} &= \begin{pmatrix} D_s^{1/2}\widetilde{S}_1^\top \\ D_{n-s}^{1/2}\widetilde{S}_2^\top \end{pmatrix} (S\boldsymbol{K}^2S^\top)^{-1} \begin{pmatrix} \widetilde{S}_1 D_s^{1/2} & \widetilde{S}_2 D_{n-s}^{1/2} \end{pmatrix} \\
&= \begin{pmatrix} A_1 & A_2 \\ A_3 & A_4 \end{pmatrix}
\end{aligned}
$$

where $A_1 = D_s^{1/2}\widetilde{S}_1^\top(S\boldsymbol{K}^2S^\top)^{-1}\widetilde{S}_1 D_s^{1/2}$. By Lemma S.5 of the eigenvalue interlacing for principal submatrices theorem, we only need to prove $\lambda A_1$ has at least $\min\{s, s_\lambda\}$ bounded eigenvalues. Using Binomial Inverse Theorem,

$$
(S\boldsymbol{K}^2S^\top)^{-1} = (\widetilde{S}_1 D_s^2 \widetilde{S}_1^T)^{-1} - (\widetilde{S}_1 D_s^2 \widetilde{S}_1^T)^{-1}\Gamma(\widetilde{S}_1 D_s^2 \widetilde{S}_1^T)^{-1}, \tag{7.1}
$$

where $\Gamma$ is a symmetric matrix defined as

$$\Gamma = \widetilde{S}_2 D_{n-s}^2 \widetilde{S}_2^\top \big(\widetilde{S}_2 D_{n-s}^2 \widetilde{S}_2^\top + \widetilde{S}_2 D_{n-s}^2 \widetilde{S}_2^\top (\widetilde{S}_1 D_1^2 \widetilde{S}_1^T)^{-1} \widetilde{S}_2 D_{n-s}^2 \widetilde{S}_2^\top\big)^{-1} \widetilde{S}_2 D_{n-s}^2 \widetilde{S}_2^\top.$$

Plugging (7.1) into $A_1$, we have

$$D_s^{1/2} \widetilde{S}_1^\top (SK^2S^\top)^{-1} \widetilde{S}_1 D_s^{1/2} = D_s^{-1} - H,$$

where $H$ is a semi-positive matrix. Based on Lemma S.6 of Weyl's inequality, the $i^{th}$ eigenvalue of $D_s^{-1}$ is greater than the $i^{th}$ eigenvalue of $A_1$. Recall $\widehat{s}_\lambda = \operatorname{argmin}\{i : \widehat{\mu}_i \leq \lambda\} - 1$, we have $\lambda/\widehat{\mu}_i \leq 1$ for $i = 1, \cdots, \widehat{s}_\lambda$. Hence, there exist at least $\min\{s, \widehat{s}_\lambda\}$ bounded eigenvalues for $\lambda A_1$. Finally, we have

$$\operatorname{tr}(\Delta^4) \geq C \min\{s, \widehat{s}_\lambda\}, \quad \text{where } C \text{ is some constant}, \tag{7.2}$$

then $\operatorname{E}_\epsilon(e^{itU}) \longrightarrow e^{-\frac{t^2}{2}}$, a.s.

We next consider $\operatorname{E}_{\boldsymbol{x},S} \operatorname{E}_\epsilon(e^{itU})$ by taking expectation w.r.t $\boldsymbol{x}, S$ on $\operatorname{E}_\epsilon(e^{itU})$. We claim $\operatorname{E}_{\boldsymbol{x},S} \operatorname{E}_\epsilon(e^{itU}) \longrightarrow e^{-\frac{t^2}{2}}$ for $t \in (-\frac{1}{2}, \frac{1}{2})$. If not, there exists a subsequence of r.v $\{\boldsymbol{x}_{n_k}, S_{n'_k}\}$, such that for $\forall \varepsilon > 0$, $|\operatorname{E}_{\boldsymbol{x}_{n_k}, S_{n_k}} \operatorname{E}_\epsilon e^{itU} - e^{-\frac{t^2}{2}}| > \varepsilon$. On the other hand, since $\operatorname{E}_\epsilon e^{itU(\boldsymbol{x}_{n_k}, S_{n_k})} \xrightarrow{P} e^{-\frac{t^2}{2}}$, which is bounded, there exists a sub-sub sequence $\{\boldsymbol{x}_{n_{k_l}}, S_{n_{k'_l}}\}$, such that

$$\operatorname{E}_\epsilon e^{itU(\boldsymbol{x}_{n_{k_l}}, S_{n_{k'_l}})} \xrightarrow{a.s} e^{-\frac{t^2}{2}}.$$

Thus by dominate convergence theorem, $\operatorname{E}_{\boldsymbol{x}_{n_{k_l}}, S_{n_{k'_l}}} \operatorname{E}_\epsilon e^{itU} \longrightarrow e^{-\frac{t^2}{2}}$, which is a contradiction. Therefore, we have $U = \frac{T_{n,\lambda} - \mu_{n,\lambda}}{\sigma_{n,\lambda}}$ asymptotically converges to a standard normal distribution. □

## 7.2 Proof of Theorem 4.5

*Proof.* Under $H_0^{\text{linear}}$, it can be shown that

$$T_{n,\lambda}^* = \frac{1}{n} \epsilon' (I - H) \Delta^2 (I - H) \epsilon$$

with $\mu_{n,\lambda}^* = \operatorname{tr}\big((I - H)\Delta^2(I - H)\big)/n$ and $\sigma_{n,\lambda}^* = \sqrt{2 \operatorname{tr}\big((I - H)\Delta^4(I - H)\big)}/n$.

Similar to Theorem 4.2, we only need to prove $\operatorname{tr}\big((I - H)\Delta^2(I - H)\big) \to \infty$ as $n \to \infty$. Notice that $\operatorname{tr}\big((I - H)\Delta^2(I - H)\big) = \operatorname{tr}\big(\Delta^2(I - H)\big) = \operatorname{tr}(\Delta^2) - \operatorname{tr}(\Delta^2 H)$. For $\operatorname{tr}(\Delta^2 H)$, we have $\operatorname{tr}(\Delta^2 H) = \operatorname{tr}(\Delta^2 H^2) = \operatorname{tr}(H \Delta^2 H)$. Since $\operatorname{rank}(H \Delta^2 H) \leq d + 1$, and $\lambda_{\max}(H \Delta H) \leq 1$, therefore $\operatorname{tr}(\Delta^2 H) \leq d + 1$. Here $d$ is the dimension of $x$. Finally we have $\operatorname{tr}\big((I - H)\Delta^2(I - H)\big) \geq \min\{s, \widehat{s}_\lambda\} - d - 1$. The last step is based on the proof of Theorem 4.2. □

### 7.3 Proof of Theorem 4.6

In this section, we prove the the testing is minimax optimal as stated in Theorem 4.6.

*Proof.*

$$
\begin{aligned}
& n\|\widehat{f}_R\|_n^2 \\
=& n\|\boldsymbol{K}S^\top(S\boldsymbol{K}^2S^\top + \lambda S\boldsymbol{K}S^\top)^{-1}S\boldsymbol{K}f + \\
& \boldsymbol{K}S^\top(S\boldsymbol{K}^2S^\top\lambda S\boldsymbol{K}S^\top)^{-1}S\boldsymbol{K}\epsilon\|_n^2 \\
=& n\|\operatorname{E}_\epsilon \widehat{f}_R\|_n^2 + n\|\boldsymbol{K}S^\top(S\boldsymbol{K}^2S^\top + \lambda S\boldsymbol{K}S^\top)^{-1}S\boldsymbol{K}\epsilon\|_n^2 \\
& + 2(\operatorname{E}_\epsilon \widehat{\boldsymbol{f}}_R)^\top \boldsymbol{K}S^\top(S\boldsymbol{K}^2S^\top + \lambda S\boldsymbol{K}S^\top)^{-1}S\boldsymbol{K}\epsilon \\
:=& T_1 + T_2 + T_3.
\end{aligned}
\tag{7.3}
$$

Lemma 4.3 shows that $\|f - \operatorname{E}_\epsilon \widehat{f}_R\|_n^2 \leq C\lambda$ with probability $1 - e^{-c_1 s} - e^{-c_2 s_\lambda}$. Set $C' = \sqrt{2C}$. Given the separation rate $\|f\|_n^2 \geq C'^2 d_{n,\lambda}^2 = 2C(\lambda + \sigma_{n,\lambda})$, we have

$$
T_1 = n\|\operatorname{E}_\epsilon \widehat{f}_R\|_n^2 \geq \frac{n}{2}\|f\|_n^2 - n\|f - \operatorname{E}_\epsilon \widehat{f}_R\|_n^2 \geq nC(\lambda + \sigma_{n,\lambda}) - nC\lambda \geq nC\sigma_{n,\lambda}
$$

with probability at least $1 - e^{-c_1 s} - e^{-c_2 s_\lambda}$, where $c_1, c_2$ is specified in Assumption A2.

Next, notice that $T_3 = (\operatorname{E}_\epsilon \widehat{\boldsymbol{f}}_R)^\top \Delta\epsilon$. Consider $a^\top\Delta^2 a$, where $a = (a_1, \cdots, a_n) \in \mathbb{R}^n$ is an arbitrary vector. Since $a^\top\Delta^2 a \leq \lambda_{\max}(\Delta^2)a^\top a$, where $\Delta^2$ has the same non-zero eigenvalue as

$$
\begin{aligned}
\widetilde{\Delta}^2 =& (S\boldsymbol{K}^2S^\top + \lambda S\boldsymbol{K}S^\top)^{-1}S\boldsymbol{K}^2S^\top(S\boldsymbol{K}^2S^\top + \lambda S\boldsymbol{K}S^\top)^{-1}S\boldsymbol{K}^2S^\top \\
=& (I + \lambda(S\boldsymbol{K}^2S^\top)^{-1}S\boldsymbol{K}S)^{-2},
\end{aligned}
$$

then we have $\|\widetilde{\Delta}^2\|_{\mathrm{op}} \leq 1$, and $\lambda_{\max}(\Delta^2) \leq 1$. Therefore,

$$
\operatorname{E}_\epsilon T_3^2 = (\operatorname{E}_\epsilon \widehat{\boldsymbol{f}}_R)^\top \Delta^2 (\operatorname{E}_\epsilon \widehat{\boldsymbol{f}}_R) \leq (\operatorname{E}_\epsilon \widehat{\boldsymbol{f}}_R)^\top (\operatorname{E}_\epsilon \widehat{\boldsymbol{f}}_R) = T_1,
$$

then

$$
\operatorname{P}\left(|T_3| \geq \varepsilon^{-\frac{1}{2}} T_1^{1/2} | \boldsymbol{x}, S\right) \leq \frac{\operatorname{E}_\epsilon T_3^2}{\varepsilon^{-1}T_1} \leq \varepsilon
$$

Define $\mathcal{E}_1 = \{T_1 \geq Cn\sigma_{n,\lambda}\}$, $\mathcal{E}_2 = \{\frac{T_2/n - \mu_{n,\lambda}}{\sigma_{n,\lambda}} \leq C_\epsilon\}$, where $C_\varepsilon$ satisfies $\operatorname{P}(\mathcal{E}_2) \geq 1 - \varepsilon$, $\mathcal{E}_3 = \{T_3 \geq$

$-\varepsilon^{-1/2}T_1^{1/2}\}$. Finally, with probability at least $1-e^{-c_1 s}-e^{-c_2 s_\lambda}$,

$$
\begin{aligned}
&\mathrm{P}_f\Big(\frac{\frac{1}{n}(T_1+T_2+T_3)-\mu_{n,\lambda}}{\sigma_{n,\lambda}}\geq z_{1-\alpha/2}\Big|\boldsymbol{x},S\Big)\\
\geq&\mathrm{P}_f\Big(\frac{T_1+T_2}{n\sigma_{n,\lambda}}+\frac{\frac{1}{n}T_2-\mu_{n,\lambda}}{\sigma_{n,\lambda}}\geq z_{1-\alpha/2},\,\mathcal{E}_1\cap\mathcal{E}_2\cap\mathcal{E}_3\Big|\boldsymbol{x},S\Big)\\
\geq&\mathrm{P}_f\Big(\frac{T_1(1-\varepsilon^{-1/2}T_1^{-1/2})}{n\sigma_{n,\lambda}}-C_\varepsilon\geq z_{1-\alpha/2},\,\mathcal{E}_1\cap\mathcal{E}_2\cap\mathcal{E}_3\Big|\boldsymbol{x},S\Big)\\
\geq&\mathrm{P}_f\Big(C(1-\frac{1}{\sqrt{Cn\sigma_{n,\lambda}\varepsilon}})-C_\varepsilon\geq z_{1-\alpha/2},\mathcal{E}_1\cap\mathcal{E}_2\cap\mathcal{E}_3\Big)\\
=&\mathrm{P}_f(\mathcal{E}_1\cap\mathcal{E}_2\cap\mathcal{E}_3)\\
\geq&1-3\varepsilon
\end{aligned}
$$

The second to the last equality is achieved by choosing $C$ to satisfy

$$
\frac{1}{\sqrt{Cn\sigma_{n,\lambda}\varepsilon}}<\frac{1}{2}\quad\text{and}\quad\frac{1}{2}C-C_\varepsilon\geq z_{1-\alpha/2}.
$$

□

### 7.4 Proof of Lemma 4.7

In this section, we analyze the orders of $\mu_{n,\lambda}$ and $\sigma_{n,\lambda}$ for PDK.

*Proof.* Recall in (7.2), we proved that $\mathrm{tr}(\Delta)\gtrsim\min\{s,\widehat{s}_\lambda\}$. Next we show with probability approaching 1,

$$
\mathrm{tr}(\Delta^4)\leq\mathrm{tr}(\Delta^2)\leq\mathrm{tr}(\Delta)\lesssim\widehat{s}_\lambda. \tag{7.4}
$$

On the other hand, when $\lambda\geq 1/n$, by Lemma S.2 (a), with probability at least $1-e^{-c_m n^{(2m-3)/(2m-1)}}$, $\widehat{s}_\lambda\asymp s_\lambda$. Combining (7.2) with (7.4), we have $\sigma^2_{n,\lambda}\asymp s_\lambda/n^2$ and $\mu_{n,\lambda}\asymp s_\lambda/n$ with probability approaching 1.

Note that $\mathrm{tr}(\Delta)=\mathrm{tr}(\widetilde{\Delta})$, where $\widetilde{\Delta}=DU^\top(S\boldsymbol{K}^2S^\top+\lambda S\boldsymbol{K}S^\top)^{-1}SUD$. $\widetilde{\Delta}$ can be written as

$$
\widetilde{\Delta}=\begin{pmatrix}\widetilde{\Delta}_1 & \widetilde{\Delta}_2\\ \widetilde{\Delta}_3 & \widetilde{\Delta}_4\end{pmatrix}
$$

with $\widetilde{\Delta}_1=D_s\widetilde{S}_1^\top(S\boldsymbol{K}^2S^\top+\lambda S\boldsymbol{K}S^\top)^{-1}\widetilde{S}_1D_s$, and $\widetilde{\Delta}_4=D_{n-s}\widetilde{S}_2^\top(S\boldsymbol{K}^2S^\top+\lambda S\boldsymbol{K}S^\top)^{-1}\widetilde{S}_2D_{n-s}$. Here $D_s=\mathrm{diag}\{\widehat{\mu}_1,\cdots,\widehat{\mu}_s\}$, $D_{n-s}=\mathrm{diag}\{\widehat{\mu}_{s+1},\cdots,\widehat{\mu}_n\}$, $\widetilde{S}_1=SU_s$, and $\widetilde{S}_2=SU_{n-s}$, where $U_s$ is the first $s$ column of $U$ and $U_{n-s}$ is the last $n-s$ column of $U$.

Let $\Lambda=\mathrm{diag}\{\Lambda_1,\Lambda_2\}=\mathrm{diag}\{D_s^2+\lambda D_s,D_{n-s}^2+\lambda D_{n-s}\}$. Then $\widetilde{\Delta}_1$ can be expressed as

$$
\begin{aligned}
\widetilde{\Delta}_1=&D_s\widetilde{S}_1^\top(\widetilde{S}_1\Lambda_1\widetilde{S}_1^\top+\widetilde{S}_2\Lambda_2\widetilde{S}_2^\top)^{-1}\widetilde{S}_1D_s\\
=&D_s^2\Lambda_1^{-1}-D_s\widetilde{S}_1^T(\widetilde{S}_1\Lambda_1\widetilde{S}_1^\top)^{-1}\big((\widetilde{S}_2\Lambda_2\widetilde{S}_2^\top)^{-1}+(\widetilde{S}_1\Lambda_1\widetilde{S}_1^\top)^{-1}\big)^{-1}(\widetilde{S}_1\Lambda_1\widetilde{S}_1^\top)^{-1}\widetilde{S}_1D_s.
\end{aligned}
$$

Therefore,

$$\operatorname{tr}\left(D_1\widetilde{S}_1^\top(S\boldsymbol{K}^2S^\top+\lambda S\boldsymbol{K}S^\top)^{-1}\widetilde{S}_1D_1\right)\le\operatorname{tr}(D_1^2\Lambda_1^{-1})\le s_\lambda. \tag{7.5}$$

The last inequality is deduced by the following step

$$\operatorname{tr}(D_1^2\Lambda_1^{-1})=\sum_{i=1}^{\widehat{s}_\lambda}\frac{\widehat{\mu}_i}{\widehat{\mu}_i+\lambda}+\sum_{i=\widehat{s}_\lambda+1}^{s}\frac{\widehat{\mu}_i}{\widehat{\mu}_i+\lambda}\le\widehat{s}_\lambda+\frac{1}{\lambda}\sum_{i=\widehat{s}_\lambda+1}^{s}\widehat{\mu}_i\le 2s_\lambda$$

with probability at least $1-e^{-c_mn^{(2m-3)/(2m-1)}}-e^{-c_1s}-e^{-c_2s_\lambda}$ by Lemma S.2(a) and Lemma 3.1.

Next, we consider $\widetilde{\Delta}_4$. Note that

$$\begin{aligned}
&\operatorname{tr}\left(D_{n-s}\widetilde{S}_2^\top(S\boldsymbol{K}^2S^\top+\lambda S\boldsymbol{K}S^\top)^{-1}\widetilde{S}_2D_{n-s}\right)\\
\le&\operatorname{tr}\left(D_{n-s}\widetilde{S}_2^\top(\lambda S\boldsymbol{K}S^\top)^{-1}\widetilde{S}_2D_{n-s}\right)\\
\le&\frac{\widehat{\mu}_s}{\lambda}\operatorname{tr}\left(\left((\widetilde{S}_2D_{n-s}\widetilde{S}_2^\top)^{-1}\widetilde{S}_1D_s\widetilde{S}_1^\top+I_{s\times s}\right)^{-1}\right)\le s\widehat{\mu}_s/\lambda.
\end{aligned}\tag{7.6}$$

We show $s\widehat{\mu}_s/(\lambda s_\lambda)\le C$ with probability approaching 1, where $C$ is some absolute constant. Then by (7.6), $\operatorname{tr}(\widetilde{\Delta}_4)\lesssim s_\lambda$. If $ds_\lambda\le s\le n^{1/(2m)}$, then by Lemma S.2(a), we have $\mu_s/2\le\widehat{\mu}_s\le 3\mu_s/2$ with probability at least $1-e^{-c_mn^{(2m-3)/(2m-1)}}$, where $m>3/2$. Then

$$\frac{s\widehat{\mu}_s}{\lambda s_\lambda}\le\frac{3s\mu_s}{2\lambda s_\lambda}\lesssim\frac{s^{1-2m}}{s_\lambda^{1-2m}}=\mathcal{O}(1).$$

If $s\gg n^{1/(2m)}$, based on the proof of Lemma S.2, in (S.3) and (S.4),

$$\mathrm{P}\left(|\widehat{\mu}_i-\mu_i|\ge\mu_i\widetilde{\varepsilon}+r^{1-2m}\right)\le 1-\exp\left(-c_m'n\widetilde{\varepsilon}^2/r^2\right).$$

Let $\widetilde{\varepsilon}=n^{-\frac{2m-1}{2m}}s^{2m-1}$ and $r=n^{\frac{1}{2m}}s^{\frac{1}{2m-1}}$, then

$$s\widehat{\mu}_s\le s\mu_s\widetilde{\varepsilon}+sr^{1-2m}=n^{-\frac{2m-1}{2m}},$$

with probability at least $1-e^{-c'n^{(2m-3)/(2m-1)}}$. The probability is obtained by calculating $n\widetilde{\varepsilon}^2/r^2$. Based on the assumption $1/n\le\lambda\le 1$, $\lambda s_\lambda\asymp\lambda^{1-1/(2m)}\ge n^{-\frac{2m-1}{2m}}$. Finally we have $\frac{s\widehat{\mu}_s}{\lambda s_\lambda}\le C$, i.e., $s\widehat{\mu}_s/\lambda\le Cs_\lambda$, where $C$ is some bounded constant.

Combining with (7.5), we have

$$\operatorname{tr}(\Delta)=\operatorname{tr}(\widetilde{\Delta})\le\operatorname{tr}(\widetilde{\Delta}_1)+\operatorname{tr}(\widetilde{\Delta}_2)\lesssim s_\lambda,$$

with probability at least $1-e^{-c'n^{(2m-3)/(2m-1)}}-e^{-c_1s}-e^{-c_2s_\lambda}$, where $c'$ is a constant only depending on $m$ and $c_1,c_2>0$ are defined in Assumption A2. □

### 7.5 Proof of Lemma 4.10

*Proof.* Following the same notation and strategy in the proof of Lemma 4.7, (7.5) also holds for EDK, with probability at least $1 - e^{-c_{\gamma,p}n(\log n)^{-2/p}} - e^{-c_1 s} - e^{-c_2 s_\lambda}$ by Lemma S.2 (a) and Lemma 3.1. For EDK, (7.6) also holds. Next we will prove that $\operatorname{tr}(\widetilde{\Delta}_2) \le s\widehat{\mu}_s/\lambda \le Cs_\lambda$, where $C$ is an absolute constant. If $ds_\lambda \le s \lesssim n^{1/2-\varepsilon} \triangleq n^a$ for any $0 < \varepsilon < 1/2$, then by Lemma S.2 (b), $\widehat{\mu}_s \le \frac{3}{2}\mu_s$. Therefore

$$\frac{s\widehat{\mu}_s}{\lambda s_\lambda} \lesssim \frac{s\mu_s}{s_\lambda \mu_{s_\lambda}} \le 1,$$

with probability at least $1-e^{-c_{\gamma,p}n(\log n)^{-2/p}}$, where the last inequality is by the fact that $s\mu_s \asymp se^{-\gamma s^p}$ is decreasing w.r.t $s$ when $\gamma p s^p - 1 > 0$. When $s \ge n^{1/2}$, by Lemma S.2 (b),

$$\frac{s\widehat{\mu}_s}{\lambda s_\lambda} \le \frac{s^2\mu_s}{\lambda s_\lambda} \lesssim \frac{ns^2}{s_\lambda} e^{-\gamma s^p} = o(1),$$

with probability at least $1-e^{-n}$. Thus, we achieve $\operatorname{tr}\left(D_{n-s}\widetilde{S}_2^\top (S\boldsymbol{K}^2S^\top + \lambda S\boldsymbol{K}S^\top)^{-1}\widetilde{S}_2 D_{n-s}\right) \le s_\lambda$ with probability at least $1 - e^{-c_{\gamma,p}n(\log n)^{-2/p}}$. Combining with (7.2) and (7.5), we have $\operatorname{tr}(\Delta) \lesssim s_\lambda$ with probability at least $1 - e^{-c_{\gamma,p}n(\log n)^{-2/p}} - e^{-c_1 s} - e^{-c_2 s_\lambda}$. □

### 7.6 Proof of Theorem 4.13

*Proof.* Notice that

$$\begin{aligned}
\|\widehat{f}_R - f_0\|_n^2 &= \|\operatorname{E}_\epsilon \widehat{f}_R - f_0\|_n^2 + \|\widehat{f}_R - \operatorname{E}_\epsilon \widehat{f}_R\|_n^2 + \frac{2}{n}\big(\widehat{f}_R - \operatorname{E}_\epsilon \widehat{f}_R\big)^\top\big(\operatorname{E}_\epsilon \widehat{f}_R - f_0\big) \\
&\equiv T_1 + T_2 + T_3.
\end{aligned}$$

We first consider $T_1$ as follows:

$$\begin{aligned}
\|\operatorname{E}_\epsilon \widehat{f}_R - f_0\|_n^2 &= \|\boldsymbol{K}S^\top(S\boldsymbol{K}^2S^\top + \lambda S\boldsymbol{K}S^\top)^{-1}S\boldsymbol{K}f_0 - f_0\|_n^2 \\
&= \|U^T\boldsymbol{K}S^\top(S\boldsymbol{K}^2S^\top + \lambda S\boldsymbol{K}S^\top)^{-1}S\boldsymbol{K}f_0 - U^T f_0\|_n^2.
\end{aligned}$$

Let $S = U_s$, where $U_s$ is the first $s$ columns of $U$. Let $f_0(\cdot) = \sum_{i=1}^n K(x_i,\cdot)w_i$ with $w = (w_1,\cdots,w_n)^\top = U\alpha$, where $\alpha \in \mathbb{R}^n$ satisfies

$$\alpha_i^2 = \begin{cases} \frac{1}{n}\frac{C}{s}\widehat{\mu}_i^{-1}, & \text{for } i = s+1,\cdots,2s \\ 0, & 0 \end{cases} \tag{7.7}$$

Then $\|f_0\|_{\mathcal{H}}^2 = n\alpha^\top D\alpha = n\sum_{i=s+1}^{2s}\alpha_i^2\widehat{\mu}_i = C$, and

$$\begin{aligned}
\|\operatorname{E}_\epsilon \widehat{f}_R - f_0\|_n^2 &= n\sum_{i=1}^{s}\alpha_i^2(\widehat{\mu}_i^2(\widehat{\mu}_i+\lambda)^{-1} - \widehat{\mu}_i)^2 + n\sum_{i=s+1}^{n}\alpha_i^2\widehat{\mu}_i^2 \\
&= \frac{1}{s}\sum_{i=s+1}^{2s}\widehat{\mu}_i \ge \widehat{\mu}_{2s} \ge \frac{1}{2}\mu_{2s} \gg \lambda^\dagger.
\end{aligned} \tag{7.8}$$

The last inequality holds with probability greater than $1 - e^{-n\delta_n}$ by Lemma S.2. On the other hand, there always exists $\lambda \gg \lambda^\dagger$, such that the corresponding $s_\lambda = s/d$. Then by (7.8), with probability greater than $1 - e^{-n\delta_n}$,

$$\| \mathrm{E}_\epsilon \widehat{f}_R - f_0\|_n^2 \leq \frac{3}{2}\mu_s \leq \frac{3}{2}\mu_{s_\lambda} \asymp \lambda,$$

i.e., $T_1 = \mathcal{O}_P(\lambda)$, based on the definition 3.1 for $s_\lambda$.

Furthermore, we have $T_2 = \mathcal{O}_P(\mu_{n,\lambda}) = \mathcal{O}_P(\frac{s_\lambda}{n})$ by the proof of Corollary 4.4. Note that $\lambda^\dagger$ satisfies $\lambda^\dagger \asymp \frac{s_\lambda^\dagger}{n}$. Then for $\lambda \gg \lambda^\dagger$, we have $T_2 = o_p(T_1)$. Therefore, $T_3 = o_P(T_1)$ due to Cauchy-Schwarz inequality $T_3 \leq T_1^{1/2} T_2^{1/2}$. Finally, with probability at least $1 - e^{-n\delta_n}$,

$$\sup_{f_0 \in \mathcal{B}} \|\widehat{f}_R - f_0\|_n^2 \gtrsim \sup_{f_0 \in \mathcal{B}} \| \mathrm{E}_\epsilon \widehat{f}_R - f_0\|_n^2 \geq C\mu_{2s} \gg C\lambda^\dagger.$$

The last step is based on the definition of $\lambda^\dagger$ and the fact that $2s \ll s^\dagger$. □

## 7.7 Proof of Theorem 4.14

*Proof.* Without loss of generality, here we consider $H_0 : f = f_0$ with $f_0 = 0$. We construct the true $f(\cdot) = \sum_{i=1}^n K(x_i, \cdot)w_i$ with $w = (w_1, \cdots, w_n)^\top = U\alpha$, where $\alpha \in \mathbb{R}^n$ satisfies

$$\alpha_i^2 = \begin{cases} \frac{1}{n}\frac{C}{s-1}\widehat{\mu}_{gs+k}^{-1} & \text{for } \ i = (gs + k) \quad k = 1, 2, \cdots, s-1 \\ 0 & 0 \end{cases} \tag{7.9}$$

Choose $g \geq 1$ to be an integer satisfying $(g + 1)s \ll s^*$. By definition,

$$\|f\|_{\mathcal{H}}^2 = n\alpha^\top D\alpha = n\sum_{k=1}^{s-1} \alpha_{gs+k}^2 \widehat{\mu}_{gs+k} = C$$

and

$$\|f\|_n^2 = n\alpha^\top D^2\alpha = n\sum_{k=1}^{s-1} \alpha_{gs+k}^2 \widehat{\mu}_{gs+k}^2 = \frac{C}{s-1}\sum_{k=1}^{s-1} \widehat{\mu}_{gs+k}.$$

Then by Lemma S.2, with probability at least $1 - e^{-n\delta_n}$,

$$\|f\|_n^2 \geq \frac{C}{2}\mu_{gs+s} = \beta_{n,\lambda}^2 d^{*2},$$

where $\beta_{n,\lambda}^2 = \frac{C}{2}\mu_{(g+1)s}/\mu_{s^*}$, and $\beta_{n,\lambda}^2 \to \infty$ as $n \to \infty$, $d^{*2} = \lambda^* \asymp \mu_{s^*}$.

Let $S = U_s$, where $U_s$ is the first $s$ columns of $U$. Then $S$ satisfies Assumption A2 with $c_1 = c_2 = +\infty$, i.e., the K-satisfiability holds almost surely. $SU = (\widetilde{S}_1\, \widetilde{S}_2)$, where $\widetilde{S}_1 = SU_s = I_{s\times s}$ and $\widetilde{S}_2 = SU_{n-s} = \mathbf{0}$. Recall in eq. (7.3) for $nT_{n,\lambda}$. Plugging $S$ and $f$ into $T_1$, we have

$$\begin{aligned} T_1 &= \boldsymbol{f}^\top \boldsymbol{K} S^\top (S\boldsymbol{K}^2 S^\top + \lambda S\boldsymbol{K}S^\top)^{-1} S\boldsymbol{K}^2 S^\top (S\boldsymbol{K}^2 S^\top + \lambda S\boldsymbol{K}S^\top)^{-1} S\boldsymbol{K}\boldsymbol{f} \\ &= \boldsymbol{f}^\top UD\widetilde{S}^\top (\widetilde{S}D^2\widetilde{S}^\top + \lambda\widetilde{S}D\widetilde{S}^\top)^{-1}\widetilde{S}D^2\widetilde{S}^\top(\widetilde{S}d\widetilde{S}^\top + \lambda\widetilde{S}D\widetilde{S}^\top)^{-1}\widetilde{S}DU^\top \boldsymbol{f} \\ &= n\sum_{i=1}^s \alpha_i^2 \frac{\widehat{\mu}_i^4}{(\widehat{\mu}_i + \lambda)^2} = 0, \end{aligned}$$

where the last step is by the construction of $\alpha$ that $\alpha_1 = \alpha_s = 0$. $T_1 = 0 \ll n\sigma_{n,\lambda}$. Furthermore, $|T_3| = T_1^{1/2} O_{P_f}(1) = o_{P_f}(n\sigma_{n,\lambda})$. Therefore

$$\begin{aligned} \frac{T_{n,\lambda} - \mu_{n,\lambda}}{\sigma_{n,\lambda}} &= \frac{T_1 + T_3}{n\sigma_{n,\lambda}} + \frac{T_2/n - \mu_{n,\lambda}}{\sigma_{n,\lambda}} \\ &= \frac{T_2/n - \mu_{n,\lambda}}{\sigma_{n,\lambda}} + o_{P_f}(\sigma_{n,\lambda}) \\ &\xrightarrow{d} N(0,1). \end{aligned}$$

Then we have, as $n \to \infty$, with probability at least $1 - e^{-cn\delta_n} - e^{-c_1 s} - e^{-c_2 s_\lambda}$,

$$\inf_{f\in\mathcal{B}, \|f\|_n \geq \beta_{n,\lambda} d^*} \mathrm{P}_f\big(\phi_{n,\lambda} = 1 | \boldsymbol{x}, S\big) \leq \mathrm{P}_f\big(\phi_{n,\lambda} = 1 | \boldsymbol{x}, S\big) \to \alpha.$$

□

**Acknowledgments.** We thank Dr. Wen-Xin Zhou for pointing out [24].

*Supplement to*

# Nonparametric Testing under Random Projection

In this document, additional proofs, simulation results and technical arguments are provided.

- Section S.1 and S.2 include some technical lemmas.
- Section S.3 includes examples to verify Assumption A1.
- Section S.4 includes the proof of Lemma 3.2.
- Section S.5 includes the proof of Lemma 3.1.
- Section S.6 includes the proof of Lemma 4.1.
- Section S.7 includes the proof of Lemma 4.3.
- Section S.8 includes the proof of Corollary 4.4.
- Section S.9 includes the proof of Theorem 5.1.
- Section S.10 includes the proof of Theorem 5.2.
- Section S.11 includes some auxiliary lemmas.
- Section S.12 includes additional simulation results using Bernoulli random projection matrices.

## S.1 A Key Lemma

We first show that $\Psi_\lambda(r)$ and $\widehat{\Psi}_\lambda(r)$ have an asymptotically equivalent expression in terms of $\mu_i$'s ($\widehat{\mu}_i$'s) for a wide-ranging $\alpha$, where $\widehat{\mu}_i$'s are the eigenvalues of $\boldsymbol{K}$ (in a descending order $\widehat{\mu}_1 \geq \widehat{\mu}_2 \geq \cdots \geq \widehat{\mu}_n \geq 0$). Recall that $\mu_i$'s are eigenvalues of the kernel function $K$; see (2.1).

**Lemma S.1.** *(a) Suppose $\mu_1 > 1/n$. For any $\lambda > 1/n$, it holds that*

$$\Psi_\lambda(r) \asymp \sqrt{\frac{1}{n}\sum_{i=1}^{\infty} \kappa_\lambda \min\{\frac{r}{\kappa_\lambda}, \mu_i\}}. \tag{S.1}$$

*(b) For any $\lambda > 0$, it holds that*

$$\widehat{\Psi}_\lambda(r) \asymp \sqrt{\frac{1}{n}\sum_{i=1}^{n} \kappa_\lambda \min\{\frac{r}{\kappa_\lambda}, \widehat{\mu}_i\}}. \tag{S.2}$$

*Proof.* We first prove (S.2). Define $\boldsymbol{x} = (x_1, \cdots, x_n)$. Let $(K, \langle\cdot,\cdot\rangle)$ be a RKHS $\mathcal{H}$. Any $f \in \mathcal{H}$ can be presented as $f(\cdot) = \sum_{i=1}^n c_i K(x_i, \cdot) + \xi(\cdot)$ with $\xi(\cdot) \perp \text{span}\{K(x_1,\cdot), \cdots, K(x_n,\cdot)\}$. Therefore $f(x_j) = \sum_{i=1}^n c_i K(x_i, x_j)$ for $j = 1, \cdots, n$. Let $\boldsymbol{f} = (f(x_1), \cdots, f(x_n))^\top$, we have $\boldsymbol{f} = n\boldsymbol{K}\boldsymbol{c}$, where $\boldsymbol{K}$ is the kernel matrix, $\boldsymbol{c} = (c_1, \cdots, c_n)^\top$. Let $\boldsymbol{K} = \Phi\Phi^\top$ with $\Phi = UD^{1/2}$. Then $n\boldsymbol{K}\boldsymbol{c} = n\Phi\Phi^\top\boldsymbol{c} = \sqrt{n}\Phi\beta$ with $\sqrt{n}\Phi^\top\boldsymbol{c} \equiv \beta$. Then we have

$$\sum_{i=1}^n f^2(x_i) = n^2\boldsymbol{c}^\top\boldsymbol{K}^2\boldsymbol{c} = n\beta^\top\Phi^\top\Phi\beta = n\beta^\top D^{1/2}U^\top UD^{1/2}\beta = n\sum_{i=1}^n \beta_i^2\widehat{\mu}_i.$$

Note that $P_nf^2 = \frac{1}{n}\sum_{i=1}^n f(x_i)^2 \le r$ is equivalent to $\sum_{i=1}^n \beta_i^2\widehat{\mu}_i \le r$. Therefore,

$$\sup_{\substack{f\in\mathcal{F}_\lambda \\ P_nf^2\le r}} |\sum_{i=1}^n \sigma_i f(x_i)|^2 = \sup_{\substack{f\in\mathcal{F}_\lambda \\ P_nf^2\le r}} n^2|\sigma^\top\boldsymbol{K}\boldsymbol{c}|^2 = \sup_{\substack{\beta\in\mathbb{R}^n, \sum_{i=1}^n\beta_i^2\le\kappa_\lambda \\ \sum_{i=1}^n\beta_i^2\widehat{\mu}_i\le r}} n|\sigma^\top\Phi\beta|^2$$

Define $\widehat{F}_\lambda = \{\beta \in \mathbb{R}^n | \sum_{i=1}^n \beta_i^2 \le \kappa_\lambda, \sum_{i=1}^n \beta_i^2\widehat{\mu}_i \le r\}$, and $\widetilde{F}_\lambda = \{\beta \in \mathbb{R}^n | \sum_{i=1}^n \widehat{d}_i\beta_i^2 \le 1\}$, where $\widehat{d}_i = (\kappa_\lambda \min\{1, r/(\kappa_\lambda\widehat{\mu}_i)\})^{-1}$. So $\widetilde{F}_\lambda \subseteq \widehat{F}_\lambda \subseteq \sqrt{2}\widetilde{F}_\lambda$. Let $\Lambda = \text{diag}\{\widehat{d}_1, \cdots, \widehat{d}_n\}$. Then we have

$$\begin{aligned}
&\mathrm{E}\left(\sup_{\substack{\beta\in\mathbb{R}^n, \sum_{i=1}^n\beta_i^2\le\kappa_\lambda \\ \sum_{i=1}^n\beta_i^2\widehat{\mu}_i\le r}} |\sigma^\top\Phi\beta|^2 \middle| \boldsymbol{x}\right) \asymp \mathrm{E}\left(\sup_{\beta\in\widetilde{F}_\lambda} |\sigma^\top\Phi\Lambda^{-1/2}\Lambda^{1/2}\beta|^2 \middle| \boldsymbol{x}\right) \\
=&\mathrm{E}\left(\sup_{\substack{d\in\mathbb{R}^n \\ d'd\le 1}} |\sigma^\top\Phi\Lambda^{-1/2}d|^2 \middle| \boldsymbol{x}\right) = \mathrm{E}\left(\|\sigma^\top\Phi\Lambda^{-1/2}\|_2^2 \middle| \boldsymbol{x}\right) = \mathrm{E}\left(\sigma^\top\Phi\Lambda^{-1}\Phi^\top\sigma \middle| \boldsymbol{x}\right)
\end{aligned}$$

Note that

$$\mathrm{E}\left(\sigma^\top\Phi\Lambda^{-1}\Phi^\top\sigma \middle| \boldsymbol{x}\right) = \mathrm{tr}\left(\Phi\Lambda^{-1}\Phi^\top\right) = \mathrm{tr}\left(\Lambda^{-1}\text{diag}\{\widehat{\mu}_1, \cdots, \widehat{\mu}_n\}\right) = \sum_{i=1}^n \frac{\widehat{\mu}_i}{\widehat{d}_i}$$

Therefore, by Kahane-Khintchine inequality, we have

$$\widehat{\Psi}_\lambda(r) \asymp \sqrt{\frac{1}{n}\sum_{i=1}^n \widehat{\mu}_i/\widehat{d}_i} = \sqrt{\frac{1}{n}\sum_{i=1}^n \kappa_\lambda \min\{\widehat{\mu}_i, \frac{r}{\kappa_\lambda}\}}.$$

Similarily, we can achieve (S.1). □

## S.2 Properties of eigenvalues

**Lemma S.2.** *(a) Suppose that $K$ has eigenvalues satisfying $\mu_i \asymp i^{-2m}$ with $m > 3/2$. Then for $i = 1, \cdots, n^{1/(2m)}$,*

$$P\left(|\widehat{\mu}_i - \mu_i| \le \frac{1}{2}\mu_i\right) \ge 1 - e^{-c_m n i^{-4m/(2m-1)}}.$$

*where $c_m$ is an universal constant depending only on $m$.*

(*b*) *Suppose that* $K$ *has eigenvalues satisfying* $\mu_i \asymp \exp(-\gamma i^p)$ *with* $\gamma > 0$, $p \geq 1$. *Then for* $i = o(n^{1/2})$,

$$P\big(|\widehat{\mu}_i - \mu_i| \leq \frac{1}{2}\mu_i\big) \geq 1 - e^{-c_{\gamma,p} n i^{-2}},$$

*where* $c_{\gamma,p}$ *is an universal constant depending only on* $\gamma$ *and* $p$.

*For* $i = O(n^{1/2})$, *we have*

$$P\left(|\widehat{\mu}_i - \mu_i| \leq i\mu_i\right) \geq 1 - e^{-c'_{\gamma,p} n},$$

*where* $c'_{\gamma,p}$ *is an universal constant depending only on* $\gamma$ *and* $p$.

*Proof.* We apply the proof of Theorem 3 in [8] to deduce our results. Recall in Theorem 3 of [8], for $1 \leq i \leq n$, $1 \leq r \leq n$,

$$|\widehat{\mu}_i - \mu_i| \leq \mu_i \|C_n^r\| + \mu_r + \Lambda_{>r}, \tag{S.3}$$

where $\Lambda_{>r} = \sum_{i=r+1}^{\infty} \mu_i$, and $\|C_n^r\|$ satisfies

$$\mathrm{P}\Big(\|C_n^r\| \geq \widetilde{\varepsilon}\Big) \leq r(r+1) e^{-\frac{n\widetilde{\varepsilon}^2}{2M^4 r^2}}, \tag{S.4}$$

based on Lemma 7 in [8]. $M$ is an absolute constant here.

First we prove Lemma S.2 (a). Consider the polynomial decaying kernel with $\mu_i \asymp i^{-2m}$. Notice that

$$\Lambda_{>r} \asymp \sum_{i=r+1}^{\infty} i^{-2m} \leq \int_r^{\infty} x^{-2m} dx = \frac{r^{1-2m}}{2m-1}$$

Consider $i = 1, \cdots, n^{\frac{1}{2m}}$, let

$$\frac{r^{1-2m}}{2m-1} = \frac{1}{4}\mu_i,$$

then $r = a_m i^{2m/(2m-1)}$, where $a_m$ is a constant only depends on $m$. Let $\widetilde{\varepsilon} = \frac{1}{4}$. Plugging $r$ into (S.4), we have

$$\mathrm{P}\Big(\|C_n^r\| \geq \frac{1}{4}\Big) \leq e^{-c_m n i^{-4m/(2m-1)}},$$

where $c_m = (64M^4 a_m^2)^{-1}$ is an universal constant depends on $m$. Then we obtain that

$$\mathrm{P}\Big(|\widehat{\mu}_i - \mu_i| \leq \frac{1}{2}\mu_i\Big) \geq 1 - e^{-c_m n i^{-4m/(2m-1)}}.$$

Next, we prove Lemma S.2 (b). Consider the exponential decaying kernel with $\mu_i \asymp e^{-\gamma i^p}$. For $1 \leq r \leq n$, when $p = 1$, then

$$\Lambda_{>r} \asymp \sum_{i=r+1}^{\infty} e^{-\gamma i} \leq \int_r^{\infty} e^{-\gamma x} dx = \frac{e^{-\gamma r}}{\gamma};$$

when $p \geq 2$, using integration by parts, we have

$$\begin{aligned}\Lambda_{>r} \asymp \sum_{i=r+1}^{\infty} e^{-\gamma i^p} &\leq \int_r^{\infty} e^{-\gamma x^p} dx \\ &= \frac{1}{\gamma p r^{p-1}} e^{-\gamma r^p} - \int_r^{\infty} \frac{p-1}{\gamma p x^p} e^{-\gamma x^p} dx \leq a_{\gamma,p} e^{-\gamma r^p}.\end{aligned}$$

For $i = o(n^{1/2})$, let $\mu_r + \Lambda_{>r} \leq (1 + a_{\gamma,p}) e^{-\gamma r^p} = \frac{1}{4}\mu_i$, we have $r = b_{\gamma,p} i$, where $b_{\gamma,p}$ is a constant only depends on $\gamma, p$. Then plugging $\widetilde{\epsilon} = \frac{1}{4}$ and $r$ into (S.4), we have

$$\mathrm{P}\Big(\|C_n^r\| \geq \frac{1}{4}\Big) \leq e^{-c_{\gamma,p} n i^{-2}}.$$

where $c_{\gamma,p} = (64 M^4 b_{\gamma,p}^2)^{-1}$ is an absolute constant only depends on $\gamma, p$. Finally, by (S.3), we have

$$\mathrm{P}\Big(|\widehat{\mu}_i - \mu_i| \leq \frac{1}{2}\mu_i\Big) \geq 1 - e^{-c_{\gamma,p} n i^{-2}}.$$

When $i \geq n^{1/2}$, we do not need a very tight bound. Let $\widetilde{\varepsilon} = i$, $r = i$, then we have

$$\mathrm{P}\Big(|\widehat{\mu}_i - \mu_i| \leq i\mu_i\Big) \geq 1 - e^{-c'_{\gamma,p} n},$$

where $c'_{\gamma,p}$ is an absolute constant only depends on $\gamma, p$. □

## S.3 Verification of Assumption A1

Let us verify Assumption A1 in PDK and EDK.

First consider PDK with $\mu_i \asymp i^{-2m}$ for a constant $m > 1/2$ which includes kernels of Sobolev space and Besov Space. An $m$-th order Sobolev space, denoted $\mathcal{H}^m([0,1])$, is defined as

$$\begin{aligned}\mathcal{H}^m([0,1]) = &\{f : [0,1] \to \mathbb{R} | f^{(j)} \text{ is abs. cont for } j = 0, 1, \cdots, m-1, \\ &\text{and } f^m \in L_2([0,1])\}.\end{aligned}$$

An $m$-order periodic Sobolev space, denoted $H_0^m(\mathbb{I})$, is a proper subspace of $\mathcal{H}^m([0,1])$ whose element fulfills an additional constraint $g^{(j)}(0) = g^{(j)}(1)$ for $j = 0, \ldots, m-1$. The basis functions $\phi_i$'s of $H_0^m(\mathbb{I})$ are

$$\phi_i(z) = \begin{cases} \sigma, & i = 0, \\ \sqrt{2}\sigma \cos(2\pi k z), & i = 2k, k = 1, 2, \ldots, \\ \sqrt{2}\sigma \sin(2\pi k z), & i = 2k-1, k = 1, 2, \ldots. \end{cases}$$

The corresponding eigenvalues are $\mu_{2k} = \mu_{2k-1} = \sigma^2 (2\pi k)^{-2m}$ for $k \geq 1$ and $\mu_0 = \infty$. In this case, $\sup_{i \geq 1} \|\phi\|_{\sup} < \infty$. For any $k \geq 1$,

$$\sum_{i=k+1}^{\infty} \mu_i \lesssim \int_k^{\infty} x^{-2m} dx = \frac{k^{1-2m}}{2m-1} \lesssim \frac{k\mu_k}{2m-1}.$$

Therefore, there exists a constant $C < \infty$, such that

$$\sup_{k\geq 1} \frac{\sum_{i=k+1}^{\infty} \mu_i}{k\mu_k} = C < \infty.$$

Hence, Assumption A1 holds true. Verification of Assumption A1 on the eigenfunctions for Sobolev space kernel can be found in [34].

Next, let us consider EDK with $\mu_i \asymp \exp(-\gamma i^p)$ for constants $\gamma > 0$ and $p > 0$. Gaussian kernel $K(x, x') = \exp\big(-(x - x')^2/\sigma^2\big)$ is an EDK of order $p = 2$, with eigenvalues $\mu_i \asymp \exp(-\pi i^2)$ as $i \to \infty$, and the corresponding eigenfunctions

$$\phi_i(x) = (\sqrt{5}/4)^{1/4}(2^{i-1} i!)^{-1/2} e^{-(\sqrt{5}-1)x^2/4} H_i((\sqrt{5}/2)^{1/2} x),$$

where $H_i(\cdot)$ is the $i$-th Hermite polynomial; see [38] for more details. Then $\sup_{i\geq 1} \|\phi_i\|_{\sup} < \infty$ trivially holds. For any $k \geq 1$,

$$\sum_{i=k+1}^{\infty} \mu_i \lesssim \int_k^{\infty} e^{-\gamma x^p} dx = \frac{1}{\gamma p k^{p-1}} e^{-\gamma k^p} - \int_k^{\infty} \frac{p-1}{\gamma p x^p} e^{-\gamma x^p} dx \leq \frac{1}{\gamma p k^{p-1}} e^{-\gamma k^p}.$$

Therefore,

$$\sup_{k\geq 1} \frac{\sum_{i=k+1}^{\infty} \mu_i}{k\mu_k} < \infty.$$

Hence, Assumption A1 holds.

## S.4 Proof of Lemma 3.2

*Proof.* We first observe that $\mathcal{F}_\lambda = \mathrm{star}(\mathcal{F}_\lambda, 0)$, the star-hull of $\mathcal{F}_\lambda$ at zero. Note that the supremum in the definitions of $\Psi_\lambda$ and $\widehat{\Psi}_\lambda$ is based on "quadratic type" constraints $Pf^2 \leq r$ and $P_n f^2 \leq r$. Then following [6], $\Psi_\lambda$ and $\widehat{\Psi}_\lambda$ are both sub-root functions, and thus have unique nonzero fixed points. Also refer to [23] for the definitions of star-hull and sub-root functions. Define

$$\widehat{\Psi}'_\lambda(r) = \widehat{\Psi}_\lambda(r) + \frac{c_1\delta}{n} = \mathrm{E}\Big\{ \sup_{\substack{f\in\mathcal{F}_\lambda \\ P_n f^2 \leq r}} \frac{1}{n}\sum_{i=1}^{n} \sigma_i f(x_i) \Big| \boldsymbol{x} \Big\} + \frac{c_1\delta}{n},$$

where $c_1$ is a constant. Then by Theorem 4.2 in [4], the fixed points $r_\lambda$ and $\widehat{r}_\lambda$ of $\Psi_\lambda(r)$ and $\widehat{\Psi}'_\lambda(r)$ satisfy: $r_\lambda \asymp \widehat{r}_\lambda$ with probability at least $1 - 4e^{-\delta}$, provided that $r_\lambda \geq c_1\delta/n$.

Let $r = r_\lambda$ in Lemma S.1, we have

$$r_\lambda \asymp \sqrt{\frac{1}{n}\sum_{i=1}^{\infty} \kappa_\lambda \min\Big\{\frac{r_\lambda}{\kappa_\lambda}, \mu_i\Big\}}. \tag{S.5}$$

Define $\eta_\lambda = \mathrm{argmin}\{i : \mu_i \le r_\lambda/\kappa_\lambda\} - 1$, then (S.5) implies

$$\frac{nr_\lambda^2}{\kappa_\lambda} \asymp \eta_\lambda \frac{r_\lambda}{\kappa_\lambda} + \sum_{i=\eta_\lambda+1}^{\infty} \mu_i. \tag{S.6}$$

Note that

$$\frac{\sum_{i=k+1}^{\infty} \mu_i}{k\mu_{k+1}} = \frac{1}{k} + \frac{k+1}{k} \cdot \frac{\sum_{i=k+2}^{\infty} \mu_i}{(k+1)\mu_{k+1}} \le 1 + 2C,$$

where $C = \sup_{k\ge 1} \frac{\sum_{i=k+1}^{\infty} \mu_i}{k\mu_k} < \infty$ by Assumption A1. Therefore, $\sum_{i=\eta_\lambda+1}^{\infty} \mu_i \lesssim \eta_\lambda \mu_{\eta_\lambda+1} \le \eta_\lambda \frac{r_\lambda}{\kappa_\lambda}$. Then by (S.6), we have $\frac{nr_\lambda^2}{\kappa_\lambda} \asymp \eta_\lambda \frac{r_\lambda}{\kappa_\lambda}$, i.e., $r_\lambda \asymp \frac{\eta_\lambda}{n}$. Note

$$\mu_{\eta_\lambda} > \frac{r_\lambda}{\kappa_\lambda} \asymp \frac{\eta_\lambda}{n\kappa_\lambda} = \frac{\lambda \eta_\lambda}{s_\lambda},$$

and recall $s_\lambda = \mathrm{argmin}\{i : \mu_i \le \lambda\} - 1$ which implies $\mu_{s_\lambda+1} \le \lambda < \mu_{s_\lambda}$. Then,

$$\frac{\mu_{\eta_\lambda}}{\eta_\lambda} \gtrsim \frac{\lambda}{s_\lambda} \ge \frac{\mu_{s_\lambda+1}}{s_\lambda} \ge \frac{\mu_{s_\lambda+1}}{s_\lambda+1}. \tag{S.7}$$

Note that $\mu_k/k$ is a decreasing function of $k$, we thus have $\eta_\lambda \lesssim s_\lambda + 1$ by (S.7), i.e., $\eta_\lambda \lesssim s_\lambda$. On the other hand,

$$\mu_{\eta_\lambda+1} \le \frac{r_\lambda}{\kappa_\lambda} \asymp \frac{\lambda\eta_\lambda}{s_\lambda} \lesssim \frac{\mu_{s_\lambda}(\eta_\lambda+1)}{s_\lambda},$$

i.e., $\frac{\mu_{\eta_\lambda+1}}{\eta_\lambda+1} < \frac{\mu_{s_\lambda}}{s_\lambda}$, and we have $\eta_\lambda + 1 \gtrsim s_\lambda$, i.e., $\eta_\lambda \gtrsim s_\lambda$. Therefore, $\eta_\lambda \asymp s_\lambda$. Then, we achieve that $r_\lambda \asymp \frac{s_\lambda}{n}$. Suppose there exists an constant $c_2$, such that $r_\lambda \ge c_2 \frac{s_\lambda}{n}$, let $\delta = c_2 s_\lambda / c_1$, then with probability greater than $1 - e^{-cs_\lambda}$, $\widehat{r}_\lambda \asymp r_\lambda \asymp s_\lambda/n$, where $c = \frac{c_2}{2c_1}$. □

## S.5 Proof of Lemma 3.1

*Proof.* Plugging the fixed point $r_\lambda$ and $\widehat{r}_\lambda$ into (S.1) and (S.2) in Lemma S.1, we have

$$r_\lambda \asymp \sqrt{\frac{1}{n}\sum_{i=1}^{\infty} \kappa_\lambda \min\left\{\frac{r_\lambda}{\kappa_\lambda}, \mu_i\right\}}, \tag{S.8}$$

$$\widehat{r}_\lambda \asymp \sqrt{\frac{1}{n}\sum_{i=1}^{n} \kappa_\lambda \min\left\{\frac{\widehat{r}_\lambda}{\kappa_\lambda}, \widehat{\mu}_i\right\}} + \frac{c_1\delta}{n}. \tag{S.9}$$

By Lemma 3.2, $r_\lambda \asymp s_\lambda/n$, and $r_\lambda/\kappa_\lambda \asymp \lambda$; for the empirical version, let $\delta = s_\lambda$, then with probability at least $1 - 4e^{-s_\lambda}$, $\widehat{r}_\lambda \asymp s_\lambda/n$ leads to $\widehat{r}_\lambda/\kappa_\lambda \asymp \lambda$. Recall that $\widehat{s}_\lambda = \mathrm{argmin}\{i : \widehat{\mu}_i \le \lambda\} - 1$. Then by (S.9), with probability at least $1 - 4e^{-s_\lambda}$,

$$\frac{1}{n}\sum_{i=\widehat{s}_\lambda+1}^{n} \widehat{\mu}_i \le \frac{1}{n}\sum_{i=1}^{n} \min\left\{\lambda, \widehat{\mu}_i\right\} \asymp \frac{1}{n}\sum_{i=1}^{n} \min\left\{\frac{\widehat{r}_\lambda}{\kappa_\lambda}, \widehat{\mu}_i\right\} \lesssim \widehat{r}_\lambda^2/\kappa_\lambda \asymp \lambda s_\lambda/n,$$

where the last step is by $\kappa_\lambda = \frac{s_\lambda}{n\lambda}$, and $\widehat{r}_\lambda \asymp s_\lambda/n$. Therefore,

$$\sum_{i=\widehat{s}_\lambda+1}^{n} \widehat{\mu}_i \lesssim \lambda s_\lambda \leq s_\lambda \mu_{s_\lambda},$$

based on the definition (3.1) that $\lambda < \mu_{s_\lambda}$. □

## S.6 Proof of Lemma 4.1

In this section, we first prove the following (S.10) and (S.11):

$$\mathrm{P}\Big(\frac{1}{\sqrt{2}} \leq \lambda_{\min}(SU_1) \leq \lambda_{\max}(SU_1) \leq \sqrt{\frac{3}{2}}\Big|\boldsymbol{x}\Big) \geq 1 - \exp\big(-c_1 s\big) \tag{S.10}$$

almost surely, where $c_1 > 0$ is an absolute constant independent of $n, s$; $\lambda_{\min}(SU_1)$ ($\lambda_{\max}(SU_1)$) is the smallest (largest) singular value of $SU_1$.

$$\mathrm{P}\Big(\mathrm{P}\big(\|SU_2D_2^{1/2}\|_{\mathrm{op}} \leq c\lambda\big|\boldsymbol{x}\big) \geq 1 - e^{-c_1' s}\Big) \geq 1 - e^{-c_2 s_\lambda}, \tag{S.11}$$

where $c$ and $c_1'$ are constants independent of $n, s$ and $c_2 = 1/2$. The result of Lemma 4.1 directly follow from (S.10) and (S.11).

*Proof.* For $\boldsymbol{K} = UDU^\top$, let $U = (U_1, U_2)$ with $U_1 \in \mathbb{R}^{n\times\widehat{s}_\lambda}$, and $U_2 \in \mathbb{R}^{n\times(n-\widehat{s}_\lambda)}$. Recall $S = \frac{1}{\sqrt{s}}S^*$, where $S^*$ is the random matrix with independent centered sub-Gaussian entries (with variance as one), then each row $S_i^*$ is independent sub-Gaussian isotropic random vectors in $\mathbb{R}^n$, i.e., $\mathrm{E}\, S_i^* S_i^{*\top} = I_{n\times n}$, $i = 1, \cdots, s$. Let $SU_1 = \frac{1}{\sqrt{s}}(\eta_1, \cdots, \eta_s)^\top$, where $\eta_i \in \mathbb{R}^{\widehat{s}_\lambda\times 1}$ with each entry $\eta_{ij} = S_i^{*\top}U_{1(j)}$, $U_{1(j)}$ is the $j$th column of $U_1$, $j = 1, \cdots, \widehat{s}_\lambda$.

Firstly, conditional on $\boldsymbol{x}$, by the definition of sub-Gaussian random vector, each entry $\eta_{ij}$ is sub-Gaussian, $\eta_i$ and $\eta_j (i \neq j)$ are independent, and $\eta_i$ is isotropic sub-Gaussian random vector due to the fact that $\mathrm{E}(\eta_i\eta_i^\top|\boldsymbol{x}) = U_1^\top(\mathrm{E}\, S_i^*S_i^{*\top})U_1 = I_{\widehat{s}_\lambda\times\widehat{s}_\lambda}$. By Theorem 5.39 in [40], for any $t > 0$,

$$\begin{aligned}&\mathrm{P}\Big(\frac{\sqrt{s} - C\sqrt{s_\lambda} - t}{\sqrt{s}} \leq \lambda_{\min}(SU_1) \leq \lambda_{\max}(SU_1) \leq \frac{\sqrt{s} + C\sqrt{s_\lambda} + t}{\sqrt{s}}|\boldsymbol{x}\Big)\\ \geq& 1 - 2e^{-ct^2}\end{aligned}$$

almost surely. Let $t = \frac{\sqrt{s}}{5}$, and $d \geq (0.02C)^2$, we have

$$\mathrm{P}\Big(\frac{1}{\sqrt{2}} \leq \lambda_{\min}(SU_1) \leq \lambda_{\max}(SU_1) \leq \sqrt{\frac{3}{2}}|\boldsymbol{x}\Big) \geq 1 - 2e^{-cs/25} \tag{S.12}$$

almost surely. Here $C, c > 0$ only depend on the sub-Gaussian norm $L := \max_i \|\eta_i\|_{\psi_2}$ conditional on $\boldsymbol{x}$. Note that $\eta_i = S_i^{*\top}U_i$,

$$\|U_i^\top S_i^*\|_{\psi_2} = \sup_{\nu\in\mathcal{S}^{s_\lambda-1}} \|\langle U_1^\top S_i^*, \nu\rangle\|_{\psi_2} = \sup_{\kappa\in\mathcal{S}^{n-1}} \|\kappa^\top S_i^*\|_{\psi_2} \quad \text{and}$$

$$\sup_{\kappa\in\mathcal{S}^{n-1}}\|\kappa^\top S_i^*\|_{\psi_2}^2=\|\sum_{j=1}^n \kappa_j S_{ij}^*\|_{\psi_2}^2\le C\sum_{j=1}^n \kappa_j^2\|S_{ij}^*\|_{\psi_2}^2\le C\max_{1\le j\le n}\|S_{ij}^*\|_{\psi_2}^2$$

Therefore, $L\le \max_{i,j}\|S_{ij}^*\|_{\psi_2}$, which is bounded. Lastly, we have

$$\mathrm{P}\Big(\mathrm{P}\Big(\frac{1}{\sqrt{2}}\le\lambda_{\min}(SU_1)\le\lambda_{\max}(SU_1)\le\sqrt{\frac{3}{2}}\Big|\boldsymbol{x}\Big)\ge 1-2e^{-cs/25}\Big)=1.$$

Set $\widetilde{c}_1=c/32$. Then (S.10) has been proved.

Next, we prove (S.11). Define $\mathcal{A}=\{\boldsymbol{x}:\boldsymbol{x}\,\text{satisfies}\,\sum_{i=\widehat{s}_\lambda+1}^n\widehat{\mu}_i\le Cs_\lambda\mu_{s_\lambda}\}$. Then $\mathrm{P}(\boldsymbol{x}\in\mathcal{A})\ge 1-4e^{-s_\lambda}$ by Lemma 3.1.

Since $(SU_2D_2^{1/2})^\top(SU_2D_2^{1/2})$ has the same non-zero eigenvalues as $SU_2D_2U_2^\top S^\top$, it is equivalent to prove $\lambda_{\max}(SU_2D_2U_2^\top S^\top)\lesssim\lambda$, where $\lambda_{\max}(\cdot)$ refers to the maximum singular value. For every $\nu\in\mathcal{S}^{s-1}$, $\nu=\kappa+w$, where $\kappa$ belongs to the 1/2-net $\mathcal{N}=\{\mu_1,\cdots,\mu_M\}$ of the set $\mathcal{S}^{s-1}$, here $M\le e^{2s}$; and $\|w\|\le 1/2$, where $\|\cdot\|$ is the Euclidean norm. Then

$$\begin{aligned}\|SU_2D_2U_2^\top S^\top\|_{\mathrm{op}}&=\sup_{\|\nu\|=1,\nu\in\mathcal{S}^{s-1}}\|SU_2D_2U_2^\top S^\top\nu\|\\&\le\sup_{\kappa\in\mathcal{S}^{s-1}}\|SU_2D_2U_2^\top S^\top\kappa\|+\sup_{w\in\mathcal{S}^{s-1}}\|SU_2D_2U_2^\top S^\top w\|\\&\le\max_{\kappa\in\mathcal{S}^{s-1}}\|SU_2D_2U_2^\top S^\top\kappa\|+\frac{1}{2}\|SU_2D_2U_2^\top S^\top\|_{\mathrm{op}},\end{aligned}$$

therefore

$$\begin{aligned}\|SU_2D_2U_2^\top S^\top\|_{\mathrm{op}}&\le 2\max_{\kappa\in\mathcal{S}^{s-1}}\|SU_2D_2U_2^\top S^\top\kappa\|\\&=2\max_{\kappa\in\mathcal{S}^{s-1}}|\langle SU_2D_2U_2^\top S^\top\kappa,\kappa\rangle|\\&=2\max_{\kappa\in\mathcal{S}^{s-1}}|\kappa' SU_2D_2U_2^\top S^\top\kappa|\\&=\frac{2}{s}\max_{\kappa\in\mathcal{S}^{s-1}}|\kappa' S^*U_2D_2U_2^\top S^{*\top}\kappa|,\end{aligned}$$

where the last equality is by the definition that $S=S^*/\sqrt{s}$.

Note that $\eta=S^{*\top}\kappa\in\mathbb{R}^n$ is a sub-Gaussian vector, and $\eta_i=\sum_{j=1}^s\kappa_jS_{ji}^*$ is independent with $\eta_j$ for $i,j\in\{1,\cdots,n\},i\ne j$; also, $\mathrm{E}(\eta_i)=0$,

$$\mathrm{Var}\left(\sum_{j=1}^s\kappa_jS_{ji}^*\right)=\sum_{j=1}^s\kappa_j^2\,\mathrm{Var}(S_{ji}^*)\le C,$$

where the last inequality is due to the fact that $\sum_{j=1}^s\kappa_j^2=1$. Let $Q=\frac{1}{s}U_2D_2U_2^\top$. By Hanson-Wright inequality (stated in Lemma S.4), we have

$$\mathrm{P}(\big|\eta'Q\eta-\mathrm{tr}(Q)\big|\ge t\big|\boldsymbol{x}\in\mathcal{A})\le 2e^{-c\min\{\frac{t^2}{K^4\|Q\|_{\mathrm{F}}^2},\frac{t}{K^2\|Q\|_{\mathrm{op}}}\}}. \tag{S.13}$$

Conditional on $\boldsymbol{x} \in \mathcal{A}$, $\operatorname{tr}(Q) = \frac{1}{s}\operatorname{tr}(D_2 U_2^\top U_2) = \frac{1}{s}\sum_{i=\widehat{s}_\lambda+1}^n \widehat{\mu}_i \le \frac{s_\lambda \lambda}{s} \le L\lambda$, where $L$ is some absolute constant by the assumption that $s \ge ds_\lambda$. The penultimate inequality is based on Lemma 3.1 and the definition of $s_\lambda$ in eq. (3.1). Also, note that

$$\|Q\|_{\mathrm{F}}^2 = \frac{1}{s^2}\operatorname{tr}(D_2^2) \le \frac{1}{s^2}\widehat{\mu}_{s_\lambda+1}\sum_{i=\widehat{s}_\lambda+1}^n \widehat{\mu}_i \le \frac{s_\lambda}{s^2}\widehat{\mu}_{s_\lambda+1}\lambda$$

the last step is based on Lemma 3.1 for $\boldsymbol{x} \in \mathcal{A}$. Let $t = L\lambda/2$, then

$$\frac{t^2}{\|Q\|_{\mathrm{F}}^2} \ge \frac{\lambda s}{\widehat{\mu}_{s_\lambda+1}} = \frac{\lambda}{\|Q\|_{\mathrm{op}}}.$$

Therefore, (S.13) can be further stated as

$$\mathrm{P}(|\eta' Q\eta - \operatorname{tr}(Q)| \ge L\lambda/2 | \boldsymbol{x} \in \mathcal{A}) \le 2e^{-cL\lambda s/(2\widehat{\mu}_{s_\lambda+1})} \le 2e^{-c'Ls},$$

where the last inequality is by the definition of $\widehat{\mu}_{s_\lambda+1}$. Finally, taking union bound over all $\mu \in \mathcal{N}$, we have

$$\begin{aligned}
&\mathrm{P}\Big(\mathrm{P}\big(\|SU_2D_2U_2^\top S^\top\|_{\mathrm{op}} \le 3L\lambda | \boldsymbol{x}\big) \ge 1 - e^{-(c'L-2)s}\Big)\\
\ge&\mathrm{P}\Big(\mathrm{P}\big(\|SU_2D_2U_2^\top S^\top\|_{\mathrm{op}} \le 3L\lambda | \boldsymbol{x} \in \mathcal{A}\big) \ge 1 - e^{-(c'L-2)s}\Big) \cdot \mathrm{P}\Big(\boldsymbol{x} \in \mathcal{A}\Big)\\
\ge&1 - 4e^{-s_\lambda}.
\end{aligned}$$

Let $c = 3L$, $c_1' = c'L - 2 > 0$ and $c_2 = 1/2$. Then, we have proved (S.11). Finally, taking $c_1 = \min\{\widetilde{c}_1, c_1'\}$, where $\widetilde{c}_1$ refers to (S.10) and $c_1'$ refers to (S.11), Lemma 4.1 have been proved. □

## S.7 Proof of Lemma 4.3

*Proof.* Suppose the true function is $f_0$, then $y_i = f_0(x_i) + \epsilon_i$ for $i = 1, \cdots, n$. Notice that $\mathrm{E}_\epsilon \widehat{f}_R(\cdot) = \sum_{i=1}^n (S\beta^\dagger)_i K(\cdot, x_i)$ with $\beta^\dagger = \frac{1}{n}(S\boldsymbol{K}^2S^\top + \lambda S\boldsymbol{K}S^\top)^{-1}S\boldsymbol{K}\boldsymbol{f}_0$, where $\boldsymbol{f}_0 = (f_0(x_1), \cdots, f_0(x_n))^\top$. It is in fact the solution of a noiseless version of quadratic program:

$$\beta^\dagger = \operatorname*{argmin}_{\beta \in \mathbb{R}^s}\left\{\frac{1}{n}\|\boldsymbol{f}_0 - n\boldsymbol{K}S^\top\beta\|_2^2 + n\lambda\beta^\top S\boldsymbol{K}S^\top\beta\right\}. \tag{S.14}$$

To prove $\|\mathrm{E}_\epsilon \widehat{f}_R - f_0\|_n^2 \le C\lambda$ with probability approaching 1, we only need to find an $\widetilde{\beta}$, such that

$$\frac{1}{n}\|\boldsymbol{f}_0 - n\boldsymbol{K}S^\top\widetilde{\beta}\|_2^2 + n\lambda\widetilde{\beta}^\top S\boldsymbol{K}S^\top\widetilde{\beta} \le C\lambda$$

with probability at least $1 - e^{-c_1 s} - e^{-c_2 s_\lambda}$, where $c_1, c_2$ are defined in Assumption A2. Note that $\boldsymbol{K} = UDU^\top$. Setting $\boldsymbol{z} = \frac{1}{\sqrt{n}}U^\top\boldsymbol{f}_0$, (S.14) is equivalent to (S.15) as follows:

$$\beta^\dagger = \operatorname*{argmin}_{\beta \in \mathbb{R}^s}\left\{\|\boldsymbol{z} - \sqrt{n}D\widetilde{S}^\top\beta\|_2^2 + n\lambda\beta^\top\widetilde{S}^\top D\widetilde{S}\beta\right\}. \tag{S.15}$$

Let $\boldsymbol{z} = (z_1, z_2)^\top$, where $z_1 \in \mathbb{R}^{\widehat{s}_\lambda}$, and $z_2 \in \mathbb{R}^{n-\widehat{s}_\lambda}$. Correspondingly, divide $D$ into $D_1, D_2$, where $D_1 = \operatorname{diag}\{\widehat{\mu}_1, \cdots, \widehat{\mu}_{\widehat{s}_\lambda}\}$ and $D_2 = \operatorname{diag}\{\widehat{\mu}_{\widehat{s}_\lambda+1}, \cdots, \widehat{\mu}_n\}$. Denote $\widetilde{S} = SU = (\widetilde{S}_1, \widetilde{S}_2)$ with $\widetilde{S}_1 \in \mathbb{R}^{s\times s_\lambda}$ as the left block and $\widetilde{S}_2 \in \mathbb{R}^{s\times(n-s_\lambda)}$ as the right block. We construct an $\widetilde{\beta}$ as

$$\widetilde{\beta} = \frac{1}{\sqrt{n}}\widetilde{S}_1(\widetilde{S}_1^\top \widetilde{S}_1)^{-1} D_1^{-1} z_1 \in \mathbb{R}^s.$$

Plugging $\widetilde{\beta}$ into (S.15), we see that

$$\begin{aligned}\|\boldsymbol{z} - \sqrt{n}D\widetilde{S}^\top \widetilde{\beta}\|_2^2 &= \|z_1 - \sqrt{n}D_1\widetilde{S}_1^\top \widetilde{\beta}\|_2^2 + \|z_2 - D_2\widetilde{S}_2^\top \widetilde{S}_1(\widetilde{S}_1^\top \widetilde{S}_1)^{-1}D_1^{-1}z_1\|_2^2 \\ &= T_1^2 + T_2^2.\end{aligned}$$

It is obvious that $T_1^2 = 0$, and next we analyze $T_2$.

Note that for any $f_0(\cdot) \in \mathcal{H}$, there exists an $n\times 1$ vector $\omega$, such that $f_0(\cdot) = \sum_{i=1}^n K(\cdot, x_i)\omega_i + \xi(\cdot)$, where $\xi(\cdot) \in \mathcal{H}$ is orthogonal to the span of $\{K(\cdot, x_i), i = 1, \cdots, n\}$. Then, $\xi(x_j) = \langle \xi, K(\cdot, x_j)\rangle = 0$, and $f_0(x_j) = \sum_{i=1}^n K(x_i, x_j)\omega_i$. Therefore $\boldsymbol{f}_0 = n\boldsymbol{K}\omega$, where $\boldsymbol{K}$ is the empirical kernel matrix. Suppose $\|f_0\|_{\mathcal{H}} \le 1$, then

$$\begin{aligned}&n\omega^\top \boldsymbol{K}\omega \le \|f_0\|_{\mathcal{H}} \le 1 \Rightarrow n\omega^\top \boldsymbol{K}\boldsymbol{K}^{-1}\boldsymbol{K}\omega^\top \le 1 \Rightarrow \frac{1}{n}\boldsymbol{f}_0{}^\top \boldsymbol{K}^{-1}\boldsymbol{f}_0 \le 1 \\ \Rightarrow& \frac{1}{n}\boldsymbol{f}_0{}^\top UD^{-1}U^\top \boldsymbol{f}_0 \le 1,\end{aligned}$$

which leads to the ellipse constrain that $\|D^{-1/2}\boldsymbol{z}\|_2 \le 1$, where $\boldsymbol{z} = \frac{1}{\sqrt{n}}U^\top \boldsymbol{f}_0$. Obviously, $\|D_1^{-1/2}z_1\|_2 \le 1$, $\|D_2^{-1/2}z_2\|_2 \le 1$. Notice that $\frac{1}{\lambda}\boldsymbol{f}^\top U_2U_2^\top \boldsymbol{f} \le \boldsymbol{f}^\top U_2D_2^{-1}U_2^\top \boldsymbol{f} < n$, then $\boldsymbol{f}^\top U_2U_2^\top \boldsymbol{f} \le n\lambda$, and

$$\begin{aligned}T_2 &\le \|z_2\|_2 + \|\sqrt{D_2}\|_{\text{op}}\|\sqrt{D_2}\widetilde{S}_2^\top\|_{\text{op}}\|\widetilde{S}_1\|_{\text{op}}\|(\widetilde{S}_1^\top \widetilde{S}_1)^{-1}\|_{\text{op}}\|D_1^{-1/2}\|_{\text{op}}\|D_1^{-1/2}z_1\|_{\text{op}} \\ &\le (1+3c)\sqrt{\lambda}\end{aligned}$$

with probability at least $1 - e^{-c_1 s} - e^{-c_2 s_\lambda}$ by Assumption A2. Therefore, we have $\|\boldsymbol{z} - D\widetilde{S}^\top \widetilde{\beta}\|_2^2 \le c'\lambda$, where $c' = (1+3c)^2$. For the penalty term,

$$\begin{aligned}&n\widetilde{\beta}^\top S\boldsymbol{K}S^\top \widetilde{\beta} \\ \le& z_1^\top D_1^{-1}z_1 + \|z_1^\top D_1^{-1/2}\|_2\|D_1^{-1/2}\|_{\text{op}}\|\widetilde{S}_2D_2^{1/2}\|_{\text{op}}\|D_2^{1/2}\widetilde{S}^\top\|_{\text{op}}\|D_1^{-1/2}\|_{\text{op}}\|D_1^{-1/2}z_1\|_{\text{op}} \\ \le& 1 + c^2,\end{aligned}$$

where $c$ is constant from the definition 4.1. Finally, we can claim that

$$\|\operatorname{E}_\epsilon \widehat{f}_R - f_0\|_n^2 \le \|\boldsymbol{z} - \sqrt{n}D\widetilde{S}^\top \widetilde{\beta}\|_2^2 + n\lambda\widetilde{\beta}^\top \widetilde{S}^\top D\widetilde{S}\widetilde{\beta} \le C\lambda$$

with probability at least $1 - e^{-c_1 s} - e^{-c_2 s_\lambda}$, where $c$ is some constant, $C = 2 + 10c^2 + 6c$ is an absolute constant. □

## S.8 Proof of Corollary 4.4

*Proof.* Denote $\mathrm{E}_\epsilon \widehat{f}_R$ as the the expectation of $\widehat{f}_R$ w.r.t $\epsilon$. Note that

$$\|\widehat{f}_R - f_0\|_n^2 \le 2\|\widehat{f}_R - \mathrm{E}_\epsilon \widehat{f}_R\|_n^2 + 2\|\mathrm{E}_\epsilon \widehat{f}_R - f_0\|_n^2$$

and $\|\widehat{f}_R - \mathrm{E}_\epsilon \widehat{f}_R\|_n^2 = \frac{\epsilon^\top}{\sqrt{n}}\Delta^2\frac{\epsilon}{\sqrt{n}}$, where $\|\frac{\epsilon}{\sqrt{n}}\|_{\psi_2} \le \frac{L}{\sqrt{n}}$ and $\|\Delta^2\|_{\mathrm{op}} \le 1$. Recall $\|\cdot\|_{\psi_2}$ is the sub-Gaussian norm. Here $\|\epsilon\|_{\psi_2} \le L$, with $L$ as an absolute constant. Then by Hanson-Wright concentration inequality ([30]) (stated in Lemma S.4), with probability greater than $1 - e^{-c_1 s} - e^{-c_2 s_\lambda}$,

$$\begin{aligned}
&\mathrm{P}\Big(\|\widehat{f}_R - \mathrm{E}_\epsilon \widehat{f}_R\|_n^2 - \mathrm{E}_\epsilon \|\widehat{f}_R - \mathrm{E}_\epsilon \widehat{f}_R\|_n^2 \ge \frac{\mathrm{tr}(\Delta^2)}{2n} \big| \boldsymbol{x}, S\Big)\\
=&\mathrm{P}\Big(\frac{1}{n}\epsilon^\top \Delta^2 \epsilon - \frac{\mathrm{tr}(\Delta^2)}{n} \ge \frac{\mathrm{tr}(\Delta^2)}{2n}\big|\boldsymbol{x}, S\Big)\\
\le&\exp\Big(-c\min\big(\frac{\mathrm{tr}^2(\Delta^2)}{4K^4\|\Delta^2\|_{\mathrm{F}}^2}, \frac{\mathrm{tr}(\Delta^2)}{\|\Delta^2\|_{\mathrm{op}}}\big)\Big)\\
\le&\exp(-c\,\mathrm{tr}(\Delta^2)),
\end{aligned}$$

where $\|\cdot\|_{\mathrm{F}}$ is the Frobenius norm. The last inequality holds by the fact that $\|\Delta^2\|_{\mathrm{F}}^2 \le \|\Delta^2\|_{\mathrm{op}}\,\mathrm{tr}(\Delta^2)$ and $\|\Delta^2\|_{\mathrm{op}} \le 1$. Lastly, by (7.2), $\mathrm{tr}(\Delta^2) \ge \min\{s, \widehat{s}_\lambda\} \ge \widehat{s}_\lambda$, which goes to $+\infty$ as $n \to \infty$, we have that, with probability approaching 1, $\|\widehat{f}_R - \mathrm{E}_\epsilon \widehat{f}_R\|_n^2 \le \frac{3}{2}\mu_{n,\lambda}$. □

## S.9 Proof of Theorem 5.1

*Proof.* Note that $\lambda_m \asymp n^{-4m/(4m+1)}(\log\log n)^{2m/(4m+1)}$. Define $s_{\lambda_m} \asymp n^{2/(4m+1)}(\log\log n)^{-1/(4m+1)}$ for $m \ge 2$. If we define the event $\mathcal{C}(s_{\lambda_m}) = \{\boldsymbol{x}, S_m : \mathrm{tr}(\Delta_m^8) \asymp s_{\lambda_m}\}$, where

$$\Delta_m = \boldsymbol{K}_m S_m^\top (S_m \boldsymbol{K}_m^2 S_m^\top + \lambda_m S_m \boldsymbol{K}_m S_m^\top)^{-1} S_m \boldsymbol{K}_m,$$

then by Lemma 4.7, we have

$$\mathrm{P}\big(\mathcal{C}(s_{\lambda_m})\big) \ge 1 - e^{-c_1 s_m} - e^{-c_2 s_{\lambda_m}}.$$

For the remainder of the proof, our arguments are conditional on the event $\mathcal{C}(s_{\lambda_m})$.

Under the null hypothesis, $\boldsymbol{y} = \epsilon := (\epsilon_1, \dots, \epsilon_n)^\top \sim N(0, I_n)$. Let $A_{n,m} = \frac{\Delta_m^2}{\sqrt{2\,\mathrm{tr}(\Delta_m^4)}}$, then

$$\tau_m = \epsilon^\top A_{n,m}\epsilon - \mathbb{E}_\epsilon(\epsilon^\top A_{n,m}\epsilon).$$

Define $Z_n = (Z_{n,1}, \cdots, Z_{n,m_n-1})^\top$ be an $(m_n - 1)-$ dimensional centered Gaussian vector with covariance matrix $I_{m_n}$, which is an identity matrix.

We first prove

$$\sup_{\zeta\in\mathbb{R}}\Big|\mathrm{P}\big(\max_{1\le m\le m_n-1}\tau_{m+1} \le \zeta | \mathcal{C}(s_{\lambda_m})\big) - P\big(\max_{1\le m\le m_n-1} Z_{n,m} \le \zeta\big)\Big| \to 0. \tag{S.16}$$

We only need to verify the conditions in Lemma S.3. By Lemma 2.2 in [28], we have

$$\mathrm{E}(\tau_m^4) - 3\,\mathrm{E}(\tau_m^2)^2 = 48\,\mathrm{tr}(A_{n,m}^4) = 12\,\mathrm{tr}(\Delta_m^8)/\big(\mathrm{tr}(\Delta_m^4)\big)^2.$$

Note that $\mathrm{tr}(\Delta_m^8) \asymp s_{\lambda_m}$ and $\mathrm{tr}(\Delta_m^4) \asymp s_{\lambda_m}$. Therefore

$$\max_{1\le m\le m_n-1} \big(\mathrm{E}(\tau_{m+1}^4) - 3\,\mathrm{E}(\tau_{m+1}^2)^2\big) \log^6(m_n-1) \le C m_n s_{\lambda_{m+1}}^{-1} \log^6(m_n-1) \to 0. \tag{S.17}$$

On the other hand,

$$\begin{aligned}
&\max_{1\le k,l\le m_n-1} |I_{m_n-1}(k,l) - \mathrm{E}(\tau_k\tau_l)| \log^2(m_n-1) \qquad \text{(S.18)}\\
\le &\sum_{k=1}^{m_n-1} |I_{m_n-1}(k,k) - \mathrm{E}(\tau_{k+1}^2)| \log^2(m_n-1)\\
&+ 2 \sum_{1\le k<l\le m_n-1} |I_{m_n-1}(k,l) - \mathrm{E}(\tau_{k+1}\tau_{l+1})| \log^2(m_n-1)\\
=&J_1 + J_2.
\end{aligned}$$

For $J_1$, notice that $\mathrm{E}(\tau_k^2) = 2\,\mathrm{tr}(A_{n,k}^2) = 1$. Then $J_1 = 0$.

Next, we consider $J_2$.

$$\begin{aligned}
\mathrm{E}(\tau_{k+1}\tau_{l+1}) = 2\,\mathrm{tr}(A_{n,k+1}A_{n,l+1}) &= \frac{\mathrm{tr}(\Delta_{k+1}^2\Delta_{l+1}^2)}{\sqrt{\mathrm{tr}(\Delta_{k+1}^4)}\sqrt{\mathrm{tr}(\Delta_{l+1}^4)}}\\
&\le C\|\Delta_{k+1}^2\|_{\mathrm{op}} \frac{s_{\lambda_{l+1}}}{\sqrt{s_{\lambda_{l+1}}}\sqrt{s_{\lambda_{k+1}}}}\\
&\le C n^{-\frac{4}{(4k+5)(4l+5)}} (\log\log n)^{\frac{1}{2(4l+5)}}.
\end{aligned}$$

Therefore, $J_2 \le C m_n^2 n^{-\frac{4}{(4k+5)(4l+5)}} (\log\log n)^{\frac{1}{2(4l+5)}} \log^2 m_n \to 0$. And in (S.18), $\max_{1\le k,l\le m_n-1} |I_{m_n-1}(k,l) - \mathrm{E}(\tau_{k+1}\tau_{l+1})| \log^2(m_n-1) \to 0$. (S.16) has been proved.

Therefore, as $n\to\infty$,

$$\begin{aligned}
&\Big|\mathrm{P}\Big(B_n\big(\max_{1\le m\le m_n-1} \tau_{m+1} - B_n\big) \le c_{\bar\alpha} | \mathcal{C}(s_{\lambda_m})\Big) - \mathrm{P}\Big(B_n\big(\max_{1\le m\le m_n-1} Z_{n,m} - B_n\big) \le c_{\bar\alpha}\Big)\Big|\\
&\to 0.
\end{aligned}$$

It follows by [18] that, as $n\to\infty$,

$$\Big|\mathrm{P}\Big(B_n\big(\max_{1\le m\le m_n-1} Z_{n,m} - B_n\big) \le c_{\bar\alpha}\Big) - (1-\bar\alpha)\Big| \to 0.$$

Therefore, we have with probability at least $\prod_{m=2}^{m_n}(1 - e^{-c_1 s_m} - e^{-c_2 s_{\lambda,m}})$,

$$\mathrm{P}\Big(B_n\big(\max_{2\le m\le m_n} \tau_m - B_n\big) \le c_{\bar\alpha} | \boldsymbol{x}, S_m\Big) \to 1-\bar\alpha.$$

We claim that $\mathrm{P}\big(\tau_{n,m_n} \le c_{\bar\alpha}\big) \to 1-\bar\alpha$. Otherwise, there exists a subsequence $\{\boldsymbol{x}_{n_k}, S_{n'_k}\}$, s.t, for any $\epsilon > 0$,

$$\Big|\mathrm{E}_{\boldsymbol{x}_{n_k}, S_{n'_k}}\Big(\mathrm{P}\big(\tau_{n,m_n} \le c_{\bar\alpha} | \boldsymbol{x}_{n_k}, S_{n'_k}\big)\Big) - (1-\alpha)\Big| > \epsilon.$$

On the other hand, $\mathrm{P}\big(\tau_{n,m_n} \le c_{\bar\alpha} | \boldsymbol{x}_{n_k}, S_{n'_k}\big) \xrightarrow{p} 1-\alpha$, which implies there exists a subsubsequence $\boldsymbol{x}_{n_{n_k}}, S_{n_{n'_k}}$, such that

$$\mathrm{P}\big(\tau_{n,m_n} \le c_{\bar\alpha} | \boldsymbol{x}_{n_{n_k}}, S_{n_{n'_k}}\big) \xrightarrow{a.s.} 1-\alpha.$$

By Bounded Convergence Theorem,

$$\mathrm{E}_{\boldsymbol{x}_{n_{n_k}}, S_{n_{n'_k}}}\Big(\mathrm{P}\big(\tau_{n,m_n} \le c_{\bar\alpha} | \boldsymbol{x}_{n_{n_k}}, S_{n_{n'_k}}\big)\Big) \to (1-\alpha),$$

which is a contradiction. Proof is completed. □

## S.10 Proof of Theorem 5.2

*Proof.* We follow the same notation as in the proof of Theorem 5.1.

Since $m_n \to \infty$, eventually $m_n \ge m$. Then it holds that

$$\begin{aligned}
&\inf_{\substack{f\in\mathcal{B}_{n,m_*}\\ \|f\|_n \ge C_\varepsilon \delta(n,m_*)}} P_f(\tau_{n,m_n} \ge c_{\bar\alpha}) \\
=& \inf_{\substack{f\in\mathcal{B}_{n,m_*}\\ \|f\|_n \ge C_\varepsilon \delta(n,m_*)}} P_f(\max_{1\le k\le m_n} \tau_k \ge B_n + c_{\bar\alpha}/B_n) \\
\ge& \inf_{\substack{f\in\mathcal{B}_{n,m_*}\\ \|f\|_n \ge C_\varepsilon \delta(n,m_*)}} P_f(\tau_{m_*} \ge B_n + c_{\bar\alpha}/B_n) \\
\ge& \inf_{\substack{f\in\mathcal{B}_{n,m_*}\\ \|f\|_n \ge C_\varepsilon \delta(n,m_*)}} P_f\Big(\frac{\boldsymbol{y}^\top \Delta_{m_*}^2 \boldsymbol{y} - \mathrm{tr}(\Delta_{m_*}^2)}{s_{n,m_*}} \ge B_n + c_{\bar\alpha}/B_n\Big).
\end{aligned}$$

For $f \in \mathcal{B}_{n,m_*}$ with $\|f\|_n \ge C_\varepsilon \delta(n,m_*)$ for $C_\varepsilon$ to be described later, note that

$$\begin{aligned}
& \boldsymbol{y}^\top \Delta_{m_*}^2 \boldsymbol{y} \\
=& \boldsymbol{f}^\top \Delta_{m_*}^2 \boldsymbol{f} + 2\boldsymbol{f}^\top \Delta_{m_*}^2 \epsilon + \epsilon^\top \Delta_{m_*}^2 \epsilon \\
\equiv& J_1 + 2J_2 + J_3.
\end{aligned}$$

In the above, $J_1 = \boldsymbol{f}^\top \Delta_{m_*}^2 \boldsymbol{f} = n\|\mathrm{E}_\epsilon \widehat{f}_R\|_n^2 \ge (n/2)\|f\|_n^2 - n\|f - \mathrm{E}_\epsilon \widehat{f}_R\|_n^2 \ge (n/2)C_\varepsilon \delta^2(n,m_*) - nC\lambda_m = n(C_\varepsilon/2 - C)\delta^2(n,m_*)$, with $C_\varepsilon/2 - C > 0$, with probability greater than $1 - e^{-c_1 s_{m_*}} - e^{-c_2 s_{\lambda,m_*}}$ by Lemma 4.3.

For any $f$, we have

$$\mathrm{E}_\epsilon\{J_2^2\} = \boldsymbol{f}^\top \Delta_{m_*}^2 \boldsymbol{f} \le \boldsymbol{f}^\top \Delta_{m_*} \boldsymbol{f} = J_1,$$

implying that

$$\mathrm{P}(|J_2| > \varepsilon^{-1/2} J_1^{1/2} | \boldsymbol{x}, S) \leq \varepsilon,$$

i.e., $2J_2 = J_1^{1/2} O_P(1)$ uniformly for $f$. Also note that

$$\frac{J_3 - \mathrm{tr}(\Delta_{m_*}^2)}{s_{n,m_*}} \xrightarrow{d} N(0,1).$$

Combining the above analysis of $J_1, J_2, J_3$, we get that uniformly for any $f \in \mathcal{B}_{n,m_*}$ with $\|f\|_n \geq C_\varepsilon \delta(n, m_*)$, with probability approaching one,

$$\begin{aligned}
& \frac{J_1 + 2J_2 + J_3 - \mathrm{tr}(\Delta_{m_*}^2)}{s_{n,m_*}} \\
= & \frac{J_1 + J_1^{1/2} O_P(1)}{s_{n,m_*}} + O_P(1) \\
= & \frac{J_1}{s_{n,m_*}} (1 + o_P(1)) + O_P(1) \\
\geq & n(C_\varepsilon/2 - C)\delta(n, m_*)^2 s_{\lambda,m_*}^{-1/2} (1 + o_P(1)) + O_P(1) \\
\geq & (\log\log n)^{1/2} (1 + o_P(1)) + O_P(1) \\
\geq & B_n + c_{\bar{\alpha}}/B_n,
\end{aligned}$$

where the last inequality follows by $B_n \asymp (\log\log n)^{1/2}$. □

## S.11 Auxiliary Lemmas

**Lemma S.3.** *([24]) For each $n \in \mathbb{N}$, let $\boldsymbol{\xi}_n$ be an $N_n$-dimensional centered Gaussian vector with covariance matrix $\Sigma_n = (\Sigma_n(k,l))_{1 \leq k,l \leq N_n}$ and $d_n \geq 2$ be an integer. Also, for each $k = 1, \ldots, d_n$, let $A_{n,k}$ be an $N_n \times N_n$ symmetric matrix and $Z_n = (Z_{n,1}, \ldots, Z_{n,d_n})^\top$ be an $d_n$-dimensional centered Gaussian vector with covariance matrix $\mathfrak{C}_n = (\mathfrak{C}_n(k,l))_{1 \leq k,l \leq d_n}$. Set $F_{n,k} := \boldsymbol{\xi}_n^\top A_{n,k} \boldsymbol{\xi}_n - E[\boldsymbol{\xi}_n^\top A_{n,k} \boldsymbol{\xi}_n]$ and suppose that the following conditions are satisfied:*

1. *There is a constant $b > 0$ such that $\mathfrak{C}_n(k,k) \geq b$ for every $n$ and every $k = 1, \ldots, d_n$.*
2. $\max_{1 \leq k,l \leq d_n} \left( E(F_{n,k}^4) - 3E(F_{n,k}^2)^2 \right) \log^6 d_n \to 0$ *as* $n \to \infty$.
3. $\max_{1 \leq k,l \leq d_n} |\mathfrak{C}_n(k,l) - E(F_{n,k} F_{n,l})| \log^2 d_n \to 0$ *as* $n \to \infty$.

*Then we have*

$$\sup_{x \in \mathbb{R}} \left| P\Big( \max_{1 \leq k \leq d_n} F_{n,k} \leq x \Big) - P\Big( \max_{1 \leq k \leq d_n} Z_{n,k} \leq x \Big) \right| \to 0$$

*and*

$$\sup_{x\in\mathbb{R}}\Big|P\Big(\max_{1\le k\le d_n}|F_{n,k}|\le x\Big)-P\Big(\max_{1\le k\le d_n}|Z_{n,k}|\le x\Big)\Big|\to 0$$

*as* $n\to\infty$.

**Lemma S.4.** *(Hanson-Wright inequality [30]) Let* $X=(X_1,\cdots,X_n)\in\mathbb{R}^n$ *be a random vector with independent components* $X_i$ *which satisfy* $\mathrm{E}\,X_i=0$ *and* $\|X_i\|_{\psi_2}\le K$*. Let* $A$ *be an* $n\times n$ *matrix. Then, for every* $t\ge 0$,

$$P\Big(|X^\top AX-\mathrm{E}\,X^\top AX|>t\Big)\le 2\exp\Big(-c\min\big(\frac{t^2}{K^4\|A\|_{HS}^2},\frac{t}{K^2||A||}\big)\Big)$$

Here $\|A\|_{HS}$ is the Hilbert-Schmidt (or Frobenius) norm of $A$.

**Lemma S.5.** *(Eigenvalue interlacing theorem) Suppose* $A\in\mathbb{R}^{n\times n}$ *is symmetric. Let* $B\in\mathbb{R}^{m\times m}$ *with* $m<n$ *be a principal submatrix (obtained by deleting both* $i-$*th row and* $i-$*th column for some values of* $i$*). Suppose* $A$ *has eigenvalues* $\lambda_1\le\cdots\lambda_n$ *and* $B$ *has eigenvalues* $\beta_1\le\cdots\le\beta_m$*. Then*

$$\lambda_k\le\beta_k\le\lambda_{k+n-m}\qquad \textit{for } k=1,\cdots,m.$$

*And if* $m=n-1$,

$$\lambda_1\le\beta_1\le\lambda_2\le\beta_2\le\cdots\le\beta_{n-1}\le\lambda_n.$$

**Lemma S.6.** *(Weyl's inequality) Let* $M,H$ *and* $P$ *are* $n\times n$ *Hermitian matrices with* $M=H+P$*, where* $M$ *has eigenvalues* $\mu_1\ge\cdots\ge\mu_n$*, and* $H$ *has eigenvalues* $\nu_1\ge\cdots\ge\nu_n$*, and* $P$ *has eigenvalues* $\rho_1 geq\cdots\rho_n$*. Then the following inequalities hold for* $i=1,\cdots,n$*:*

$$\nu_i+\rho_n\le\mu_i\le\nu_i+\rho_1$$

*If* $P$ *is positive definite, then this implies*

$$\mu_i>\nu_i,\quad \forall i=1,\cdots,n.$$

## S.12 Some Simulation Results

In this section, we provide some simulation results to show the performance of the proposed testing procedure using Bernoulli random projection matrix. The simulation settings are as the same as in Simulation Study 6.1 and 6.2, except using Bernoulli random matrices.

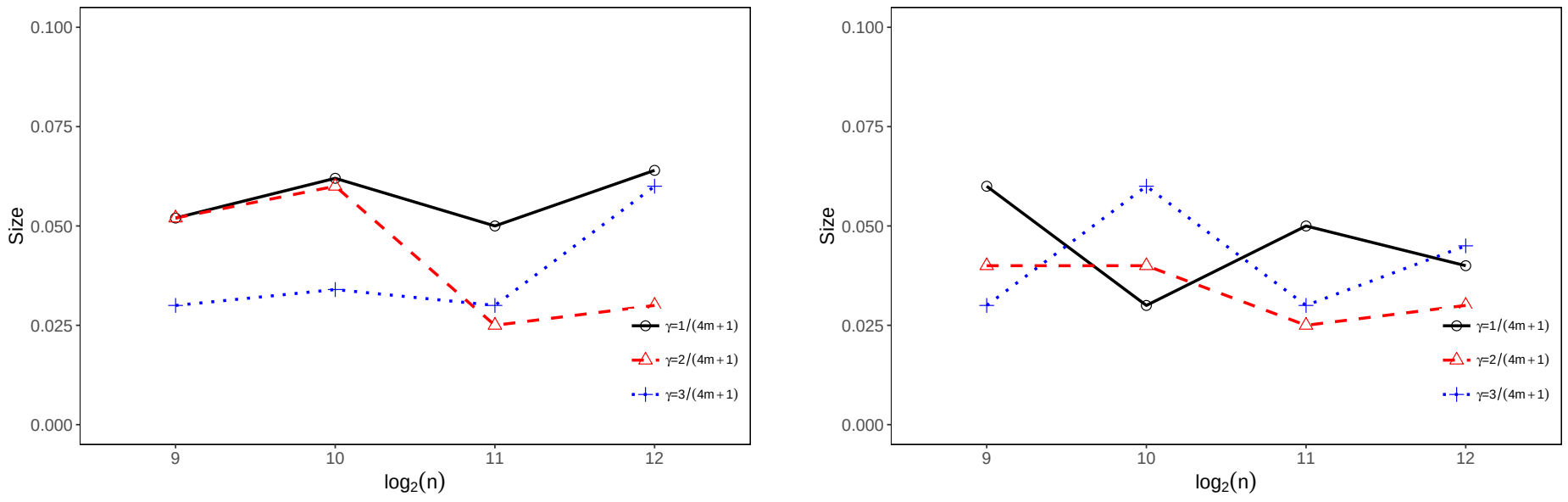

$(a)$DT $(b)$AT

Figure 8: Bernoulli random matrix. Size for $(a)$ DT and $(b)$ AT with projection dimension varies. Signal strength $c = 0$.

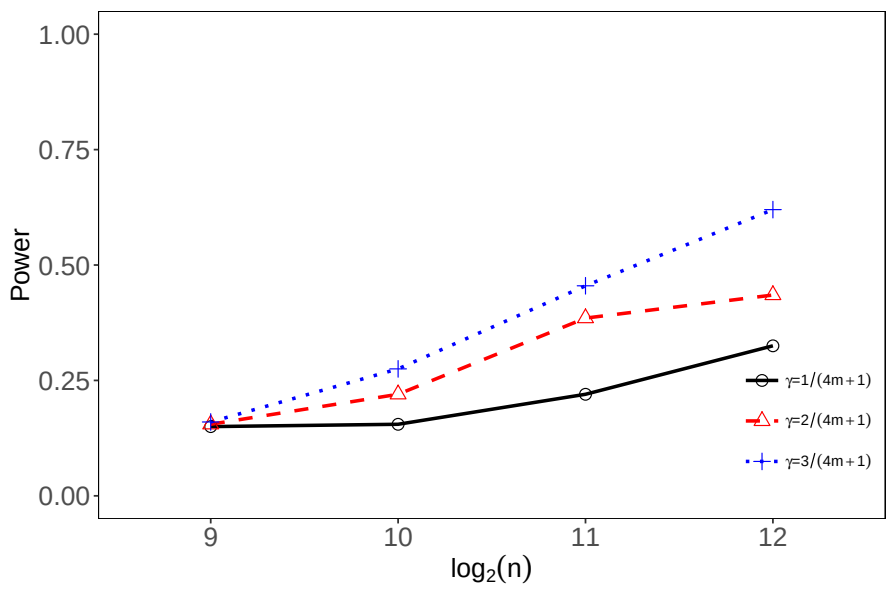

$(a)$ DT

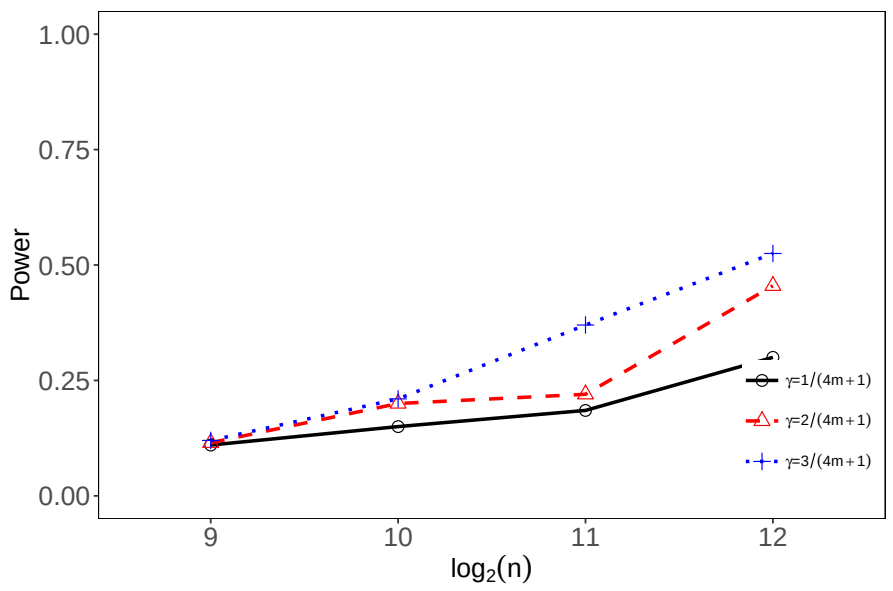

$(b)$AT

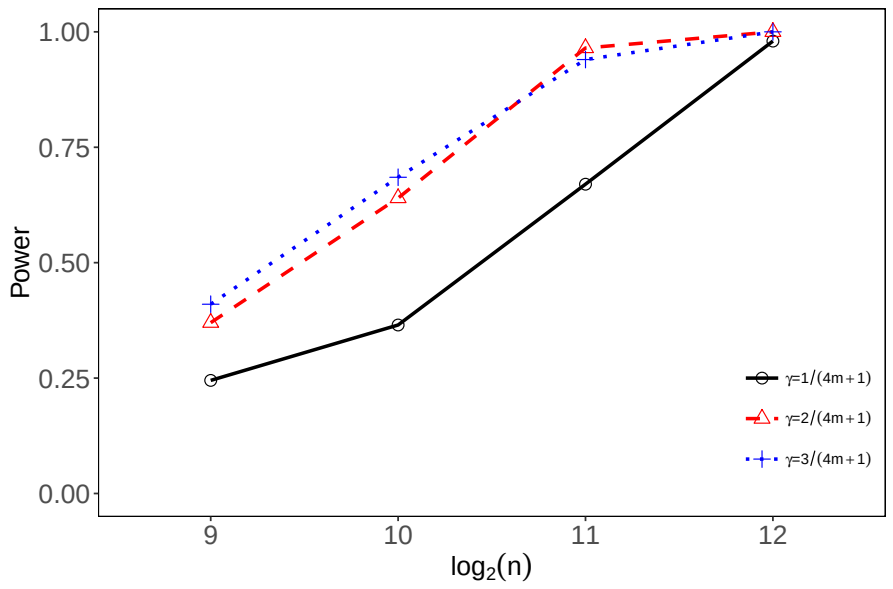

$(c)$DT

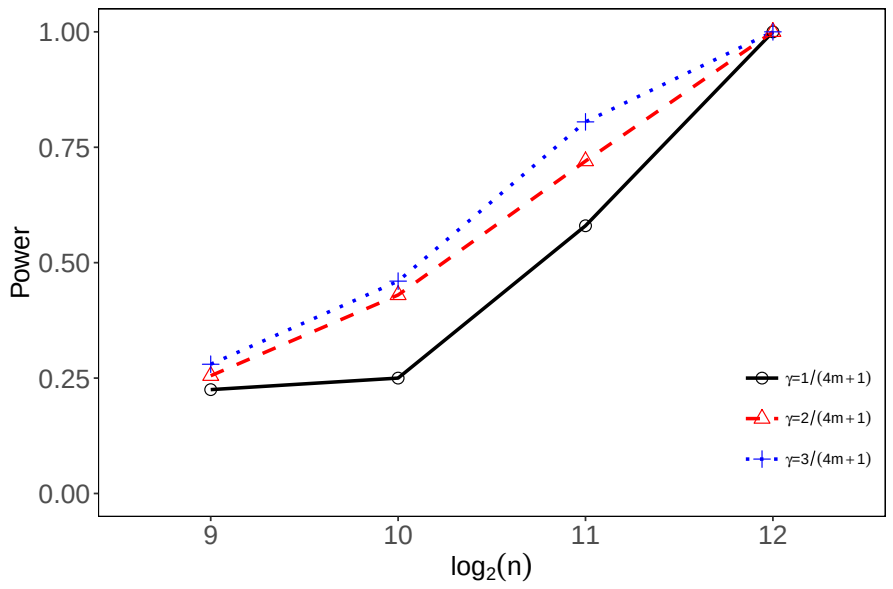

$(d)$AT

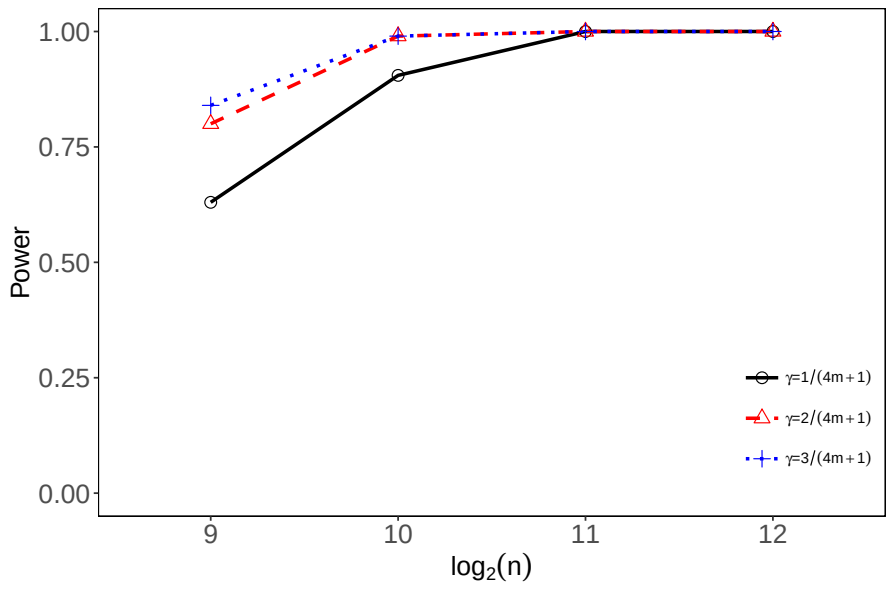

$(e)$DT

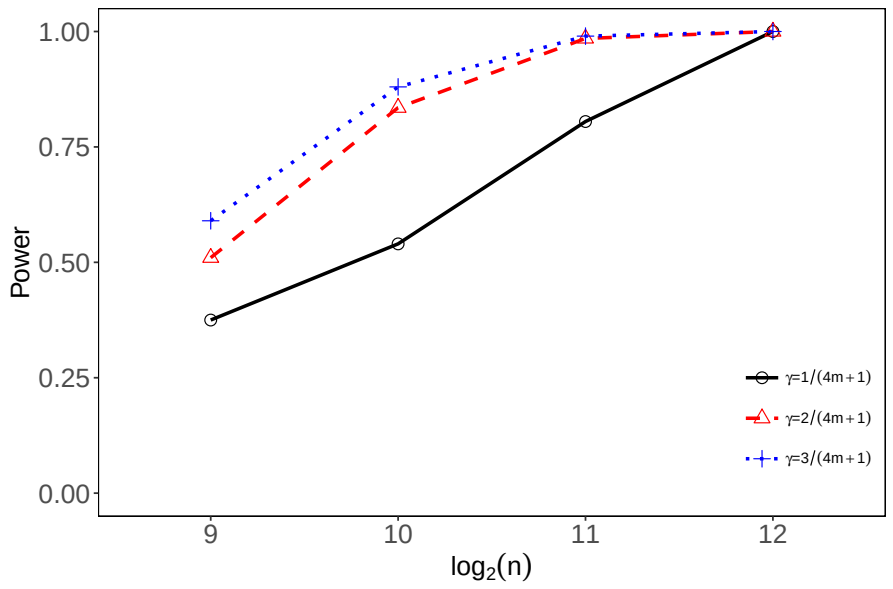

$(f)$AT

Figure 9: Bernoulli random matrix. Power for DT and AT with projection dimension varies. Signal strength $c = 0.01$ for $(a)$ and $(b)$; $c = 0.02$ for $(c)$ and $(d)$; $c = 0.03$ for $(e)$ and $(f)$.

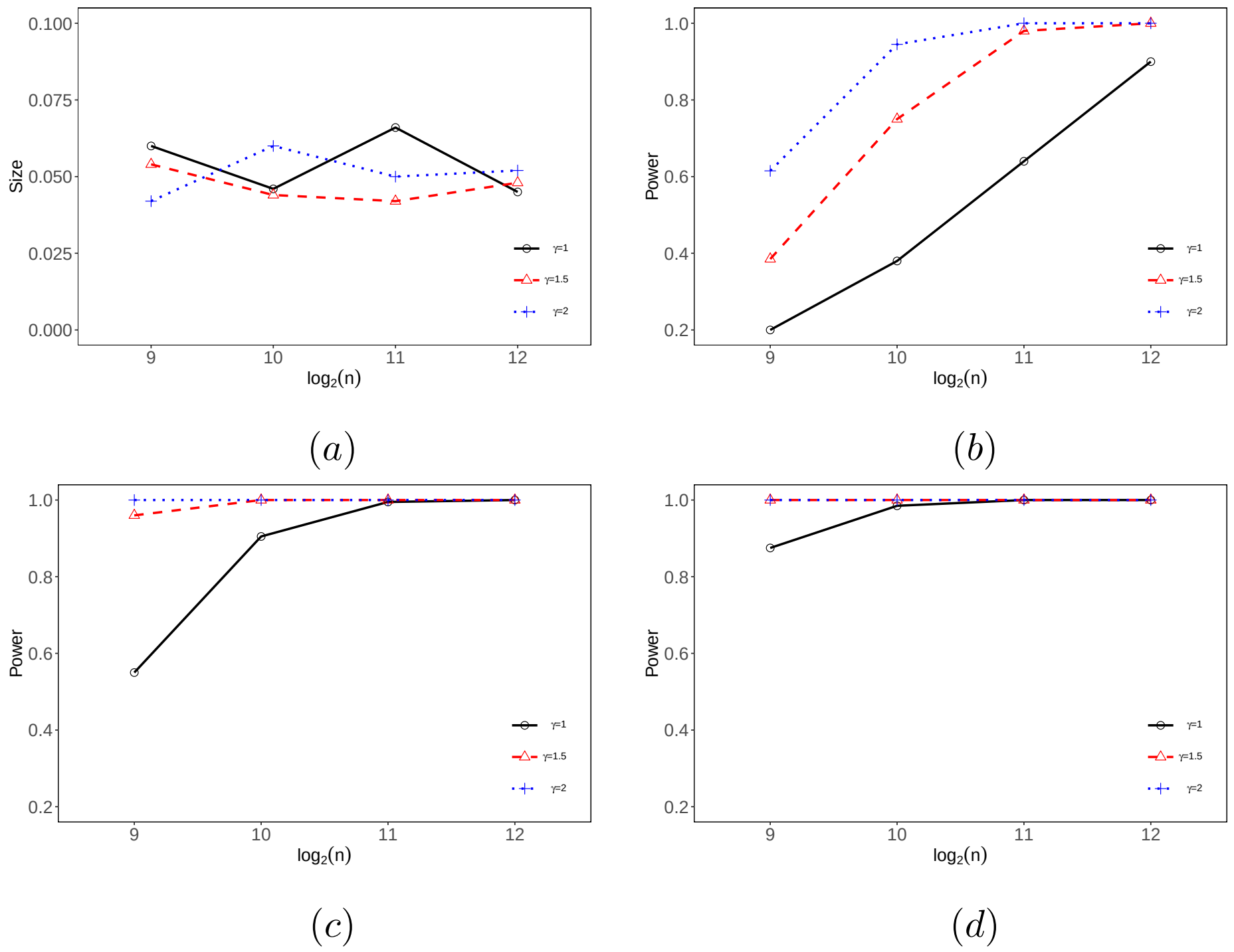

Figure 10: Bernoulli random matrix. Size and power for DT with projection dimension varies. Signal strength $c = 0$ for $(a)$; $c = 0.05$ for $(b)$; $c = 0.1$ for $(c)$; $c = 0.15$ for $(d)$.